\documentclass[10pt]{articleHJ}
\usepackage[english]{babel}
\usepackage{amsmath, amsthm, amsfonts, mathrsfs, amsfonts,bm,array}
\usepackage{mathtools}
\usepackage{bbm}
\usepackage{amssymb}
\usepackage{graphicx}
\DeclareMathAlphabet\bmcal{OMS}{cmsy}{b}{n}

\newcolumntype{L}{>{$}l<{$} }

\usepackage[colorlinks=true]{hyperref}
\hypersetup{urlcolor=blue, citecolor=red}

\usepackage[text={6.0in,8.6in},centering,letterpaper]{geometry}
\usepackage[numbers,comma,square,sort&compress]{natbib}

\usepackage{tabularx}
\usepackage{enumitem}

\newlength{\defbaselineskip}
\newcommand{\setlinespacing}[1]%
           {\setlength{\baselineskip}{#1 \defbaselineskip}}

\newcommand{\singlespacing}{\setlength{\baselineskip}{\defbaselineskip}}

\newcommand{\rmax}{\ensuremath{r_{\mathrm{max}}}}
\newcommand{\rmin}{\ensuremath{r_{\mathrm{min}}}}

\newcommand{\D}{\ensuremath{\mathcal{D}}}
\newcommand{\1}{\ensuremath{\mathbbm{1}}}

\newcommand{\N}{\ensuremath{\mathbb{N}}}
\newcommand{\Z}{\ensuremath{\mathbb{Z}}}

\newcommand{\R}{\ensuremath{\mathbb{R}}}
\newcommand{\C}{\ensuremath{\mathbb{C}}}

\renewcommand{\epsilon}{\ensuremath{\varepsilon}}

\renewcommand{\ker}{\ensuremath{\mathrm{Ker}}}
\newcommand{\Range}{\ensuremath{\mathrm{Range}}}

\newcommand{\codim}{\ensuremath{\mathrm{codim}}}
\renewcommand{\span}{\ensuremath{\mathrm{span}}}
\renewcommand{\min}{\ensuremath{\mathrm{min}}}
\newcommand{\ind}{\ensuremath{\mathrm{ind}}}

\DeclarePairedDelimiter\abs{\lvert}{\rvert}

\DeclarePairedDelimiter{\ip}\langle\rangle
\DeclarePairedDelimiter{\nrm}\lVert\rVert

\theoremstyle{plain}
\newtheorem{theorem}{Theorem}[section]
\newtheorem{proposition}[theorem]{Proposition}
\newtheorem{lemma}[theorem]{Lemma}
\newtheorem{corollary}[theorem]{Corollary}

\newtheorem{definition}[theorem]{Definition}

\theoremstyle{definition}

\newtheorem{remark}[theorem]{Remark}

\newtheorem*{assumption*}{\assumptionnumber}
\providecommand{\assumptionnumber}{}
\makeatletter
\newenvironment{assumption}[2]
 {%
  \renewcommand{\assumptionnumber}{Assumption ($#2$#1)}%
  \begin{assumption*}%
  \protected@edef\@currentlabel{#1}%
 }
 {%
  \end{assumption*}
 }
\makeatother
\newcommand{\asref}[2]{$#2$\ref{#1}}

\numberwithin{equation}{section}

\begin{document}
\begin{frontmatter}

\title{Exponential Dichotomies for Nonlocal Differential Operators with Infinite Range Interactions}
\journal{...}
\author[LDA]{W. M. Schouten-Straatman\corauthref{coraut}},
\corauth[coraut]{Corresponding author. }
\author[LDB]{H. J. Hupkes}
\address[LDA]{
  Mathematisch Instituut\textemdash Universiteit Leiden \\
  P.O. Box 9512; 2300 RA Leiden; The Netherlands \\ Email:  {\normalfont{\texttt{w.m.schouten@math.leidenuniv.nl}}}
}
\address[LDB]{
  Mathematisch Instituut\textemdash Universiteit Leiden \\
  P.O. Box 9512; 2300 RA Leiden; The Netherlands \\ Email:  {\normalfont{\texttt{hhupkes@math.leidenuniv.nl}}}
}

\date{\today}

\begin{abstract}
\singlespacing
We show that MFDEs with infinite range discrete and/or continuous interactions
admit exponential dichotomies, 
building on the Fredholm theory developed by Faye and Scheel for such systems.
For the half line, we refine the earlier approach by Hupkes and Verduyn Lunel.
For the full line,
we construct these splittings
by generalizing the finite-range results obtained by Mallet-Paret and Verduyn Lunel. 
The finite dimensional space that is `missed' by these splittings can be characterized using the Hale inner product, but the resulting degeneracy issues raise subtle questions that are much harder to resolve than in the finite-range case. Indeed, there is no direct analogue for the standard 'atomicity' condition that is typically used to rule out degeneracies, since it explicitly references the smallest and largest shifts.

We construct alternative criteria that exploit finer information on the structure of the MFDE. Our results are optimal when the coefficients are cyclic with respect to appropriate shift semigroups or when the standard positivity conditions 
typically associated to comparison principles are satisfied. We illustrate these results with explicit examples and counter-examples that involve the Nagumo equation.
\end{abstract}


\end{frontmatter}

\section{Introduction}
\label{sec:int}

Many physical, chemical and biological systems
feature nonlocal interactions that
can have a fundamental impact on the underlying dynamical
behaviour.
A typical mechanism to generate such nonlocality is to
include dependencies
on spatial averages of model components, often as part of a multi-scale approach.
For example, plants take up water from the surrounding soil through their spatially-extended root network, which can be modelled by nonlocal logistic growth terms \cite{gilad2004ecosystem,gilad2007mathematical}.
The propagation of cancer cells depends on the orientation of the surrounding  extracellular matrix fibres, which leads naturally to nonlocal  
flux terms \cite{Shuttleworth2019}.
Additional examples can be found in the fields of population dynamics \cite{shen2010spreading,grinfeld2005non,shen2016spectral,britton1990spatial,sun2017periodic},
material science
\cite{feireisl2004non,bates1997global,bates1997traveling,wang2002metastability} and many others.\\

A second fundamental route that leads to nonlocality is the
consideration of spatial domains that
feature some type of  \textit{discreteness}.
The
 broken translational and rotational symmetries often lead to highly complex and surprising behaviour that disappears in the continuum limit.
 For example, recent experiments have established that light waves can be trapped in well-designed photonic lattices \cite{PhysRevLett.114.245503,PhysRevLett.114.245504}. 
 Other settings where
 discrete topological effects 
 play an essential role 
 include
 the movement of domain walls \cite{Dmitriev2000domain}, the propagation of dislocations through crystals \cite{Celli1970motion} and the development of fractures in elastic bodies \cite{Slepyan2012models}. In fact, even the simplest discretizations of standard scalar reaction-diffusion systems are known to have far richer properties than their continuous local counterparts \cite{chow1996dynamics,chow1995pattern,hupkestraveling}.\\

\paragraph{Myelinated nerve fibres}
A commonly used modelling prototype to illustrate these issues concerns the propagation of electrical signals through nerve fibres. These nerve fibres are insulated by segments of myelin coating
that are separated by periodic gaps at the so-called nodes of Ranvier \cite{RANVIER1878}. Signals travel quickly through the coated regions, but lose strength rapidly.
The movement
through the gaps is much slower, but the signal
is chemically reinforced
in preparation for the next segment \cite{LILLIE1925}.\\

One of the first mathematical models proposed to capture this propagation
was the FitzHugh-Nagumo partial differential equation (PDE) \cite{FITZHUGH1966}.
This model is able to reproduce the travelling pulses observed in nature 
\cite{FITZHUGH1968} and has been studied extensively as a consequence. These studies have led to the development of many important mathematical techniques
in areas such as singular perturbation theory \cite{CARP1977,HAST1976,JONES1984,JONESKOPLAN1991,CARTER2016,CARTER2015} 
variational calculus \cite{chen2015traveling},
Maslov index theory
\cite{beck2018instability,chen2014stability,howard2016maslov,cornwell2017opening,cornwell2018existence}
and stochastic dynamics \cite{Hamster2017stability,hamster2018stability,hamster2019travelling}.
However, as a fully local equation it is unable to incorporate the discrete structure
in a direct fashion.\\

In order to repair this, Keener and Sneed \cite{EVVPD18} proposed
to replace the FitzHugh-Nagumo PDE by its discretized counterpart \begin{equation}\label{eq:fhnLDE}\begin{array}{lcl}\dot{u}_j&=&u_{j+1}+u_{j-1}-2u_j+g(u_j;a)-w_j,
\\[0.2cm]
\dot{w}_j&=&\rho[u_j-  w_j],\end{array}\end{equation}
indexed on the spatial lattice $j \in \mathbb{Z}$.
Here the variable $u_j$ describes the potential on the $j^{\text{th}}$ node of Ranvier, while $w_j$ describes a recovery component. The nonlinearity
can be taken as the bistable
cubic
$g(u;a)=u(1-u)(u-a)$
for some $a \in (0, 1)$ and $0 < \rho \ll 1$ is a small parameter.
Such an infinite system of coupled ODEs is referred to as a lattice differential equation (LDE)\textemdash a class of equations that arises naturally when discretizing the spatial derivatives in PDEs.\\

 Since we are mainly interested in
 the propagation of electrical pulses, we introduce the travelling wave Ansatz
\begin{equation}\label{eq:trvwvansatz}
    \begin{array}{lcl}
         (u_j,w_j)(t)&=&(\overline{u},\overline{w})(j+ct),
         \qquad 
         \qquad
       (\overline{u},\overline{w})(\pm \infty) = 0.
    \end{array}
\end{equation}
Here $c$ is the speed of the wave and the smooth functions
$(\overline{u},\overline{w}): \R \to \R^2$ represent the two waveprofiles.
Plugging (\ref{eq:trvwvansatz}) into the LDE (\ref{eq:fhnLDE}) yields the differential equation 
\begin{equation}\label{eq:fhnMFDE}\begin{array}{lcl}c\overline{u}'( \sigma)&=&\overline{u}( \sigma +1 )+\overline{u}( \sigma -1)-2\overline{u}( \sigma)+g(\overline{u}( \sigma); a)-\overline{w}( \sigma) ,
\\[0.2cm]
c\overline{w}'( \sigma)&=&\rho[\overline{u}( \sigma)-  \overline{w}( \sigma)]\end{array}\end{equation}
in which $ \sigma=j+ct$. Since this system 
contains both advanced (positive) and retarded (negative) shifts, such an equation is called a functional differential equation of mixed type (MFDE).\\

In \cite{HJHSTBFHN,HJHFZHNGM} Hupkes and Sandstede
established the existence and nonlinear stability of such pulses, under a `non-pinning' condition
for the associated Nagumo LDE
\begin{equation}\label{eq:nagumoldefin}
    \begin{array}{lcl}
         \dot{u}_j&=&
          u_{j+1}+u_{j-1}-2u_j  +g(u_j; a).
    \end{array}
\end{equation}
This LDE arises when considering the first component of 
(\ref{eq:fhnLDE}) with $w = 0$. It admits travelling front solutions
\begin{equation}
\label{eq:int:nag:front}
 u_j(t)=\overline{u}_*(j+c_* t)   ,
 \qquad \qquad
 \overline{u}_*(-\infty) = 0, 
 \qquad
 \overline{u}_*(+\infty) = 1
\end{equation}
that necessarily satisfy the MFDE
\begin{equation}\label{eq:nagumomfdeintrofin}
    \begin{array}{lcl}
     c_* \overline{u}_*'( \sigma)&=&\overline{u}_*( \sigma+1)+\overline{u}_*( \sigma-1)-2\overline{u}_*( \sigma)+g(\overline{u}_*( \sigma); a).
    \end{array}
\end{equation}
The `non-pinning' condition  mentioned above demands that the wavespeed $c_*$\textemdash which depends uniquely on $ a$ \cite{MPB}\textemdash does not vanish. In the PDE case this is automatic for $ a \neq \frac{1}{2}$, but in the discrete setting
this is a non-trivial demand due to the energy barriers caused by the lattice \cite{VL28,HOFFMPcrys,BEYN1990,mallet2001crystallographic,32,EVV2005Nonlin}.\\

The main idea behind the approach developed in  \cite{HJHSTBFHN,HJHFZHNGM} is to use Lin's method \cite{LIN1993,HJHLIN} to combine the fronts (\ref{eq:int:nag:front})
and their reflections to form so-called quasi-front and quasi-back solutions
to \eqref{eq:fhnMFDE}.
Such solutions admit gaps in predetermined finite-dimensional subspaces
that can be closed by choosing the correct wavespeed. The existence of these subspaces is directly related to the construction of exponential dichotomies
for the linear MFDE
\begin{equation}\label{eq:linearnagumomfde}
    \begin{array}{lcl}
         cu'( \sigma)&=&u( \sigma+1)+u( \sigma-1)-2u( \sigma)+g_u(\overline{u}_*( \sigma); a)u( \sigma),
    \end{array}
\end{equation}
which arises as the linearization of
(\ref{eq:nagumomfdeintrofin}) around the front
solutions \eqref{eq:int:nag:front}.

\paragraph{Exponential dichotomies for ODEs}
Roughly speaking, a linear differential equation is said to admit an
\textit{exponential dichotomy} if the space of initial conditions can be written as a direct sum of a stable and an unstable subspace. Initial conditions in the former can be continued as solutions that decay exponentially in forward time, while initial conditions in the latter admit this property in backward time. In order to be more specific, we first
restrict our attention to the ODE 
\begin{equation}\label{eq:generalode}
    \begin{array}{lcl}
         \frac{d}{d\sigma}u&=&A(\sigma)u,
    \end{array}
\end{equation}
referring to  
the review paper by Sandstede \cite{sandstede2002stability} for further details. Here  $u(\sigma)\in\C^M$ and $A(\sigma)$ is an $M\times M$ matrix for any $\sigma\in \R$.
Let us write
$\Phi(\sigma,\tau)$ for the evolution operator associated to (\ref{eq:generalode}), which maps $u(\tau)$ to $u(\sigma)$.\\

Suppose first that the system (\ref{eq:generalode}) is autonomous and hyperbolic, i.e. $A(\sigma) =A$ for some matrix $A$ that has no spectrum on the imaginary axis.
Writing
$E^s_0$ and $E^u_0$ for the generalized stable respectively unstable eigenspaces of $A$, we subsequently obtain the decomposition
\begin{equation}
\label{eq:int:split:cm:aut}
    \begin{array}{lcl}
     \C^M&=&E^s_0\oplus E^u_0.
    \end{array}
\end{equation}
In addition, each of these subspaces is invariant under the action of $\Phi(\sigma, \tau) = \mathrm{exp}[A(\sigma - \tau)]$,
which decays exponentially  on $E^s_0$ for $\sigma > \tau$ and on $E^u_0$ for $\sigma < \tau$.\\

 In order to
 generalize such decompositions to non-autonomous settings, the splitting
 \eqref{eq:int:split:cm:aut} will need to vary with the base time $\tau \in I$. Here we pick $I$ to be one of the three intervals $\R^-$, $\R^+$ or  $\R$. In particular, 
 (\ref{eq:generalode}) is said to be \textit{exponentially dichotomous} on $I$ if the following properties hold.
\begin{itemize}
    \item There exists a family of projection operators  $\{P(\tau)\}_{\tau \in I}$ on $\C^M$ that commute with the evolution $\Phi(\sigma,\tau)$.
    \item The restricted operators $\Phi^s(\sigma,\tau):=\Phi(\sigma,\tau)P(\tau)$ and $\Phi^u(\sigma,\tau):=\Phi(\sigma,\tau)\big(\mathrm{id}-P(\tau)\big)$ decay exponentially for $\sigma\geq\tau$ respectively $\sigma\leq \tau$.
\end{itemize}
 
Many important features concerning these dichotomies were first described by Palmer in \cite{palmer1984exponential,palmer1988exponential}. For example, the well-known roughness theorem states that exponential dichotomies persist under small perturbations of the matrices $A(\sigma)$. In addition, there is a close connection with the Fredholm properties of the associated linear operators. Consider for example the family of linear operators 
 \begin{equation}\label{eq:def:mathcalT}
     \begin{array}{lclcl}
          \Lambda(\lambda):H^1(\R;\C^M)&\rightarrow &L^2(\R;\C^M),\qquad u&\mapsto &\frac{d}{d\sigma}u-A(\sigma)u-\lambda u,
     \end{array}
 \end{equation}
 defined for $\lambda\in\C$.
 Then $\Lambda(\lambda)$ is a Fredholm operator if and only if the system
 \begin{equation}\label{eq:generalodelambda}
    \begin{array}{lcl}
     \frac{d}{d\sigma}u&=&A(\sigma)u+\lambda u
    \end{array}
\end{equation}
 admits exponential dichotomies on both $\R^+$ and $\R^-$.
 In addition, 
 $\Lambda(\lambda)$ is invertible if and only if \eqref{eq:generalodelambda} admits exponential dichotomies on $\R$. Since systems of the form \eqref{eq:generalodelambda} arise frequently when considering the spectral properties of wave solutions to nonlinear PDEs, exponential dichotomies have a key role to play in this area. In fact, the well-known Evans function \cite{palmer1984exponential,palmer1988exponential,pego1992eigenvalues,Evans1972} detects precisely when the dichotomies on $\R^-$ and $\R^+$ can be patched together to form a dichotomy on $\R$.

\paragraph{Exponential dichotomies for MFDEs}

Several important points need to be addressed before the concepts above
can be extended to linear MFDEs such as \eqref{eq:linearnagumomfde}. The first issue is that MFDEs are typically ill-posed \cite{RUSHA1989}, preventing a natural analogue of the evolution operator $\Phi$ to be defined. The second issue is that $\C^M$ is no longer an appropriate state space. For example, computing $u'(0)$ in \eqref{eq:linearnagumomfde} requires knowledge of $u$ on the interval $[-1, 1]$. These issues were resolved independently and simultaneously by Mallet-Paret and Verduyn Lunel
in \cite{MPVL} and by H\"arterich, Scheel and Sandstede in \cite{harterich}
by decomposing suitable function spaces into separate parts that individually do admit (exponentially decaying) semiflows.\\

Applying the results in 
\cite{MPVL} to \eqref{eq:linearnagumomfde}, we obtain the decomposition
\begin{equation}
\label{eq:int:decomp:mfde}
\begin{array}{lcl}
    C([-1,1];\R)
    & = &P(\tau) + Q(\tau) + \Gamma(\tau)
    \end{array}
\end{equation}
for each $\tau \in \R$.
Here $\Gamma(\tau)$ is finite dimensional,
while functions in $P(\tau)$ and $Q(\tau)$ can be extended to exponentially decaying solutions of the MFDE \eqref{eq:linearnagumomfde} on the intervals $(-\infty,\tau]$
respectively $[\tau, \infty)$. In particular,
the intersection $P(\tau) \cap Q(\tau)$ contains segments of functions that belong to the kernel of the associated linear operator
\begin{equation}
\label{eq:int:def:L:mfde}
\begin{array}{lcl}
    [\mathcal{L} v]( \sigma)
   & = &- cv'( \sigma) + v( \sigma + 1) + v ( \sigma - 1)  - 2 v( \sigma) + g_u(\overline{u}_*( \sigma); a) v( \sigma).
   \end{array}
\end{equation}
After dividing these segments out from either $P$ or $Q$, the decomposition \eqref{eq:int:decomp:mfde} becomes a direct sum.
Similar results were obtained in \cite{harterich}, but here the authors use 
the augmented statespace
$\mathbb{C}^M \times L^2([-1,1];\R)$.\\

In many applications, it is crucial to understand the dimension of $\Gamma(\tau)$. A key tool to achieve this is the so-called Hale inner product \cite{HLVL}, which in the present context is given by
\begin{equation}
\begin{array}{lcl}
    \langle  \psi , \phi \rangle_\tau
    & = &\frac{1}{c}\big[ \psi(0)\phi(0)+\int\limits_{-1}^0\psi(s+1)\phi(s)ds-\int\limits_0^1 \psi(s-1)\phi(s)ds\big]
    \end{array}
\end{equation}
for two functions $\phi, \psi \in C([-1,1];\R)$. Indeed, one of the main results achieved in \cite{MPVL} is the identification
\begin{equation}\label{eq:PplusQHale}
    \begin{array}{lcl}
        P(\tau) + Q(\tau) 
        &=& \big\{ \phi \in C([-1,1];\R) : \langle b(\tau + \cdot) , \phi \rangle_{\tau} = 0 
        \hbox{ for every } b \in \mathrm{ker}\, \mathcal{L}^* \big\}.
    \end{array}
\end{equation}
Here $\mathcal{L}^*$ stands for the formal adjoint of $\mathcal{L}$,
which arises by switching the sign of $c$ in \eqref{eq:int:def:L:mfde}.\\

There are two potential issues that can impact the usefulness of this result. The first is that the Hale inner product could be degenerate, the second is that kernel elements of $\mathcal{L}^*$ could vanish on large intervals. For instance, \cite[Ex. V.4.8]{VLDELAY} features
an example system that admits compactly supported kernel elements, which are often referred to as \textit{small solutions}. Fortunately, both types of degeneracies can be ruled out by imposing an invertibility condition on the coefficients related to the smallest and largest shifts in the MFDE. This is easy to check and obviously satisfied for \eqref{eq:linearnagumomfde}. \\

These results from \cite{MPVL,harterich}
have been used in a variety of settings by now. These include the construction of travelling waves \cite{HJHFZHNGM,IOOSS2005}, the stability analysis of such waves \cite{HJHNLS,HJHSTBFHN},
the study of homoclinic bifurcations 
 \cite{GEORGI2008Reversible,HJHLIN},
the analysis of pseudospectral approximations
\cite{BREMASVER2006}
and the detection of indeterminacy in economic models \cite{HJHECON}.
Partial extensions of these results for MFDEs taking values in Banach spaces can be found in \cite{HJHPOLLUTION},
but only for autonomous
systems at present.

\paragraph{Infinite-range interactions}

In recent years, an active interest has arisen in systems that feature interactions that can take place over arbitrarily large distances. For example,
diffusion models based on L\'evy processes lead naturally to fractional Laplacians in the underlying PDE \cite{applebaum2009levy,bertoin1996levy}. These operators are inherently nonlocal 
and often feature infinitely many terms in their discretization schemes \cite{Ciaurri}. Systems of this type have been used for example to describe amorphous semiconductors \cite{GU1996}, liquid crystals \cite{Ciuchi}, porous media \cite{bonforte2014existence}
and game theory \cite{bjorland2012nonlocal};
see \cite{bucur2016nonlocal} for an accessible introduction. Examples  featuring 
other types  of infinite-range interactions
include Ising models to describe the behaviour of
magnetic spins on a grid \cite{BatesInfRange}
and SIR models to capture the spread of infectious diseases \cite{li2014traveling}.\\

Returning to the study of nerve axons, 
let us now consider large networks of neurons. These neurons interact with each other over large distances through their connecting fibres \cite{BRESS2014,BRESS2011,PINTO2001,SNEYD2005}. Such systems generally have a very complex structure and finding effective equations to describe their behaviour is highly challenging. One candidate that has been proposed \cite{BRESS2014} involves  FitzHugh-Nagumo type models
such as
\begin{equation}\begin{array}{lcl}\label{eq:fhninfrange}\dot{u}_j&=& h^{-2}\sum\limits_{k\in\Z_{>0}}e^{-k^2}[u_{j+k}+u_{j-k}-2u_j]+g(u_j; a)-w_j,\\[0.2cm]
\dot{w}_j&=&\rho[u_j-  w_j].
\end{array}\end{equation}
Here the constant $h>0$ represents the (scaled) discretization distance. Alternatively,
one can replace or supplement the sum in \eqref{eq:fhninfrange} by including a convolution with a smooth kernel.\\

The travelling wave Ansatz 
\begin{equation}\label{eq:trvwvansatzinf}
    \begin{array}{lcl}
         (u_j,w_j)(t)&=&(\overline{u}_h,\overline{w}_h)(hj+c_ht),
         \qquad 
         \qquad
       (\overline{u}_h,\overline{w}_h)(\pm \infty) = 0
    \end{array}
\end{equation}
now yields the MFDE 
\begin{equation}\label{eq:fhnMFDEinfrange}\begin{array}{lcl}c_h\overline{u}_h'( \sigma)&=&h^{-2}\sum\limits_{k\in\Z_{>0}}e^{-k^2}[\overline{u}_h( \sigma +hk )+\overline{u}_h( \sigma -hk)-2\overline{u}_h( \sigma)]+g(\overline{u}_h( \sigma); a)-\overline{w}_h( \sigma)\\[0.2cm]
c_h\overline{w}_h'( \sigma)&=&\rho[\overline{u}_h( \sigma)-  \overline{w}_h( \sigma)],
\end{array}\end{equation}
which includes infinite-range interactions. In particular, it is no longer possible to apply the exponential splitting results from \cite{MPVL,harterich}. Nevertheless,
Faye and Scheel obtained an existence result for such waves in \cite{Faye2015}, pioneering a new approach to analyze spatial dynamics that circumvents the use of a state space. Extending the spectral convergence technique developed by Bates, Chen and Chmaj \cite{BatesInfRange}, we were able to show that such waves are nonlinearly stable \cite{HJHFHNINFRANGE}, but only for small $h > 0$. In any case, at present there is no clear mechanism that allows finite-range results to be easily extended to settings with infinite-range interactions. \\

\paragraph{Infinite-range MFDEs}

In this paper we take a step towards building such a bridge by constructing exponential dichotomies
for the non-autonomous, integro-differential MFDE
\begin{equation}\label{eq:basissysteem}
    \begin{array}{lcl}
         \dot{x}(\sigma)&=&\sum\limits_{j=-\infty}^\infty A_j(\sigma)x(\sigma+r_j) +\int_{\R}\mathcal{K}(\xi;\sigma)x(\sigma+\xi)d\xi,
    \end{array}
\end{equation}
which is allowed to have infinite-range interactions.
Here, we have $x(\sigma)\in \C^M$ for $t\in\R$ and the scalars $r_j$ for $j\in\Z$ are called the \textit{shifts}. Typically, we use $C_b(\R)$ as our state space, but whenever this is possible we use smaller spaces to formulate sharper results. This allows us to consider settings where the shifts are unbounded in one direction only. This occurs for example when considering delay equations. \\ 

The Fredholm properties of the linear operator associated to (\ref{eq:basissysteem}) have been described by Faye and Scheel in \cite{Faye2014}. We make heavy use of these properties here, 
continuing the program initiated in the bachelor thesis of Jin \cite{BachJin}, who considered autonomous
versions of \eqref{eq:basissysteem}. In such settings
it is possible to extend the techniques developed by Hupkes and Augeraud-V\'eron in \cite{HJHPOLLUTION} for MFDEs posed on Banach spaces. However, it is unclear at present how to generalize these methods to non-autonomous systems. \\

\paragraph{Splittings on the full line}

In \S\ref{section:existence}-\ref{section:fredholm}
we construct exponential splittings
for (\ref{eq:basissysteem}) on the full line.
Our main result essentially states that the decomposition
\eqref{eq:int:decomp:mfde} and the characterization 
\eqref{eq:PplusQHale}
remain valid for the state space $C_b(\R)$.
In addition, we explore
the Fredholm and continuity properties of the projection operators associated to the splitting
\eqref{eq:int:decomp:mfde}.
Our arguments in these sections are heavily based on the framework developed by  Mallet-Paret and Verduyn Lunel in \cite{MPVL}. However, the unbounded shifts raise some major technical challenges.\\

The primary complication is that the iteration scheme
used in \cite{MPVL} to establish the exponential decay of functions in $P(\tau)$ and $Q(\tau)$ breaks down. Indeed, the authors show that there exist $L > 0$
so that supremum of
the former solutions on half-lines $(-\infty, \tau_*]$ is halved each time one makes the replacement $\tau_* \mapsto \tau_* - L$. To achieve this, they exploit the fact that the behaviour of solutions on the latter interval does not `see' the behaviour at $\tau_*$.
This is no longer true for unbounded shifts
and required us to develop a novel iteration scheme
that is able to separate short-range from long-range effects.\\

A second major complication arises whenever continuous functions are approximated by $C^1$-functions. Indeed,
in \cite{MPVL} these approximations automatically have bounded derivatives, but in our case we can no longer assume that these functions live in $W^{1,\infty}(\R)$.
This prevents a direct application of the Fredholm theory in \cite{Faye2014}, forcing us to take a more involved approach to carefully isolate the regions where the unbounded derivatives occur.\\

The final obstacle is caused 
by the frequent use of the Ascoli-Arzela theorem
 in \cite{MPVL}. Indeed, in our setting
 we only obtain convergence on compacta instead of full uniform convergence. Fortunately, this can be circumvented relatively easily by using the exponential decay to provide the missing compactness at infinity. \\

 \paragraph{Splittings on the half line}
 We proceed in \S \ref{section:halfline}
 by constructing exponential dichotomies 
 for \eqref{eq:basissysteem}
 on the half-line $\R^+$. In particular, 
 for any $\tau \ge 0$ we establish the decomposition 
\begin{equation}
    \begin{array}{lcl}
         C_b(\R)&=&Q(\tau)\oplus R(\tau).
    \end{array}
\end{equation}
Here $Q(\tau)$ contains (shifted) exponentially decaying functions that satisfy \eqref{eq:basissysteem} on $[\tau, \infty)$, while (shifts of) functions in $R(\tau)$ satisfy \eqref{eq:basissysteem} on $[0, \tau]$. This generalizes the finite-range results
obtained by Hupkes and Verduyn Lunel in \cite{HJHLIN},
which we achieve by following a very similar strategy.\\

Besides the general complications discussed above, the main technical obstruction here is that the construction of half-line solutions to inhomogeneous versions of 
\eqref{eq:basissysteem} becomes rather delicate. Indeed,
the approach taken in \cite{HJHLIN} modifies the inhomogeneous terms \textit{outside} the `influence region' of the half-line of interest. However, in our setting here this region encompasses the whole line, forcing us to revisit the problem in a more elaborate\textemdash and technical\textemdash fashion. \\

\paragraph{Degeneracies}

In order to successfully exploit the characterization
(\ref{eq:PplusQHale}) in applications, it is essential to revisit the degeneracy issues related to the Hale inner product and the kernel elements of $\mathcal{L}^*$. Unfortunately, the absence of a `smallest' and `largest' shift in the infinite-range setting prevents an easy generalization of the invertibility criterion discussed above. 
We explore this crucial issue at length in {\S}\ref{section:Hale}.\\

In order to sketch some of the issues involved, 
we discuss the MFDE 
\begin{equation}\label{eq:linearnagumointro}
    \begin{array}{lcl}
         cu'(\sigma)&=&\sum\limits_{k=1}^\infty \gamma_k[u( \sigma+k)+u( \sigma-k)-2u( \sigma)]+g_u(\overline{u}_*( \sigma); a)u(\sigma),
    \end{array}
\end{equation}
which can be interpreted
as an infinite-range version of the 
MFDE (\ref{eq:linearnagumomfde}) that arose
by linearizing the Nagumo LDE around a travelling wave
$\overline{u}_*$. In particular, we again
assume the limits \eqref{eq:int:nag:front}. 
This MFDE fits into our framework provided that the coefficients $\gamma_k$ decay exponentially.\\

For the case $\gamma_k = e^{-k}$, we construct an explicit non-trivial function $\psi$ that satisfies $\langle \psi, \phi \rangle_{\tau} = 0$
 for each $\phi \in C_b(\R)$, where $\langle \cdot, \cdot \rangle_{\tau}$ denotes the appropriate Hale inner product for our setting. In particular, even for strictly positive coefficients there is no guarantee that the Hale inner product is nondegenerate. We also provide such examples for systems featuring convolution kernels.\\

 One way to circumvent this problem is to focus specifically on the kernel
 elements in \eqref{eq:PplusQHale}. If these can be chosen to be nonnegative along with the coefficients $\gamma_k$, then we are able to recover the
 relation between the dimension of $\Gamma(\tau)$
 in \eqref{eq:int:decomp:mfde} and the dimension of the kernel of the operator $\mathcal{L}^*$ associated
 to the adjoint of \eqref{eq:basissysteem}. Fortunately, such positivity conditions follow naturally for systems that admit a comparison principle.\\
 
 We also explore a second avenue that can be used without sign restrictions on the coefficients $\gamma_k$. This requires us to borrow some abstract functional analytic results.
 In particular, whenever the collection of sequences $\{\gamma_k\}_{k\geq N}$ obtained by taking $N \in \mathbb{N}$ spans an infinite dimensional subset of $\ell^2(\N;\C)$, we show that the Hale inner product is nondegenerate in a suitable sense. Fortunately,
 this rather abstract condition can often be made concrete. For example, we show that it can be enforced by imposing the Gaussian decay rate $\gamma_k \sim \mathrm{exp}[-k^2]$.

\paragraph{Acknowledgements.}
Both authors acknowledge support from the Netherlands Organization for Scientific Research (NWO) (grant 639.032.612).

\section{Main Results}\label{section:mainresults}

Our main results consider the integro-differential MFDE\footnote{
In the interest of readability
we
use $t$ as our main variable
throughout the remainder of this paper,
departing from the notation $\sigma$ that we used in {\S}\ref{sec:int}.
However, the reader should keep in mind
that this variable is related to a spatial quantity for most applications.}
\begin{equation}\label{ditishetprobleem}
    \begin{array}{lcl}
         \dot{x}(t)&=&\sum\limits_{j=-\infty}^\infty A_j(t)x(t+r_j)+\int\limits_{\R}\mathcal{K}(\xi;t)x(t+\xi)d\xi,
    \end{array}
\end{equation}
where we take $x\in \C^M$ for some integer $M\geq 1$. The set of scalars $\mathcal{R}:=\{r_j:j\in\Z\}\subset\R$ and the support of
$\mathcal{K}(\cdot;t)$ need not be bounded. 
In fact, we pick two constants
\begin{equation}
    -\infty\leq \rmin \le 0 \le \rmax\leq\infty,\qquad \rmin<\rmax
\end{equation}
in such a way that 
\begin{equation}\label{eq:def:rminrmax}
    \begin{array}{rcll}
         r_j&\in&\overline{(\rmin,\rmax)},\enskip &\text{for all }j\in\Z,\\[0.2cm]
         \mathrm{supp}\big(\mathcal{K}(\cdot;t)\big)&\subset &\overline{(\rmin,\rmax)},\enskip&\text{for all }t\in\R,
    \end{array}
\end{equation}
while $\abs{r_{\min}}$ and $\abs{r_{\max}}$ are as small as possible.
One readily sees that potential solutions to (\ref{ditishetprobleem}) must be defined on intervals that have a minimal length of $\rmax-\rmin$.\\

Naturally, one can always artificially increase the quantities $\abs{r_{\min}}$ and $\abs{r_{\max}}$ by adding
matrices $A_j= 0$ to (\ref{ditishetprobleem}) with large associated shifts $\abs{r_j} \gg 1$. However, we will see that
this only weakens the predictive power of our results by needlessly enlarging
the relevant state spaces. \\



A more general version of (\ref{ditishetprobleem}) might take the form
\begin{equation}\label{eq:moregeneralversion}
    \begin{array}{lcl}
         \dot{x}(t)&=&\int\limits_{\rmin}^{\rmax} d_\theta (t,\theta) x(t+\theta), 
    \end{array}
\end{equation}
where  $d_\theta(t,\theta)$ is an $M\times M$ matrix of finite Lebesgue-Stieltjes measures on $\overline{(\rmin,\rmax)}$ for each $t\in\R$. However, the adjoint of the system (\ref{eq:moregeneralversion}) is not always a system of similar type, so to avoid technical complications we will restrict ourselves to the system (\ref{ditishetprobleem}).\\

We now formulate our two main conditions on the coefficients
in (\ref{ditishetprobleem}),
which match those used in \cite{Faye2014}. As a preparation, 
we define the exponentially weighted space
\begin{equation}
\begin{array}{lcl}
L_{\eta}^1(\R;\C^{M\times M})&:=&\Big\{\mathcal{V}\in L^1(\R;\C^{M\times M})\Big|\nrm{e^{\eta|\cdot|}\mathcal{V}(\cdot)}_{L^1(\R;\C^{M\times M})}<\infty\Big\}
\end{array}
\end{equation}
for any $\eta>0$,
with its natural norm
\begin{equation}
\begin{array}{lcl}
\nrm{\mathcal{V}}_{\eta}&:=&\nrm{e^{\eta|\cdot|}\mathcal{V}(\cdot)}_{L^1(\R;\C^{M\times M})}.
\end{array}
\end{equation}
We note that the conditions on $\mathcal{R}$ below are not actual restrictions as long as the closure $\overline{\mathcal{R}}$ is countable. Indeed,
one can simply add the missing shifts to $\mathcal{R}$ and write $A_j=0$ for the associated matrix.

\begin{assumption}{A}{\text{H}}\label{aannamesconstanten} For each $j\in\Z$ the map $t\mapsto A_j(t)$ is bounded and belongs to $C^1(\R;\C^{M\times M})$. Moreover, there exists a constant $\tilde{\eta}>0$ for which the bound
\begin{equation}\label{expdecayscalars1}
    \begin{array}{lcl}
    \sum\limits_{j=-\infty}^\infty \nrm{A_j(\cdot)}_\infty e^{\tilde{\eta}|r_j|}&<&\infty
    \end{array}
\end{equation}
holds. In addition, the set $\mathcal{R}$ is closed with $0\in\mathcal{R}$.
\end{assumption}

\begin{assumption}{K}{\text{H}}\label{aannamesconvoluties} There exists a constant $\tilde{\eta}>0$
so that the 
following properties hold.
\begin{itemize}
\item The map $t\mapsto \mathcal{K}(\cdot;t)$ belongs to $C^1\big(\R;L_{\tilde{\eta}}^1(\R;\C^{M\times M}) \big)$.
\item The kernel $\mathcal{K}$ is localized
in the sense that
\begin{equation}\label{expdecayconvkernel}
\begin{array}{lcl}
\sup\limits_{t\in\R}\nrm{\mathcal{K}(\cdot;t)}_{\tilde{\eta}}+\sup\limits_{t\in\R}\nrm{\frac{d}{dt}\mathcal{K}(\cdot;t)}_{\tilde{\eta}}&<&\infty,\\[0.3cm]
\sup\limits_{t\in\R}\nrm{\mathcal{K}(\cdot;t-\cdot)}_{\tilde{\eta}}+\sup\limits_{t\in\R}\nrm{\frac{d}{dt}\mathcal{K}(\cdot;t-\cdot)}_{\tilde{\eta}}&<&\infty .
\end{array}
\end{equation}
\end{itemize}
\end{assumption}

Our third structural condition involves the behaviour
of the coefficients in (\ref{ditishetprobleem}) as $ t \to \pm \infty$.
Following \cite{MPA,Faye2014}, we say that the system (\ref{ditishetprobleem}) is asymptotically hyperbolic if the limits
\begin{equation}
    \begin{array}{lclcl}
        A_j(\pm\infty)&:=&\lim\limits_{t\rightarrow\pm\infty} A_j(t),
        \qquad \qquad 
        \mathcal{K}(\xi;\pm\infty)&:= &\lim\limits_{t\rightarrow\pm\infty} \mathcal{K}(\xi;t)
    \end{array}
\end{equation}
exist for each $j \in \Z$ and $\xi \in \R$, while the characteristic functions
\begin{equation}
    \begin{array}{lcl}
         \Delta^{\pm}(z)&=&z I-\int_{\R}\mathcal{K}(\xi;\pm\infty)e^{z\xi}d\xi-\sum\limits_{j=-\infty}^\infty A_j(\pm\infty)e^{zr_j}
    \end{array}
\end{equation}
associated to the limiting systems
\begin{equation}\label{autlimitsystems}
    \begin{array}{lcl}
         \dot{x}(t)&=&\sum\limits_{j=-\infty}^\infty A_j(\pm\infty)x(t+r_j)+\int\limits_{\R}\mathcal{K}(\xi;\pm\infty)x(t+\xi)d\xi
    \end{array}
\end{equation}
satisfy
\begin{equation}
    \begin{array}{lcl}
         \det\Delta^{\pm}(iy)&\neq &0
    \end{array}
\end{equation}
for all $y \in \R$. In fact, we require that these
limiting systems are approached in a summable fashion.

\begin{assumption}{H}{\text{H}}\label{aannameslimit} 
The system (\ref{ditishetprobleem}) is asymptotically
hyperbolic and satisfies the limits
\begin{equation}\label{eq:sterkelimitestimates}
\begin{array}{lcl}
\lim\limits_{t\rightarrow\pm\infty}\sum\limits_{j=-\infty}^\infty |A_j(t)-A_j(\pm\infty)|e^{\tilde{\eta}|r_j|}&=&0,\\[0.2cm]
\end{array}
\end{equation}
together with
\begin{equation}
    \begin{array}{lclcl}
        \lim\limits_{t\rightarrow\pm\infty}\nrm{\mathcal{K}(\cdot;t)-\mathcal{K}(\cdot;\pm\infty)}_{\tilde{\eta}}
        &=&0, \qquad \qquad 
        \lim\limits_{t\rightarrow\pm\infty}\nrm{\mathcal{K}(\cdot;t-\cdot)-\mathcal{K}(\cdot;\pm\infty)}_{\tilde{\eta}}
        &=& 0.
    \end{array}
\end{equation}
\end{assumption}

Bounded solutions to the system (\ref{ditishetprobleem}) can be interpreted
as kernel elements of the linear operator $\Lambda:W^{1,\infty}(\R)\rightarrow L^\infty(\R)$ that acts as
\begin{equation}\label{def:Lambda}
    \begin{array}{lcl}
         (\Lambda x)(t)&=&\dot{x}(t)-\sum\limits_{j=-\infty}^\infty A_j(t)x(t+r_j)-\int\limits_{\R}\mathcal{K}(\xi;t)x(t+\xi)d\xi.
    \end{array}
\end{equation}
We will write $\Lambda^*:W^{1,\infty}(\R)\rightarrow L^\infty(\R)$ 
for the formal adjoint of this operator, which is given by
\begin{equation}\label{def:Lambda^*}
    \begin{array}{lcl}
         (\Lambda^* y)(t)&=&-\dot{y}(t)-\sum\limits_{j=-\infty}^\infty A_j(t-r_j)^\dagger y(t-r_j)-\int\limits_{\R}\mathcal{K}(\xi;t-\xi)^\dagger y(t-\xi)d\xi,
    \end{array}
\end{equation}
using 
$\dagger $ 
to denote
the conjugate transpose of a matrix.
Indeed, one may readily verify the identity
\begin{equation}
    \begin{array}{lcl}
        \langle y, \Lambda x \rangle_{L^2(\R)} &=& \langle \Lambda^* y, x \rangle_{L^2(\R)}
    \end{array}
\end{equation}
whenever $x,y \in H^1(\R)$.\\

For convenience, we borrow the notation from \cite{MPVL,HJHLIN} and write
\begin{equation}
    \begin{array}{lclcl}
         \mathcal{B} 
         &=&\ker(\Lambda), 
         \qquad \qquad 
         \mathcal{B}^* &=& \ker(\Lambda^*).
    \end{array}
\end{equation}
The following result obtained by Faye and Scheel  describes
several useful Fredholm properties that link these kernels to the ranges
of the operators $\Lambda$ and $\Lambda^*$.

\begin{proposition}[{\cite[Thm. 2]{Faye2014}}]\label{prop:fredholmproperties} Assume that (\asref{aannamesconstanten}{\text{H}}), (\asref{aannamesconvoluties}{\text{H}}) and (\asref{aannameslimit}{\text{H}}) are satisfied. Then both the operators $\Lambda$ and $\Lambda^*$ are Fredholm operators. Moreover, the kernels and ranges satisfy the identities
\begin{equation}\label{eq:RangeLambda}
    \begin{array}{lcl}
         \Range(\Lambda)&=&\{h\in L^\infty(\R)\text{ }|\text{ }\int\limits_{-\infty}^\infty y(t)^\dagger h(t)dt=0\text{ for every }y\in \mathcal{B}^*\},\\[0.2cm]
         \Range(\Lambda^*)&=&\{h\in L^\infty(\R)\text{ }|\text{ }\int\limits_{-\infty}^\infty x(t)^\dagger h(t)dt=0\text{ for every }x\in \mathcal{B}\}
    \end{array}
\end{equation}
and the Fredholm indices can be computed by
\begin{equation}
    \begin{array}{lclcl}
         \mathrm{ind}(\Lambda)&=&-\mathrm{ind}(\Lambda^*)&=&\dim\mathcal{B}-\dim\mathcal{B}^*.
    \end{array}
\end{equation}
Finally, there exist constants $C>0$ and $0<\alpha\leq\tilde{\eta}$ 
so that the estimate
\begin{equation}
    \begin{array}{lcl}
         |b(t)|&\leq & Ce^{-\alpha |t|}\nrm{b}_\infty
    \end{array}
\end{equation}
holds for any $b\in \mathcal{B} \cup \mathcal{B}^*$ and any $t\in\R$.
\end{proposition}

\subsection{State spaces}

Let us introduce the intervals
\begin{equation}
    \begin{array}{lclcl}
        D_X &=& \overline{(\rmin,\rmax)},
        \qquad \qquad
        D_Y &=&\overline{(-\rmax,-\rmin)},
    \end{array}
\end{equation}
together with the state spaces
\begin{equation}
    \begin{array}{lclcl}
        X&=&C_b(D_X), \qquad \qquad
        Y&=&C_b(D_Y),
    \end{array}
\end{equation}
which contain bounded continuous functions that we measure with the supremum norm.
Suppose now that $x$ and $y$
are two bounded continuous functions that are defined
on (at least) the interval $t + D_X$ respectively $t + D_Y$.
We then write $x_t \in X$ and $y^t \in Y$ for the segments
\begin{equation}
    \begin{array}{lclcl}
        x_t(\theta)&=&x(t+\theta),
        \qquad \qquad 
        y^t(\theta)&=&y(t+\theta),
    \end{array}
\end{equation}
in which $\theta \in D_X$ respectively $\theta \in D_Y$.
This allows us to introduce the kernel segment spaces
\begin{equation}\label{eq:defPspaces0}
    \begin{array}{lcl}
         B(\tau)&=&\{\phi\in X\text{ }|\text{ } \phi=x_\tau\text{ for some }x\in\mathcal{B}\},\\[0.2cm]
         B^*(\tau)&=&\{\psi\in Y\text{ }|\text{ }\psi=y^{\tau}\text{ for some }y\in \mathcal{B}^*\}
    \end{array}
\end{equation}
for every $\tau \in \R$.
Observe that  $B(\tau)$ and
$B^*(\tau)$ are just shifted versions of $\mathcal{B}$ and $\mathcal{B}^*$
if $\rmin=-\infty$ and $\rmax=\infty$ both hold.\\

The Hale inner product \cite{HLVL} provides a useful coupling between $X$ and $Y$. The natural definition in the current setting is given by
\begin{equation}\label{eq:defHale}
    \begin{array}{lcl}
         \ip{\psi,\phi}_t&=&\psi(0)^\dagger \phi(0)-\sum\limits_{j=-\infty}^\infty \int\limits_0^{r_j} \psi(s-r_j)^\dagger A_j(t+s-r_j)\phi(s)ds\\[0.2cm]
         &&\qquad -\int\limits_{\R}\int\limits_0^r \psi(s-r)^\dagger \mathcal{K}(r;t+s-r)\phi(s)dsdr
    \end{array}
\end{equation}
for any pair $( \phi, \psi) \in X \times Y$. Note that, by decreasing $\tilde{\eta}$ if necessary, we can strengthen (\ref{expdecayscalars1}) to obtain
\begin{equation}\label{expdecayscalars2}
    \begin{array}{lcl}
         \sum\limits_{j=-\infty}^\infty \nrm{A_j(\cdot)}_\infty |r_j|e^{\tilde{\eta}|r_j|}&<&\infty.
    \end{array}
\end{equation}
Together with (\ref{expdecayconvkernel}), 
this ensures that the Hale inner product is well-defined. In 
Lemma \ref{lemma:diffhale} below we verify the identity
\begin{equation}\label{eq:fundHale}
    \begin{array}{lcl}
         \frac{d}{dt}\ip{y^t,x_t}_t&=&y^\dagger (t)[\Lambda x](t)+[\Lambda^*y](t)^\dagger x(t)
    \end{array}
\end{equation}
for $x,y\in W^{1, \infty}(\R)$,
which indicates that the Hale inner product can be seen as the duality pairing 
between $\Lambda$ and $\Lambda^*$.\\

An important role in the sequel
is reserved for the subspaces
\begin{equation}\label{eq:xperptau}
\begin{array}{lcl}
    X^\perp(\tau) &= &
    \{\phi\in X\text{ }|\text{ }\ip{\psi,\phi}_\tau=0\text{ for every }\psi\in B^*(\tau)\},
    \end{array}
\end{equation}
which have finite codimension
\begin{equation}\label{eq:4.12}
    \begin{array}{lclclcl}
         \beta(\tau)
          &:=& \mathrm{codim}_{X} X^\perp(\tau)
          &\le& \dim B^*(\tau) &\le& \dim \mathcal{B}^*.
    \end{array}
\end{equation}
In the ODE case $r_{\min} = r_{\max} = 0$, so
one readily concludes that $\beta(\tau) = \dim \mathcal{B}^*$. However, in the present setting it is possible for the Hale inner product to be degenerate or for kernel elements to vanish on large intervals. In these cases the first respectively second inequality in (\ref{eq:4.12}) could become strict.\\

In the finite range setting of \cite{MPVL}, the authors ruled out these degeneracies
by imposing an atomic condition on the matrices $\{A_j\}$ corresponding to the shifts $\rmin$ and $\rmax$. However, there
is no obvious way to generalize this condition
when $\abs{r_{\min}}$ or $r_{\max}$ are infinite. As an alternative,
some of our results require the following technical assumption.

\begin{assumption}{Ker}{\text{H}}\label{aannames0ophalflijn} 
Consider any non-zero $d\in\mathcal{B}\cup\mathcal{B}^*$ and $\tau\in\R$. Then $d$ does not vanish on $(-\infty,\tau]$ and also does not vanish on $[\tau,\infty)$.
\end{assumption}

A similar assumption was used in \cite[Assumption H3(iii)]{HJHNLS}, where the authors remove the $\abs{r_{\min}} = r_{\max}$ restriction from the exponential dichotomy constructions
in \cite{harterich}. However, this condition
is naturally much harder to verify than the previous atomicity condition. We explore this issue at length in 
{\S}\ref{section:Hale}, 
where we present several scenarios under which (\asref{aannames0ophalflijn}{\text{H}}) can be verified.\\

We highlight one of these scenarios in the result below, which requires
sign conditions on elements of $\mathcal{B}$ and $\mathcal{B}^*$. Fortunately, for a large class of systems\textemdash including the linearization (\ref{eq:linearnagumointro}) of the Nagumo LDE\textemdash these are known consequences of the comparison principle.

\begin{proposition}[{see Prop. \ref{assumptionimplications}}]\label{prop:4.16:positiveweaker} Assume that (\asref{aannamesconstanten}{\text{H}}), (\asref{aannamesconvoluties}{\text{H}}) and (\asref{aannameslimit}{\text{H}}) are satisfied. Assume furthermore that there exists $K_{\mathrm{const}} \in \Z_{\geq 1}$ for which the following structural conditions are satisfied.
\begin{itemize}
    \item[(a)]{
      We have $r_j=j$ for $j\in\Z$, which implies $r_{\min} = - \infty$ and $\rmax=\infty$.
    }
    \item[(b)]{
      The function $A_j(\cdot)$ is constant and positive definite whenever $|j|\geq K_{\mathrm{const}}$.
    }
    \item[(c)]{
      For any $|\xi|\geq K_{\mathrm{const}}$ 
      the function $\mathcal{K}(\xi; \cdot)$
      is constant and positive definite.
    }
    \item[(d)]{
      We either have $\mathcal{B} = \{0 \}$
      or $\mathcal{B} = \mathrm{span}\{b\}$ for some nonnegative function $b$.
      The same holds for $\mathcal{B}^*$.
    }
\end{itemize}
Then the non-triviality condition (\asref{aannames0ophalflijn}{\text{H}}) is satisfied.
\end{proposition}

In {\S}\ref{section:halfline}-\ref{section:Hale} we explore some of the consequences of (\asref{aannames0ophalflijn}{\text{H}}). In addition, we propose weaker conditions under which equality holds for one
or both of the inequalities in
(\ref{eq:4.12}). However, for now we simply state the following result.

\begin{corollary}[{cf. \cite[Cor. 4.7]{MPVL}, see \S \ref{section:Hale}}]\label{cor:4.7:a} Assume that (\asref{aannamesconstanten}{\text{H}}), (\asref{aannamesconvoluties}{\text{H}}), (\asref{aannameslimit}{\text{H}}) 
and (\asref{aannames0ophalflijn}{\text{H}})
are all satisfied. Then the identities
\begin{equation}\label{eq:4.21}
    \begin{array}{lclclcl}
         \dim B(\tau)&=&\dim \mathcal{B},
         \qquad \qquad
         \beta(\tau)& = &\dim B^*(\tau)& =& \dim\mathcal{B}^*
    \end{array}
\end{equation}
hold for every $\tau\in\R$.
\end{corollary}

\subsection{Exponential dichotomies on \texorpdfstring{$\R$}{R}}

We now set out to describe our exponential splittings for \eqref{ditishetprobleem} on the full line $\R$. To this end,
we introduce the intervals
\begin{equation}
    \begin{array}{lclcl}
         D_\tau^\ominus &=&\overline{(-\infty,\tau+\rmax)},
\qquad \qquad 
D_\tau^\oplus 
&=&\overline{(\tau+\rmin,\infty)}
    \end{array}
\end{equation}
for each $\tau\in\R$. Following
the notation in
\cite{MPVL,HJHLIN}, this allows us to define
the solution spaces
\begin{equation}\label{eq:defPspaces1}
    \begin{array}{lcl}
         \mathcal{P}(\tau)&=&\{x\in C_b( D_{\tau}^\ominus ) \text{ }|\text{ }x\text{ is a bounded solution of (\ref{ditishetprobleem}) on }(-\infty,\tau]\},\\[0.2cm]
         \mathcal{Q}(\tau)&=&\{x\in C_b( D_{\tau}^\oplus )\text{ }|\text{ }x\text{ is a bounded solution of (\ref{ditishetprobleem}) on }[\tau,\infty)\},
    \end{array}
\end{equation}
together with the associated initial segments
\begin{equation}\label{eq:defPspaces2}
    \begin{array}{lcl}
         P(\tau)&=&\{\phi\in X\text{ }|\text{ } \phi=x_\tau\text{ for some }x\in\mathcal{P}(\tau)\},\\[0.2cm]
         Q(\tau)&=&\{\phi\in X\text{ }|\text{ } \phi=x_\tau\text{ for some }x\in\mathcal{Q}(\tau)\}.\\[0.2cm]
    \end{array}
\end{equation}
For $\tau\in\R$ we call $x\in \mathcal{P}(\tau)$ a left prolongation of an element $\phi=x_\tau\in P(\tau)$, with a similar definition for right prolongations. Note that, if $\rmin=-\infty$, each $\phi\in P(\tau)$ is simply a translation of a function in $\mathcal{P}(\tau)$. The corresponding result holds for $Q(\tau)$ and $\mathcal{Q}(\tau)$ if $\rmax=\infty$.\\

Again following \cite{MPVL}, we also work with the spaces
\begin{equation}\label{defhatmathcalspaces}
    \begin{array}{lcl}
        \widehat{\mathcal{P}}(\tau)&=&\{x\in\mathcal{P}(\tau)\text{ }|\text{ }\int\limits_{-\infty}^{\tau+\rmax}y(t)^\dagger x(t)dt=0\text{ for every }y\in\mathcal{B}\},\\[0.2cm]
        \widehat{\mathcal{Q}}(\tau)&=&\{x\in\mathcal{Q}(\tau)\text{ }|\text{ }\int\limits_{\tau+\rmin}^\infty y(t)^\dagger x(t)dt=0\text{ for every }y\in\mathcal{B}\},
    \end{array}
\end{equation}
together with
\begin{equation}
    \begin{array}{lcl}
        \widehat{P}(\tau)&=&\{\phi\in X\text{ }|\text{ } \phi=x_\tau\text{ for some }x\in\widehat{\mathcal{P}}(\tau)\},\\[0.2cm]
         \widehat{Q}(\tau)&=&\{\phi\in X\text{ }|\text{ } \phi=x_\tau\text{ for some }x\in\widehat{\mathcal{Q}}(\tau)\}.
    \end{array}
\end{equation}
The integrals in (\ref{defhatmathcalspaces}) convergence since functions in $\mathcal{B}$ decay exponentially. Finally, we write
\begin{equation}
    \begin{array}{lclcl}
        S(\tau)&=&P(\tau)+Q(\tau),\qquad \widehat{S}(\tau)=\widehat{P}(\tau)+\widehat{Q}(\tau).
    \end{array}
\end{equation}

Our first two results here
provide exponential decay estimates
for functions in $\widehat{\mathcal{P}}(\tau)$ and $\widehat{\mathcal{Q}}(\tau)$, together with
a direct sum decomposition for $S(\tau)$. In addition, we show that  the latter space 
can be identified with $X^\perp(\tau)$ from (\ref{eq:xperptau}).
We remark that the structure of these results
matches their counterparts from \cite{HLVL}
almost verbatim.\\

\begin{theorem}[{cf. \cite[Thm. 4.2]{MPVL}, see \S \ref{section:existence}}]\label{thm:expdecay} Assume that (\asref{aannamesconstanten}{\text{H}}), (\asref{aannamesconvoluties}{\text{H}}) and (\asref{aannameslimit}{\text{H}}) are satisfied and choose a sufficiently large $\tau_*>0$. Then there exist constants $ K_{\mathrm{dec}} >0$ and $\alpha>0$
so that for any 
$\tau\leq -\tau_*$
and $p\in\mathcal{P}(\tau)$
we have the bound
\begin{equation}\label{eq:4.8equiv:p}
    \begin{array}{lcl}
        |p(t)|+|\dot{p}(t)|&\leq &
        K_{\mathrm{dec}} e^{\alpha (t-\tau)}\nrm{p_\tau}_\infty, \qquad \qquad t\leq \tau,\\[0.2cm]
     \end{array}
\end{equation}
while for any $\tau\geq \tau_*$
and $q\in\mathcal{Q}(\tau)$
we have the corresponding estimate
\begin{equation}\label{eq:4.8equiv:q}
    \begin{array}{lcl}
         |q(t)|+|\dot{q}(t)|&\leq & K_{\mathrm{dec}}e^{-\alpha(t-\tau)}\nrm{q_\tau}_\infty,\qquad \qquad t\geq \tau .
    \end{array}
\end{equation}


In addition, the bounds  (\ref{eq:4.8equiv:p})-(\ref{eq:4.8equiv:q}) also hold
for any $p\in \widehat{\mathcal{P}}(\tau)$ and $q\in \widehat{\mathcal{Q}}(\tau)$,
now without any restriction on the value of $\tau\in\R$,
but with possibly different values of $K_{\mathrm{dec}}$ and $\alpha$.
\end{theorem}

\begin{theorem}[{cf. \cite[Thm. 4.3]{MPVL}, see \S \ref{section:existence}}]\label{thm:4.3equiv} Assume that (\asref{aannamesconstanten}{\text{H}}), (\asref{aannamesconvoluties}{\text{H}}) and (\asref{aannameslimit}{\text{H}}) are satisfied. For each $\tau\in\R$ the spaces $P(\tau)$, $Q(\tau)$, $S(\tau)$ and their counterparts $\widehat{P}(\tau)$, $\widehat{Q}(\tau)$, $\widehat{S}(\tau)$ are all closed subspaces of $X$. Moreover, we have the identities

\begin{equation}\label{eq:4.10}
    \begin{array}{lclcl}
         P(\tau)&=&\widehat{P}(\tau)\oplus B(\tau),\qquad
         Q(\tau)&=&\widehat{Q}(\tau)\oplus B(\tau),\\[0.2cm]
         \widehat{S}(\tau)&=&\widehat{P}(\tau)\oplus \widehat{Q}(\tau),\qquad
         S(\tau)&=&\widehat{S}(\tau)\oplus B(\tau)\\[0.2cm]
         &&&=&\widehat{P}(\tau)\oplus \widehat{Q}(\tau)\oplus B(\tau).
    \end{array}
\end{equation}
Finally, we have the identification 
\begin{equation}\label{eq:4.11}
    \begin{array}{lcl}
         S(\tau)&=&X^\perp(\tau),
    \end{array}
\end{equation}
where $X^\perp(\tau)$ is defined in (\ref{eq:xperptau}).
\end{theorem}

However, these theorems provide no information on how the spaces $P(\tau)$ and $Q(\tau)$ depend on $\tau$.
In order to address this issue,
we need to study the projections from the state space $X$ onto the factors $\widehat{P}(\tau)$ and $\widehat{Q}(\tau)$ using the decomposition in (\ref{eq:4.10}). To be more precise,
for a fixed $\tau_0 \in \R$ we write
\begin{equation}\label{eq:4.14}
    \begin{array}{lcl}
         X&=&\widehat{P}(\tau_0)\oplus \widehat{Q}(\tau_0)\oplus \Gamma
    \end{array}
\end{equation}
for a suitable finite dimensional subspace $\Gamma \subset X$. This allows us define projections $\Pi_{\widehat{P}}$ and $\Pi_{\widehat{Q}}$ onto the factors $\widehat{P}(\tau_0)$ respectively $\widehat{Q}(\tau_0)$.\\

In addition, we are interested in the limiting behaviour as $ \tau \to \pm \infty$. To this end, we apply 
Theorem \ref{thm:4.3equiv} to the two limiting systems (\ref{autlimitsystems}), 
which leads to the decompositions
\begin{equation}
    \begin{array}{lclcl}
         X&=&P(-\infty)\oplus Q(-\infty)&=&P(\infty)\oplus Q(\infty).
    \end{array}
\end{equation}
We write $\overleftarrow{\Pi}_P$ and $\overleftarrow{\Pi}_Q$ for the projections onto the factors $P(-\infty)$ and $Q(-\infty)$ respectively, together with $\overrightarrow{\Pi}_P$ and $\overrightarrow{\Pi}_Q$ for the projections onto the factors $P(\infty)$ and $Q(\infty)$.\\

\begin{theorem}[{cf. \cite[Thm. 4.6]{MPVL}, see \S \ref{section:fredholm}}]\label{theorem4.6} Assume that (\asref{aannamesconstanten}{\text{H}}), (\asref{aannamesconvoluties}{\text{H}}) and (\asref{aannameslimit}{\text{H}}) are satisfied. Then the spaces $\widehat{P}(\tau)$, $\widehat{Q}(\tau)$ and $\widehat{S}(\tau)$ vary upper semicontinuously with $\tau$,
while the quantities $\dim B(\tau)$ and $\beta(\tau)$ vary lower semicontinuously with $\tau$.

In particular, fix $\tau_0\in\R$ and consider any $\tau$ sufficiently close to $\tau_0$. Then the restrictions
\begin{equation}\label{eq:4.16}
    \begin{array}{lcl}
         \Pi_{\widehat{P}}:\widehat{P}(\tau)&\rightarrow &\Pi_{\widehat{P}}\big(\widehat{P}(\tau)\big)\subset \widehat{P}(\tau_0),\\[0.2cm]
         \Pi_{\widehat{Q}}:\widehat{Q}(\tau)&\rightarrow &\Pi_{\widehat{Q}}\big(\widehat{Q}(\tau)\big)\subset \widehat{Q}(\tau_0)
    \end{array}
\end{equation}
of the projections associated
to the decomposition \eqref{eq:4.14}
are isomorphisms onto their ranges, which are closed.
 Moreover, the norms satisfy
\begin{equation}\label{eq:4.17}
    \begin{array}{lclcl}
         \lim\limits_{\tau\rightarrow\tau_0}\nrm{I-\Pi_{\widehat{P}}|_{\widehat{P}(\tau)}}&=&0,\qquad \qquad \lim\limits_{\tau\rightarrow\tau_0}\nrm{I-\Pi_{\widehat{Q}}|_{\widehat{Q}(\tau)}}&=&0,
    \end{array}
\end{equation}
in which $I$ denotes the inclusion of $\widehat{P}(\tau)$ or $\widehat{Q}(\tau)$ into $X$. 

In addition, 
we have the identities
\begin{equation}\label{eq:4.20}
    \begin{array}{lclcllcl}
         \overleftarrow{\Pi}_P\big(P(\tau)\big)&=&P(-\infty),\\[0.2cm]
         \overrightarrow{\Pi}_Q\big(Q(\tau)\big)&=&Q(\infty),
    \end{array}
\end{equation}
for sufficiently negative values of $\tau$ in the first line of (\ref{eq:4.20}) and for sufficiently positive values of $\tau$ in the second line of (\ref{eq:4.20}). The associated norms satisfy the limits
\begin{equation}\label{eq:4.19}
    \begin{array}{lclcl}
         \lim\limits_{\tau\rightarrow-\infty}\nrm{I-\overleftarrow{\Pi}_P|_{P(\tau)}}&=&0,\qquad \qquad \lim\limits_{\tau\rightarrow\infty}\nrm{I-\overrightarrow{\Pi}_Q|_{Q(\tau)}}&=&0.
    \end{array}
\end{equation}
\end{theorem}

These results can be strengthened if we also assume that (\asref{aannames0ophalflijn}{\text{H}}) holds. Indeed, Corollary 
\ref{cor:4.7:a} implies that the codimension of $S(\tau)$ remains constant. This can be leveraged to obtain the following continuity properties.

\begin{corollary}[{cf. \cite[Cor. 4.7]{MPVL}, see \S \ref{section:Hale}}]\label{cor:4.7:b} Assume that (\asref{aannamesconstanten}{\text{H}}), (\asref{aannamesconvoluties}{\text{H}}),  (\asref{aannameslimit}{\text{H}}) 
and (\asref{aannames0ophalflijn}{\text{H}})
are all satisfied. Then the spaces
$\widehat{P}(\tau)$ and $\widehat{Q}(\tau)$ vary continuously with $\tau$, i.e. the projections $\Pi_{\widehat{P}}$ and $\Pi_{\widehat{Q}}$ from (\ref{eq:4.16}) are isomorphisms onto $\widehat{P}(\tau_0)$ and $\widehat{Q}(\tau_0)$ respectively. The same conclusion holds
for their counterparts $P(\tau)$ and $Q(\tau)$.
\end{corollary}

\subsection{Exponential dichotomies on half-lines}
In many applications it is useful to consider exponential dichotomies on half-lines such as $[0,\infty)$, instead of the full line. Our main goal here is to show to prove the natural generalisation of Theorem \ref{thm:4.3equiv} to this half-line setting, along
the lines of the results in \cite{HJHLIN}.\\

In particular, we set out to obtain  decompositions of the form
\begin{equation}
\label{eq:mr:decomp:half:line:X}
    \begin{array}{lcl}
         X&=&Q(\tau)\oplus R(\tau),
    \end{array}
\end{equation}
where $Q(\tau)$ is defined in (\ref{eq:defPspaces2}) and segments in $R(\tau)$ should be `extendable' to solve \eqref{ditishetprobleem} on $[0, \tau]$.  Since this is a finite interval however there is no longer a `canonical' definition for $R(\tau)$. In fact, we define these spaces in a indirect fashion, by constructing appropriate subsets
\begin{equation}
\label{eq:mr:incl:cal:r:tau}
    \mathcal{R}(\tau) \subset
    \{ r \in C_b( D^\ominus_\tau ) \text{ }|\text{ } r \hbox{ is a bounded solution of \eqref{ditishetprobleem} on } [0, \tau] \}
\end{equation}
and writing
\begin{equation}
\label{eq:mr:def:r:tau}
    \begin{array}{lcl}
        R(\tau)&=&\{\phi\in X\text{ }|\text{ } \phi=x_\tau\text{ for some }x\in\mathcal{R}(\tau)\}.
    \end{array}
\end{equation}

In order to achieve this, we exploit
continuity properties for the projection
operators that are stronger than those
obtained in Theorem \ref{theorem4.6}. In particular, we again impose the non-triviality
condition (\asref{aannames0ophalflijn}{\text{H}}). However, we explain in {\S}\ref{section:halfline} how this condition can be weakened slightly. For example, we need less information concerning the kernel space $\mathcal{B}$ to apply our construction.

\begin{theorem}[{cf. \cite[Thm. 4.1]{HJHLIN}, see \S \ref{section:halfline}}]\label{thm:lin:4.1}
Assume that (\asref{aannamesconstanten}{\text{H}}), (\asref{aannamesconvoluties}{\text{H}}), (\asref{aannameslimit}{\text{H}}) and (\asref{aannames0ophalflijn}{\text{H}}) are satisfied. Then for every $\tau\geq 0$ there exists a closed subspace $\mathcal{R}(\tau) \subset C_b(D_{\tau}^\ominus)$
that satisfies the inclusion \eqref{eq:mr:incl:cal:r:tau} together
with the following properties.
\begin{enumerate}[label=(\roman*)]
    \item Recalling the spaces
    (\ref{eq:defPspaces2})
    and \eqref{eq:mr:def:r:tau},
    the splitting \eqref{eq:mr:decomp:half:line:X}
    holds for every $\tau \ge 0$.
%
 
    \item There exist constants $K_{\mathrm{dec}}>0$ and $\alpha>0$ so that the exponential estimate
\begin{equation}
    \begin{array}{lclll}
         |x(t)|&\leq &K_{\mathrm{dec}}e^{-\alpha |t-\tau|}\nrm{x_\tau}_\infty
    \end{array}
\end{equation}
holds for every $x\in \mathcal{R}(\tau)$ and every pair $0\leq t\leq \tau$.
    \item The spaces $R(\tau)$ are invariant, in the sense that $x_t \in R(t)$ holds whenever
    $x \in \mathcal{R}(\tau)$ and $0\leq t\leq \tau$.
    The corresponding statement holds for the spaces $Q(\tau)$.
    \item The projections $\Pi_{Q(\tau)}$ and $\Pi_{R(\tau)}$ associated to the splitting
    \eqref{eq:mr:decomp:half:line:X}
    depend continuously on $\tau\geq 0$.
    In addition, there exists a constant $C\geq 0$ so that the uniform bounds $\nrm{\Pi_{Q(\tau)}}\leq C$ and $\nrm{\Pi_{R(\tau)}}\leq C$ hold for all $\tau\geq 0$.
\end{enumerate}
\end{theorem}

\section{The existence of exponential dichotomies}\label{section:existence}

Our goal in this section is to establish Theorems \ref{thm:expdecay}-\ref{thm:4.3equiv}. The strategy that we follow is heavily based on \cite{MPVL}, allowing us to simply refer to the results there  from time to time. However, the unbounded shifts force us to develop an alternative approach at several key points in the analysis. We have therefore structured this section in such a way that these modifications are highlighted.\\

The first main task is to show 
that functions in the spaces $\mathcal{P}(\tau)$ and $\mathcal{Q}(\tau)$, together with their derivatives, decay exponentially in a uniform fashion. When the shifts are unbounded, the methods developed in \cite{MPVL} can no longer be used to establish this exponential decay. In particular, 
the bound \eqref{eq:4.24equivalt} below was obtained
in \cite{MPVL}, but one cannot 
simply make the replacement
$r_{\max} \to \infty$ and still recover the desired exponential decay of solutions.
Indeed, the iterative scheme
in \cite{MPVL} breaks down, forcing us to use a different approach.\\

The key ingredient is to show that the cumulative influence of the large shifts decays exponentially. The following preliminary estimate will help us to quantify this. 

\begin{lemma}\label{thm:expdecaybewijs0} Assume that (\asref{aannamesconstanten}{\text{H}}), (\asref{aannamesconvoluties}{\text{H}}) and (\asref{aannameslimit}{\text{H}}) are satisfied.
Then there exist three constants $(p, K_{\mathrm{exp}},\alpha)\in\R_{>0}^3$ for which
the bound
\begin{equation}\label{eq:defp}
    \begin{array}{lcl}
         \sum\limits_{r_j\geq |t|}\abs{A_j(s)} e^{\alpha |r_j|} 
         + 
         \int\limits_{|t|}^{\infty}
         \abs{\mathcal{K}(\xi;s)}e^{\alpha|\xi|}d\xi
         &\leq &  K_{\mathrm{exp}} e^{-2\alpha |t|}
    \end{array}
\end{equation}
holds for all $t<-p$ and all $s\in\R$. In addition, if $\rmax<\infty$, then we can pick $p=\rmax$.
\end{lemma}
\textit{Proof.} Suppose first that $\rmax=\infty$. Setting $\alpha=\frac{\tilde{\eta}}{3}$,  
we can derive from (\ref{expdecayscalars1}) that
\begin{equation}
  \label{eq:exp:bnds:A:K:prlm}
    \begin{array}{lcl}
         \sum\limits_{r_j\geq |t|}\nrm{A_j(\cdot)}_\infty e^{\alpha |r_j|}&\leq &e^{-2\alpha t}\sum\limits_{r_j\geq |t|}\nrm{A_j(\cdot)}_\infty e^{\tilde{\eta} |r_j|}\\[0.2cm]
         &\leq &e^{-2\alpha t}\sum\limits_{j=-\infty}^\infty\nrm{A_j(\cdot)}_\infty e^{\tilde{\eta} |r_j|}
    \end{array}
\end{equation}
for $|t|$ sufficiently large. The second term in (\ref{eq:defp}) can be bounded
in the same fashion using (\ref{expdecayconvkernel}). If $\rmax<\infty$ then (\ref{eq:defp}) follows trivially for $p=\rmax$, since the left-hand side is always zero for $t<-p$ and $s\in\R$.\qed\\

Our first main result generalizes the bound 
(\ref{eq:4.24equivalt}) to the setting
where $r_{\max} = \infty$. This is achieved by splitting the relevant interval 
$[\tau,\infty)$ into  two parts
$[\tau,\tau+p]$ and $[\tau+p,\infty)$
that we analyze separately. We use the 
ideas from \cite{MPVL} to study the first part,
while careful estimates involving
\eqref{eq:defp} allow us to control
the contributions from the unbounded second interval.

\begin{proposition}
\label{thm:expdecaybewijs1} Assume that (\asref{aannamesconstanten}{\text{H}}), (\asref{aannamesconvoluties}{\text{H}}) and (\asref{aannameslimit}{\text{H}}) are satisfied, recall the 
constants $(p, K_{\mathrm{exp}},\alpha)\in\R_{>0}^3$
from Lemma \ref{thm:expdecaybewijs0}
and pick a sufficiently negative $\tau_- \ll -1$. Then there exists a constant $\sigma > 0$
so that for each $\tau\leq \tau^-$ and each $x\in\mathcal{P}(\tau)$ we have the bound
\begin{equation}\label{eq:4.24equiv}
    \begin{array}{lcl}
         |x(t)|&\leq &\max\Big\{\frac{1}{2}\sup\limits_{s\in(-\infty,\tau+p]}|x(s)|,
         \, \, K_{\mathrm{exp}}\sup\limits_{s\in[p+\tau,\infty)}e^{-\alpha(s-t)}|x(s)|\Big\},\qquad t\leq -\sigma +\tau
    \end{array}
\end{equation}
when $r_{\max} = \infty$, or alternatively
\begin{equation}\label{eq:4.24equivalt}
    \begin{array}{lcl}
         |x(t)|&\leq &\frac{1}{2}\sup\limits_{s\in(-\infty,\tau+r_{\max}]}|x(s)|,
         \qquad t\leq -\sigma +\tau
    \end{array}
\end{equation}
when $r_{\max} < \infty$.
The same\footnote{Naturally, one may need
to change the value of the constant $\sigma > 0$.} bounds hold for $x\in \widehat{\mathcal{P}}(\tau)$, but now any $\tau\in\R$ is permitted.
\end{proposition}

The second main complication occurs 
when one tries to mimic the approach
in \cite{MPVL} to study the properties of 
$S(\tau)$. Although it is relatively
straightforward to show that this space
is closed and has finite codimension in $X$,
the explicit description (\ref{eq:4.11})
for $S(\tau)$ is much harder to obtain.
The arguments in \cite{MPVL} approximate
elements of $X^\perp(\tau)$ by $C^1$-smooth functions and apply the Fredholm operator $\Lambda$ to (extensions of) these approximants. However, when $D_X$ is unbounded 
this approach breaks down, because
$C^1$-smooth functions in $X$ need not have
a bounded derivative. One can hence no longer
directly appeal to the useful Fredholm properties of $\Lambda$.\\

Our second main result provides an alternative approach that circumvents these difficulties.
The novel idea is that we split such problematic functions into two parts that both
confine the regions where the derivatives are unbounded to a half-line. This turns out to be sufficient to allow the main spirit of the analysis in \cite{MPVL} to proceed.

\begin{proposition}\label{thm:4.3equivbewijs3} Assume that (\asref{aannamesconstanten}{\text{H}}), (\asref{aannamesconvoluties}{\text{H}}) and (\asref{aannameslimit}{\text{H}}) are satisfied. Fix $\tau\in\R$ and let $X^\perp(\tau)$ be given by (\ref{eq:xperptau}). Then there exists a dense subset $D\subset X^\perp(\tau)$ with $D\subset S(\tau)$.
\end{proposition}

Besides these two main obstacles, we encounter smaller technical issues at many points during our analysis. For example, the lack of full uniform convergence on unbounded intervals from the Ascoli-Arzela theorem requires significant attention. In addition, manipulations
involving the Hale inner product on unbounded domains raise subtle convergence issues that must be addressed.

\subsection{Preliminaries}

In this subsection, we collect several
preliminary properties satisfied by the spaces introduced in (\ref{eq:defPspaces0}), (\ref{eq:defPspaces1}) and (\ref{eq:defPspaces2}). In particular,
we discuss whether functions
in $P(\tau)$ or $Q(\tau)$ have unique extensions in $\mathcal{P}(\tau)$ and $\mathcal{Q}(\tau)$ and study the intersection
$P(\tau) \cap Q(\tau)$.

\begin{lemma}\label{lemma:elementaryproperties}  Assume that (\asref{aannamesconstanten}{\text{H}}), (\asref{aannamesconvoluties}{\text{H}}) and (\asref{aannameslimit}{\text{H}}) are satisfied and fix $\tau\in\R$. Then the spaces defined
in {\S}\ref{section:mainresults}  have the following properties.
\begin{enumerate}[label=(\roman*)]
    \item 
    We have the inequalities $\dim B(\tau)\leq \dim\mathcal{B} < \infty$ and $\dim B^*(\tau)\leq \dim \mathcal{B}^* < \infty$. In addition, if $|\rmin|=\rmax=\infty$, then $\dim B(\tau)= \dim\mathcal{B}$ and $\dim B^*(\tau)=\dim \mathcal{B}^*$.
    \item The inclusions $\widehat{\mathcal{P}}(\tau)\subset \mathcal{P}(\tau)$, $\widehat{\mathcal{Q}}(\tau)\subset \mathcal{Q}(\tau)$, $\widehat{P}(\tau)\subset P(\tau)$ and $\widehat{Q}(\tau)\subset Q(\tau)$ have finite codimension of at most $\dim\mathcal{B}$.
    \item We have $B(\tau)=P(\tau)\cap Q(\tau)$.
\end{enumerate}
\end{lemma}
\textit{Proof.} Items (i) and (ii) are clear from their definition and Proposition \ref{prop:fredholmproperties}. For item (iii) we note that the inclusion $B(\tau)\subset P(\tau)\cap Q(\tau)$ is trivial. Conversely, for $\phi\in P(\tau)\cap Q(\tau)$ we pick $x\in \mathcal{P}(\tau)$ and $y\in \mathcal{Q}(\tau)$ with $\phi=x_\tau=y_\tau$, so that $x=y$ on $D_X+\tau$. This allows us to consider the function $z$ that is defined 
on the real line by
\begin{equation}
    \begin{array}{lclcl}
        z(t)&=&\begin{cases} x(t), \enskip & t\leq \rmax+\tau\\[0.2cm] y(t),\enskip & t\geq \rmin+\tau.\end{cases}
    \end{array}
\end{equation}
It is now easy to see that $z\in\mathcal{B}$, which implies $\phi\in B(\tau)$.\qed\\

\begin{lemma}
\label{prop:uniqueleftprol}  Assume that (\asref{aannamesconstanten}{\text{H}}), (\asref{aannamesconvoluties}{\text{H}}) and (\asref{aannameslimit}{\text{H}}) are satisfied. Then there exists $\mu^-\in(-\infty,\infty]$ such that every $\phi\in P(\tau)$ with $\tau < \mu^-$ has a unique left prolongation in $\mathcal{P}(\tau)$. Similarly, there exists $\mu^+\in[-\infty,\infty)$ such that  every $\phi\in Q(\tau)$ with $\tau>\mu^+$ has a unique right prolongation in $\mathcal{Q}(\tau)$.
On the other hand, any element of $\widehat{P}(\tau)$ and $\widehat{Q}(\tau)$
has a unique left respectively right prolongation, this time for any $\tau\in\R$.
\end{lemma}
\textit{Proof.} We only consider the left prolongations. If $\rmin=-\infty$, then both results are trivial with $\mu^-=\infty$. If, on the other hand, $\rmin>-\infty$, then we can follow the proof of \cite[Props. 4.8 and 4.10]{MPVL} to arrive at the desired conclusion.\qed\\

\subsection{Exponential Decay}

Our task here is to furnish a proof for Proposition \ref{thm:expdecaybewijs1}
and to use this result to establish
Theorem \ref{thm:expdecay}. Our approach
consists of three main steps: constructing
a uniform limit for a sequence
that contradicts (\ref{eq:4.24equiv}),
showing that this limit satisfies
one of the asymptotic systems
(\ref{autlimitsystems}) and subsequently
concluding that this violates the
hyperbolicity assumption (\asref{aannameslimit}{\text{H}}).
The main technical novelties with respect to \cite{MPVL} are contained in the first two steps, where we need to take special care to handle the tail contributions arising from the unbounded shifts.

\begin{lemma}\label{thm:expdecaybewijs1bewijs1} 
Consider the setting of Proposition \ref{thm:expdecaybewijs1} and let $\{\sigma_n\}_{n\geq 1}$, $\{x_n\}_{n\geq 1}$ and $\{\tau_n\}_{n\geq 1}$ be sequences with the following properties.
\begin{enumerate}[label=(\alph*)]
    \item We have $\sigma_n > 0$ 
    for each $n$, together with $\sigma_n\uparrow \infty$.
    \item We either have $x_n\in\mathcal{P}(\tau_n)$ and $\tau_n\leq \tau_-$ for each $n$ or $x_n\in\widehat{\mathcal{P}}(\tau_n)$ and $\tau_n\in\R$ for each $n$.
    
    \item For each $n \ge 1$ we have the bound
    \begin{equation}\label{eq:4.251/2}
    \begin{array}{lclclcl}
         |x_n(-\sigma_n+\tau_n)|&\geq &\frac{1}{2},\end{array}\end{equation}
together with the normalization
\begin{equation}\label{eq:expdecaynormalization}
    \begin{array}{lclclcl}
    \sup\limits_{s\in(-\infty,\tau_n+p]}|x_n(s)|&=&1.
    \end{array}
\end{equation}
\item If $\rmax=\infty$, then we have the additional bound
\begin{equation}\label{eq:4.25}
    \begin{array}{lclclcl}
         |x_n(-\sigma_n+\tau_n)|&\geq & K_{\mathrm{exp}}e^{\alpha(-\sigma_n+\tau_n)}\sup\limits_{s\in[p+\tau_n,\infty)}e^{-\alpha s}|x_n(s)|.\end{array}\end{equation}
\end{enumerate}
Then upon defining $z_n(t)=x_n(t-\sigma_n+\tau_n)$ and passing to a subsequence, we have $z_n\rightarrow z$ uniformly on compact subsets of $\R$. Moreover, we have $z\neq 0$ and $|z|\leq 1$ on $\R$.
\end{lemma}
\textit{Proof.} 
We first
consider the case  $\rmax=\infty$
and treat the two 
possibilities $x_n\in\mathcal{P}(\tau_n)$ and $x_n\in\widehat{\mathcal{P}}(\tau_n)$
simultaneously. In particular, we establish the desired uniform convergence on the compact interval
$I_L=[-L,L]$ for some arbitrary $L \ge 1$, which is contained in $(-\sigma_N, \sigma_N)$ for some sufficiently large $N$.\\

 For $n\geq N$ and $t \in I_L$ we have $|z_n(t)|\leq 1$. In addition, upon writing
 \begin{equation}
 \begin{array}{lcl}
     \overline{A}_{j,n}(t)
      &= &
     A_j(t-\sigma_n+\tau_n)x_n(t-\sigma_n+\tau_n+r_j),
     \\[0.2cm]
     %
     \overline{\mathcal{K}}_n(\xi;t)
     & = & \mathcal{K}(\xi;t-\sigma_n+\tau_n)x_n(t-\sigma_n+\tau_n+\xi),
\end{array}
 \end{equation}
 we obtain
 \begin{equation}
\label{eq:znaccentafschatten}
    \begin{array}{lclcl}
         |\dot{z}_n(t)|&=&|\dot{x}_n(t-\sigma_n+\tau_n)|
          &\le& \sum\limits_{j=-\infty}^\infty \big|\overline{A}_{j,n}(t)\big|
         +\int\limits_{\R}\big|\overline{\mathcal{K}}_n(\xi;t)\big|d\xi .\\[0.2cm]
    \end{array}
\end{equation}
We now split the sum above over the two sets
\begin{equation}
\begin{array}{lclcl}
    J^-_n(t) & = & \{j\in\Z \text{ }|\text{ }r_j\leq p+\sigma_n-t\}
    & \subset & 
    \{j\in\Z\text{ }|\text{ }r_j\leq p\}
    ,
    \\[0.2cm]
    J^+_n(t) & = & \{j\in\Z \text{ }|\text{ }r_j > p+\sigma_n-t\}
    & \subset & 
    \{j\in\Z\text{ }|\text{ }r_j\geq -L+\sigma_n+p\}
    .
\end{array}
\end{equation}

For $j \in J_n^-(t)$ we have
$t-\sigma_n+\tau_n+r_j\leq \tau_n+p$,
which in view of the normalization
(\ref{eq:expdecaynormalization})
allows us to write
\begin{equation}
    \begin{array}{lcl}
         \big| \overline{A}_{j,n}(t)\big|
         &\leq &\nrm{A_j(\cdot)}_\infty.
    \end{array}
\end{equation}
On the other hand, for $j \in J_n^+(t)$ we may use
(\ref{eq:expdecaynormalization})-\eqref{eq:4.25}
to obtain
\begin{equation}
    \begin{array}{lcl}
         \big| \overline{A}_{j,n}(t)\big|
         &\leq &\nrm{A_j(\cdot)}_\infty  K_{\mathrm{exp}}^{-1}e^{\alpha(\sigma_n-\tau_n)}e^{\alpha(t-\sigma_n+\tau_n+r_j)}|x_n(-\sigma_n+\tau_n)|\\[0.2cm]
         &\leq &\nrm{A_j(\cdot)}_\infty  K_{\mathrm{exp}}^{-1} e^{\alpha t}e^{\alpha r_j}\\[0.2cm]
         &\leq &\nrm{A_j(\cdot)}_\infty  K_{\mathrm{exp}}^{-1} e^{\alpha L}e^{\alpha r_j}.
    \end{array}
\end{equation}
In particular, we may use 
(\ref{eq:defp}) to estimate
\begin{equation}\label{eq:znaccentafschatten1}\begin{array}{lcl}
\sum\limits_{j=-\infty}^\infty \big| \overline{A}_{j,n}(t)\big|        &\leq &\sum\limits_{j\in J_n^-(t)}\nrm{A_j(\cdot)}_\infty +\sum\limits_{j\in J_n^+(t)}\nrm{A_j(\cdot)}_\infty  K_{\mathrm{exp}}^{-1}e^{\alpha L}e^{\alpha r_j}
\\[0.4cm]
         &\leq &\sum\limits_{r_j\leq p}\nrm{A_j(\cdot)}_\infty+e^{-2\alpha|L-\sigma_n-p|+\alpha L}\\[0.4cm]
         &= &\sum\limits_{r_j\leq p}\nrm{A_j(\cdot)}_\infty+
         e^{ \alpha(3L - 2p - 2 \sigma_n)}.
    \end{array}
\end{equation}
In a similar fashion, we obtain
the corresponding bound 
\begin{equation}\label{eq:znaccentafschatten2}
    \begin{array}{lcl}
    \int_\R \big| \overline{\mathcal{K}}_{n}(\xi;t)  \big| d \xi       
        &\leq &\sup\limits_{s\in\R}\int\limits_{-\infty}^p|\mathcal{K}(\xi;s)|d\xi+e^{ \alpha(3L - 2p - 2 \sigma_n)}.
    \end{array}
\end{equation}
We hence see that both $\{z_n\}_{n\geq N}$ and $\{\dot{z}_n\}_{n\geq N}$ are uniformly bounded on $I_L$.\\

Using the Ascoli-Arzela theorem, we can now pass over to some subsequence to obtain the convergence $z_n\rightarrow z$ uniformly on compact subsets of $\R$. Moreover, since $z_n(0)\geq \frac{1}{2}$ for each $n$, we obtain $z(0)\geq \frac{1}{2}$ and thus $z\neq 0$. The bound on $z_n(t)$ obtained above implies that also $|z|\leq 1$ on $\R$.\\

If $\rmax<\infty$ then this procedure can be repeated, but now one does not need the second terms
in (\ref{eq:znaccentafschatten1}) and (\ref{eq:znaccentafschatten2}).
In particular, the argument reduces to the one in \cite{MPVL}.
 \qed\\

\begin{lemma}\label{thm:expdecaybewijs1bewijs2}  
Consider the setting of Proposition \ref{thm:expdecaybewijs1} and Lemma \ref{thm:expdecaybewijs1bewijs1}. If the sequence $\{-\sigma_n+\tau_n\}_{n\geq 1}$ is unbounded, then the limiting function $z$ satisfies one of the limiting equations (\ref{autlimitsystems}). If, on the other hand, the sequence $\{-\sigma_n+\tau_n\}_{n\geq 1}$ is bounded, then there exists $\beta\in\R$ in such a way that the function $x(t)=z(t-\beta)$ satisfies $x\in\mathcal{B}$.

\end{lemma}
\textit{Proof.} Without loss of generality we assume that $-\sigma_n+\tau_n\rightarrow\infty$ if the sequence $\{-\sigma_n+\tau_n\}_{n\geq 1}$ is unbounded or $-\sigma_n+\tau_n\rightarrow\beta$ if the sequence $\{-\sigma_n+\tau_n\}_{n\geq 1}$ is bounded.
For convenience, we (re)-introduce the expressions
\begin{equation}
 \begin{array}{lcl}
     \overline{A}_{j,n}(s)
      &= &
     A_j(s-\sigma_n+\tau_n)z_n(s+r_j),
     \\[0.2cm]
     %
     \overline{\mathcal{K}}_n(\xi;s)
     & = & \mathcal{K}(\xi;s-\sigma_n+\tau_n)z_n(s+\xi)
\end{array}
 \end{equation}
 and use the integrated form of (\ref{ditishetprobleem})
to write
\begin{equation}
    \begin{array}{lcl}\label{eq:integralest1}
         z(t_2)-z(t_1)&=&\lim\limits_{n\rightarrow\infty} z_n(t_2)-z_n(t_1)\\[0.2cm]
         &=&\lim\limits_{n\rightarrow\infty}\int\limits_{t_1}^{t_2}\sum\limits_{j=-\infty}^\infty \overline{A}_{j,n}(s) ds
         +\lim\limits_{n\rightarrow\infty}\int\limits_{t_1}^{t_2} \int\limits_{\R}\overline{\mathcal{K}}_{n}(\xi;s)d\xi ds\\[0.2cm]
         &:=& \mathcal{J}_A+ \mathcal{J}_{\mathcal{K}}
    \end{array}
\end{equation}
for an arbitrary pair $t_1 < t_2$ that we fix. Upon introducing the tail expression
\begin{equation}
    \begin{array}{lcl}
        \mathcal{E}_{A;N}
        &=& \lim\limits_{n\rightarrow\infty} \int\limits_{t_1}^{t_2}\sum\limits_{|j|=N+1}^\infty \overline{A}_{j,n}(s) \, ds
    \end{array}
\end{equation}
for any $N \ge 0$,
we readily observe that
\begin{equation}
    \begin{array}{lcl}
         \mathcal{J}_A
         &=&\lim\limits_{n\rightarrow\infty}\int\limits_{t_1}^{t_2}\sum\limits_{j=-N}^N \overline{A}_{j,n}(s) ds
         + \mathcal{E}_{A;N}
         \\[0.2cm]
         &=&\int\limits_{t_1}^{t_2}\sum\limits_{j=-N}^N A_j(\infty)z(s+r_j)ds+
         \mathcal{E}_{A;N}
    \end{array}
\end{equation}
if the sequence $\{-\sigma_n+\tau_n\}_{n\geq 1}$ is unbounded, while
\begin{equation}
    \begin{array}{lcl}
         \mathcal{J}_A
         
         &=&\int\limits_{t_1}^{t_2}\sum\limits_{j=-N}^N A_j(s+\beta)z(s+r_j)ds+
         \mathcal{E}_{A;N}
    \end{array}
\end{equation}
if the sequence $\{-\sigma_n+\tau_n\}_{n\geq 1}$ is bounded. Here we evaluated the limit using
the convergence $-\sigma_n+\tau_n\rightarrow\infty$ or $-\sigma_n+\tau_n\rightarrow\beta$.
Slightly adapting the estimate \eqref{eq:znaccentafschatten1} with $L = \max\{\abs{t_1}, \abs{t_2}\}$, we find
\begin{equation}
\begin{array}{lcl}
|\mathcal{E}_{A;N}|
&\le &(t_2 - t_1) \sum\limits_{\abs{j} > N }
  \nrm{A_j(\cdot)}_{\infty}
  + \lim_{n \to \infty}
      (t_2 - t_1) e^{-2 \alpha \sigma_n} e^{\alpha (3L - 2p)}\\[0.2cm]
      &
      = &(t_2 - t_1) \sum\limits_{\abs{j} > N }
  \nrm{A_j(\cdot)}_{\infty},\end{array}
\end{equation}
which yields $\mathcal{E}_{A;N} \to 0$ as $N \to \infty$.
Since  $|z|\leq 1$ on $\R$, we can now use the dominated convergence theorem to conclude that 
\begin{equation}
    \begin{array}{lcl}
        \mathcal{J}_A &=&
        \int\limits_{t_1}^{t_2}\sum\limits_{j=-\infty}^\infty A_j(\infty)z(s+r_j)ds
    \end{array}
\end{equation}
if the sequence $\{-\sigma_n+\tau_n\}_{n\geq 1}$ is unbounded, while
\begin{equation}
    \begin{array}{lcl}
        \mathcal{J}_A &=&
        \int\limits_{t_1}^{t_2}\sum\limits_{j=-\infty}^\infty A_j(s+\beta)z(s+r_j)ds
    \end{array}
\end{equation}
if the sequence $\{-\sigma_n+\tau_n\}_{n\geq 1}$ is bounded. A similar argument for $\mathcal{J}_{\mathcal{K}}$
hence shows that 
$z$ is a solution of the limiting system (\ref{autlimitsystems}) at $+\infty$.
\qed\\

\textit{Proof of Proposition \ref{thm:expdecaybewijs1}.} Arguing by contradiction, we assume that (\ref{eq:4.24equiv}) or (\ref{eq:4.24equivalt}) fails. We can then construct sequences $\{\sigma_n\}_{n\geq 1}$, $\{x_n\}_{n\geq 1}$ and $\{\tau_n\}_{n\geq 1}$ that satisfy properties (i)-(iv) of Lemma \ref{thm:expdecaybewijs1bewijs1}. If the sequence $\{-\sigma_n+\tau_n\}_{n\geq 1}$ is also unbounded, then Lemma \ref{thm:expdecaybewijs1bewijs2} yields that $z$ is a non-trivial, bounded solution of one of the limiting equations (\ref{autlimitsystems}), contradicting the hyperbolicity of these systems.\\

If on the other hand the sequence $\{-\sigma_n+\tau_n\}_{n\geq 1}$ is bounded, we can assume that $-\sigma_n+\tau_n\rightarrow\beta$ for some $\beta\in\R$. Since necessarily $\tau_n\rightarrow \infty$, this can only happen if $x_n\in\widehat{\mathcal{P}}(\tau_n)$ for each $n$. Lemma \ref{thm:expdecaybewijs1bewijs2} yields that $x_n\rightarrow x$ uniformly on compact subsets of $\R$ and that $0\neq x\in\mathcal{B}$. On account of Proposition \ref{prop:fredholmproperties} we find that $x$ decays exponentially. By definition of $\widehat{\mathcal{P}}$ we therefore obtain
\begin{equation}\label{eq:4.26equiv}
    \begin{array}{lclcl}
         0&=&\int\limits_{-\infty}^{\infty}x(t)^\dagger x_n(t)dt
         &\rightarrow & \int\limits_{-\infty}^\infty |x(t)|^2dt,
    \end{array}
\end{equation}
which yields a contradiction since $x\neq 0$.\qed\\

We now shift our attention to the proof of Theorem \ref{thm:expdecay}. In particular, we set up an iteration scheme
to leverage the bound
\eqref{eq:4.24equiv} and show that
solutions in $\mathcal{P}(\tau)$ decay exponentially. As a preparation,
we provide a uniform bound on the supremum of such solutions.

\begin{lemma}
\label{thm:expdecaybewijs2} Assume that (\asref{aannamesconstanten}{\text{H}}), (\asref{aannamesconvoluties}{\text{H}}) and (\asref{aannameslimit}{\text{H}}) are satisfied. Recall the constant $\mu^-$ from Lemma \ref{prop:uniqueleftprol} and fix $\tau^-<\mu^-$. Then there exists $C>0$ in such a way for each $\tau\leq \tau^-$ and each $x\in\mathcal{P}(\tau)$ we have the bound
\begin{equation}\label{eq:4.25equiv}
    \begin{array}{lcl}
         \nrm{x}_{C_b(D_\tau^\ominus)}&\leq & C\nrm{x_{\tau}}_\infty.
    \end{array}
\end{equation}
The same bound holds for any $x\in \widehat{\mathcal{P}}(\tau)$, with a possibly different value of $C$, where now any $\tau\in\R$ is permitted.
\end{lemma}
\textit{Proof.} The bound (\ref{eq:4.25equiv}) is in fact an equality with $C=1$ if $\rmin=-\infty$, so we assume that $\rmin>-\infty$. If $\rmax<\infty$ the final part of the proof of \cite[Thm. 4.2]{MPVL} can be repeated, hence we also assume that $\rmax=\infty$.\\

Arguing by contradiction, we consider sequences $\{x_n\}_{n\geq 1},\{\tau_n\}_{n\geq 1}$ and
$\{C_n\}_{n\geq 1}$ with $C_n\rightarrow\infty$ and
\begin{equation}
    \begin{array}{lclcl}
         \nrm{x_n}_{C_b(D_{\tau_n}^\ominus)}&=&C_n\nrm{(x_n)_{\tau_n}}_\infty&=&1,
    \end{array}
\end{equation}
with either $x_n\in \mathcal{P}(\tau_n)$ and $\tau_n\leq \tau^-$ for each $n$ or $x_n\in\widehat{\mathcal{P}}(\tau_n)$ and $\tau_n\in\R$.\\

We want to emphasize that due to the lack of a natural choice for the sequence $\{\sigma_n\}_{n\geq 1}$ which satisfies (a) of Lemma \ref{thm:expdecaybewijs1bewijs1}, we cannot immediately apply this result. However, we will follow more or less the same procedure to arrive at a slightly weaker conclusion. Note that the function $z_n(t)=x_n(t+\tau_n)$ is a solution of (\ref{ditishetprobleem}) on the interval $(-\infty,0]$ for each value of $n$. In addition, we note that
\begin{equation}\label{eq:znon0infty}
    \begin{array}{lclcl}
         \sup\limits_{t\in[\rmin,\infty)}|z_n(t)|&\leq & \nrm{(x_n)_{\tau_n}}_\infty
         &=&
          C_n^{-1}.
    \end{array}
\end{equation}

We can now follow the proof of Lemma \ref{thm:expdecaybewijs1bewijs1}, using (\ref{eq:znon0infty}) to control the behaviour of $z_n$ on $[0,\infty)$, and pass to a subsequence to obtain $z_n\rightarrow z$ uniformly on compact subsets of $(-\infty,0]$. In addition, 
(\ref{eq:znon0infty})
allows us to extend this convergence 
to all compact subsets of $\R$, with $z_0=0$. 
For each $n\geq 1$ we pick $s_n$ in such a way that  $|x_n(-s_n+\tau_n)|=1$.
On account of Proposition \ref{thm:expdecaybewijs1} the set 
$\{s_n\}_{n\geq 1}$ is bounded,
which means that 
$z$ is not identically zero.\\

Suppose first that the sequence $\{\tau_n\}_{n\geq 1}$ is unbounded. Since each function $z_n$ is a solution of (\ref{ditishetprobleem}) on $(-\infty,0]$,
we can follow the proof of Lemma \ref{thm:expdecaybewijs1bewijs2} to conclude that $z$ is a bounded solution of one of the limiting equations (\ref{autlimitsystems}) on $(-\infty,0]$. Moreover, since $z_0=0$ it follows that $z$ is also a solution of the limiting equation (\ref{autlimitsystems}) on $[0,\infty)$. Hence $z$ is a nontrivial, bounded solution on $\R$ of one of the limiting equations (\ref{autlimitsystems}), which yields a contradiction. \\

Suppose now that $\{\tau_n\}_{n\geq 1}$ is in fact a bounded sequence. Then 
after passing to a subsequence
we obtain  $\tau_n\rightarrow \tau_0$.
Following the proof of Lemma \ref{thm:expdecaybewijs1bewijs2},
we see that the function
$x(t)=z(t-\tau_0)$
is a nontrivial, bounded solution of (\ref{ditishetprobleem}) on $(-\infty,\tau_0]$. Since $z_0=0$, we get that $x_{\tau_0}=0$ and therefore $x$ is a nontrivial, bounded left prolongation of the zero solution from the starting point $\tau_0$. If $\tau_0<\mu^-$, this gives an immediate contradiction to Lemma \ref{prop:uniqueleftprol}. If on the other hand $\tau_0\geq\mu^->\tau^-$, then
our assumptions allow us to conclude that
$x_n\in\widehat{\mathcal{P}}(\tau_n)$
for all $n$.
A  computation similar to (\ref{eq:4.26equiv}) shows that $x\in\widehat{\mathcal{P}}(\tau_0)$, which contradicts Lemma \ref{prop:uniqueleftprol}. This establishes (\ref{eq:4.25equiv}).\qed\\

\begin{lemma}\label{thm:expdecaybewijs3}
Assume that (\asref{aannamesconstanten}{\text{H}}), (\asref{aannamesconvoluties}{\text{H}}) and (\asref{aannameslimit}{\text{H}}) are satisfied. Recall the constant $\mu^-$ from Lemma \ref{prop:uniqueleftprol} and fix $\tau^-<\mu^-$. Then there exist constants $\tilde{K}>0$ and $\tilde{\alpha}>0$ so that the bound
\begin{equation}
    \begin{array}{lcl}
         |x(t)|&\leq & \tilde{K}e^{\alpha(t-\tau)}\nrm{x}_{C_b(D_\tau^\ominus)}
    \end{array}
\end{equation}
holds for all $\tau\leq\tau^-$, all $x\in\mathcal{P}(\tau)$ and all $t\leq \tau+p$.
\end{lemma}
\textit{Proof.}
The proof of \cite[Thm. 4.2]{MPVL}
can be used to handle the case $\rmax<\infty$, so we assume here
that $\rmax=\infty$. Pick any $x\in\mathcal{P}(\tau)$, which we normalize to have $\nrm{x}_{C_b(D_\tau^\ominus)}=1$. Recalling the constants 
from Proposition \ref{thm:expdecaybewijs1},
we assume without loss of generality that 
\begin{equation}
\label{eq:exp:dec:prlm:bnd:k:exp}
    \begin{array}{lclcl}
        K_{\mathrm{exp}}&\geq &1, 
        \qquad \qquad
        K_{\mathrm{exp}}e^{-\alpha(\sigma+p)}&\leq &\frac{1}{4}.
    \end{array}
\end{equation}
For $t\leq -\sigma+\tau$, 
this allows us to estimate
\begin{equation}\label{eq:m=0}
    \begin{array}{lcl}
         |x(t)|&\leq &\max\Big\{\frac{1}{2}\sup\limits_{s\in(-\infty,\tau+p]}|x(s)|, K_{\mathrm{exp}}\sup\limits_{s\in[p+\tau,\infty)}e^{-\alpha(s-t)}|x(s)|\Big\}\\[0.4cm]
         &\leq &\max\big\{\frac{1}{2}, K_{\mathrm{exp}}e^{-\alpha(p+\tau+\sigma-\tau)}\big\}\\[0.2cm]
         &=&
         \frac{1}{2}.
    \end{array}
\end{equation}
We aim to show, by induction, that for each integer $m \ge 0$ we have the bound
\begin{equation}\label{eq:indhyp}
    \begin{array}{lclcl}
    
         |x(t)|&\leq &2^{-(m+1)},\qquad\qquad t&\leq &t_m,
    \end{array}
\end{equation}
where we have introduced
\begin{equation}
    \begin{array}{lcl}
         t_m&:=&-m(\sigma+p)-\sigma+\tau.
    \end{array}
\end{equation}
Indeed, if (\ref{eq:indhyp}) holds for each $m\in\Z_{\geq 0}$, then we obtain the desired estimate
\begin{equation}\label{eq:indunif}
    \begin{array}{lcl}
         |x(t)|&\leq &\tilde{K}e^{\tilde{\alpha}(t-\tau)} 
    \end{array}
\end{equation}
for any $t\leq\tau$ with $\tilde{\alpha}=\frac{\ln(2)}{\sigma+p}$ and $\tilde{K}=e^{\tilde{\alpha}(\sigma+p)}$, which concludes the proof.  \\

The case $m=0$ follows from (\ref{eq:m=0}), so we pick $M\geq 1$ and assume that (\ref{eq:indhyp}) holds for each value of $0\leq m\leq M-1$.
Since $x\in\mathcal{P}(\tau)$ and since $\sigma>0$ and $p>0$, we must have $x\in\mathcal{P}(t_M+\sigma)$ as well. Fix $t\leq t_M$. Then Proposition \ref{thm:expdecaybewijs1} yields the bound
\begin{equation}\label{eq:indbewijs0}
    \begin{array}{lcl}
         |x(t)|&\leq &\max\Big\{\frac{1}{2}\sup\limits_{s\in(-\infty,t_M+\sigma+p]}|x(s)|, K_{\mathrm{exp}}\sup\limits_{s\in[t_M+\sigma+p,\infty)}e^{-\alpha(s-t)}|x(s)|\Big\}.
    \end{array}
\end{equation}
Since $t_M + \sigma + p = t_{M-1}$,
we may apply 
(\ref{eq:indhyp}) with $m=M-1$
to obtain
\begin{equation}\label{eq:indbewijs1}
    \begin{array}{lclcl}
         \frac{1}{2}\sup\limits_{s\in(-\infty,t_M+\sigma+p]}|x(s)|
         &\le &\frac{1}{2} 2^{-M}
         &= &2^{-(M+1)}.
    \end{array}
\end{equation}
In addition, 
we may use \eqref{eq:exp:dec:prlm:bnd:k:exp} and (\ref{eq:indhyp})
to estimate
\begin{equation}\label{eq:indbewijs3}
    \begin{array}{lcl}
          K_{\mathrm{exp}}\sup\limits_{s\in[t_m,t_{m-1}]}e^{-\alpha(s-t)}|x(s)|&\leq & K_{\mathrm{exp}} e^{-\alpha\big(t_m-t_M\big)}2^{-m}\\[0.2cm]
         &= & K_{\mathrm{exp}}e^{-\alpha(M-m)(p+\sigma)}2^{-m}\\[0.2cm]
         &\leq &\big(\frac{1}{4}\big)^{M-m}2^{-m}\\[0.2cm]
         &\leq &2^{-(M+1)},
    \end{array}
\end{equation}
for $0\leq m\leq M-1$. Finally, we can estimate
\begin{equation}\label{eq:indbewijs4}
    \begin{array}{lcl}
          K_{\mathrm{exp}}\sup\limits_{s\in[\tau-\sigma,\infty)}e^{-\alpha(s-t)}|x(s)|&\leq & K_{\mathrm{exp}} e^{-\alpha(\tau-\sigma-t_M)}\\[0.2cm]
         &= & K_{\mathrm{exp}}e^{-\alpha M(p+\sigma)}\\[0.2cm]
         &\leq &2^{-(M+1)}.
    \end{array}
\end{equation}
Combining (\ref{eq:indbewijs0}) with (\ref{eq:indbewijs1})-(\ref{eq:indbewijs4}) now yields the bound
\begin{equation}
    \begin{array}{lcl}
         |x(t)|&\leq &2^{-(M+1)},
    \end{array}
\end{equation}
as desired. 
\qed\\

\textit{Proof of Theorem \ref{thm:expdecay}.} 
We only show the result for the $\mathcal{P}$-spaces; the result for the $\mathcal{Q}$-spaces follows analogously. If $\rmax<\infty$, the proof of \cite[Thm. 4.2]{MPVL} can be repeated, so we assume that $\rmax=\infty$. Pick any $x\in\mathcal{P}(\tau)$. From Lemma \ref{thm:expdecaybewijs3} we obtain the bound
\begin{equation}
    \begin{array}{lcl}
         |x(t)|&\leq & \tilde{K}e^{\alpha(t-\tau)}\nrm{x}_{C_b(D_\tau^\ominus)}
    \end{array}
\end{equation}
for all $t\leq \tau+p$. Since $x\in\mathcal{P}(\tau)$, we can write
\begin{equation}
    \begin{array}{lcl}
         \dot{x}(t)&=&\sum\limits_{j=-\infty}^\infty A_j(t)x(t+r_j)+\int_{\R}\mathcal{K}(\xi;t)x(t+\xi)d\xi
    \end{array}
\end{equation}
for $t\leq\tau$. Lemma \ref{thm:expdecaybewijs0} allows us to estimate 

\begin{equation}
    \begin{array}{lcl}
         \big|\sum\limits_{j=-\infty}^\infty A_j(t)x(t+r_j)d\xi\big|
         &\leq &\sum\limits_{t+r_j\leq \tau+p}\nrm{A_j(\cdot)}_\infty \tilde{K}e^{\alpha(t+r_j-\tau)}\nrm{x}_\infty+\sum\limits_{t+r_j>\tau+p}\nrm{A_j(\cdot)}_\infty \nrm{x}_\infty\\[0.4cm]
         
         &\leq &\sum\limits_{t+r_j\leq \tau}\nrm{A_j(\cdot)}_\infty e^{\alpha|r_j|} \tilde{K}e^{\alpha(t-\tau)}\nrm{x}_\infty+ K_{\mathrm{exp}}e^{2\alpha(t-\tau)}\nrm{x}_\infty\\[0.4cm]
         &\leq &\sum\limits_{j=-\infty}^\infty \nrm{A_j(\cdot)}_\infty e^{\alpha|r_j|} \tilde{K}e^{\alpha(t-\tau)}\nrm{x}_\infty+ K_{\mathrm{exp}}e^{2\alpha(t-\tau)}\nrm{x}_\infty
    \end{array}
\end{equation}
for any $t\leq \tau$. Using a similar estimate for the convolution kernel, we obtain the bound
\begin{equation}\label{eq:finalexp1}
    \begin{array}{lcl}
        |\dot{x}(t)|&\leq &\tilde{K}e^{\alpha(t-\tau)}\nrm{x}_\infty\sum\limits_{j=-\infty}^\infty \nrm{A_j(\cdot)}_\infty e^{\alpha|r_j|}+ K_{\mathrm{exp}}e^{2\alpha(t-\tau)}\nrm{x}_\infty\\[0.2cm]
         &&\qquad +\tilde{K}e^{\alpha(t-\tau)}\nrm{x}_\infty\sup\limits_{s\in\R}\nrm{\mathcal{K}(\cdot;s)}_{\tilde{\eta}}+ K_{\mathrm{exp}}  e^{2\alpha(t-\tau)}\nrm{x}_\infty       
    \end{array}
\end{equation}
for any $t\leq \tau$. Since $\rmax=\infty$, we can derive from Lemma \ref{thm:expdecaybewijs2} that
\begin{equation}\label{eq:finalexp2}
    \begin{array}{lclcl}
         \nrm{x}_\infty&=&\nrm{x}_{C_b(\D_\tau^\ominus)}&\leq & C\nrm{x_\tau}_\infty.
    \end{array}
\end{equation}
The bounds (\ref{eq:finalexp1})-(\ref{eq:finalexp2}) together establish the desired result.
\qed\\

\subsection{The restriction operators \texorpdfstring{$\pi^+$}{pi+} and \texorpdfstring{$\pi^-$}{pi-}}
\label{sec:ex:restr}

It is often convenient to split the domain $D_X$ into the two parts
\begin{equation}
    \begin{array}{lclcl}
        D_X^- &=&\overline{(\rmin,0)},
        \qquad \qquad
        D_X^+&=&\overline{(0,\rmax)}
    \end{array}
\end{equation}
and study the restriction of functions in $X$ to the spaces
\begin{equation}
    \begin{array}{lclcl}
         X^- &=& C_b(D_X^-),
        \qquad \qquad
        X^+ &=& C_b(D_X^+).
    \end{array}
\end{equation}
In particular, we introduce the 
operators $\pi^+:X\rightarrow X^+$ and $\pi^-:X\rightarrow X^-$ that act as
\begin{equation}\label{eq:defpi}
    \begin{array}{lclcl}
         (\pi^\pm f)(t)&=&f(t),\qquad \qquad &t\in D_X^\pm.\\[0.2cm]
    \end{array}
\end{equation}
Moreover, for a subspace $E\subset X$ we let $\pi^+_E$ and $\pi^-_E$ denote the restrictions of $\pi^+$ and $\pi^-$ to $E$. 
We obtain some preliminary compactness results below, leaving a more detailed analysis
of these operators to {\S}\ref{section:fredholm}.

\begin{proposition}[{cf. \cite[Thm. 4.4]{MPVL}}]\label{compactoperator} Assume that (\asref{aannamesconstanten}{\text{H}}), (\asref{aannamesconvoluties}{\text{H}}) and (\asref{aannameslimit}{\text{H}}) are satisfied. Then
for every $\tau\in\R$,
the operators $\pi_{P(\tau)}^-$, $\pi_{Q(\tau)}^+$, 
$\pi_{\widehat{P}(\tau)}^-$ and $\pi_{\widehat{Q}(\tau)}^+$ are all compact.
\end{proposition}
\textit{Proof.} Suppose first that $\rmin=-\infty$ and fix $\tau\in\R$. Let $\{\phi_n\}_{n\geq 1}$ be a bounded sequence in $\widehat{P}(\tau)$ and write $\{x_n\}_{n\geq 1}$ 
for the corresponding sequence in $\widehat{\mathcal{P}}(\tau)$ that has $(x_n)_\tau=\phi_n$ for each $n\geq 1$. After passing to a subsequence,
the exponential bound (\ref{eq:4.8equiv:p}) allows us to obtain the convergence
$x_{n}\rightarrow x$ uniformly on compact subsets of $(-\infty,0]$.
For any $\epsilon > 0$, 
we can use (\ref{eq:4.8equiv:p}) to pick $L \gg 1$ in such a way that
$|x_n(t)| < \frac{\epsilon}{2}$
and hence $|x(t)| \le \frac{\epsilon}{2}$ holds for all $t \le - L$. The uniform convergence on $[-L, 0]$ now allows us to pick $N \gg 1$ so that $|x_n(t) - x(t)| \le \epsilon$ for all $t \le 0$ and $n \ge N$. In particular,  $\{x_{n}\}_{n\geq 1}$ converges in $X^-$, which shows that $\pi^-_{\widehat{P}(\tau)}$ is compact.
\\

The case where $\rmin>-\infty$ 
can be treated as in the proof of \cite[Thm. 4.4]{MPVL} and will be omitted. The compactness of $\pi^+_{\widehat{Q}(\tau)}$
follows by symmetry. Finally,
the operators $\pi_{P(\tau)}^-$ and $\pi_{Q(\tau)}^+$ are compact since they are finite-dimensional extensions of $\pi_{\widehat{P}(\tau)}^-$ and $\pi_{\widehat{Q}(\tau)}^+$ respectively.\qed\\

The second part of Corollary \ref{cor:4.11equiv} below references the subpaces $P(\pm\infty)\subset X$ and $Q(\pm\infty)\subset X$, being the spaces corresponding the limiting equations (\ref{autlimitsystems}) with the decomposition given in (\ref{eq:4.10}). Since the systems (\ref{autlimitsystems}) also satisfy the conditions (\asref{aannamesconstanten}{\text{H}}), (\asref{aannamesconvoluties}{\text{H}}) and (\asref{aannameslimit}{\text{H}}), we can apply the results from the previous sections to the subspaces $P(\pm\infty)$ and $Q(\pm\infty)$.

\begin{corollary}[{cf. \cite[Cor. 4.11]{MPVL}}]\label{cor:4.11equiv} 
Assume that (\asref{aannamesconstanten}{\text{H}}), (\asref{aannamesconvoluties}{\text{H}}) and (\asref{aannameslimit}{\text{H}}) are satisfied
and let $\{\phi_n\}_{n\geq 1}$ and $\{\psi_n\}_{n \ge 1}$ be bounded sequences, with $\phi_n\in \widehat{P}(\tau_n)$ and $\psi_n\in \widehat{P}(\tau_0)$ for each $n\geq 1$. 
Suppose furthermore that  $\tau_n\rightarrow\tau_0$ and that the sequence $\{\pi^+(\phi_n-\psi_n)\}_{n\geq 1}$ converges in $X^+$. Then after passing to a subsequence,
the differences $\{\phi_n-\psi_n\}_{n\geq 1}$ converge to some $\phi\in \widehat{P}(\tau_0)$, uniformly on compact subsets of $D_X$.\\

The conclusion above remains valid after
making the replacements
\begin{equation}
\{ \widehat{P}(\tau_n), \widehat{P}(\tau_0), \tau_0 \}
\mapsto \{ P(\tau_n), P(-\infty), -\infty \}.
\end{equation}
In addition,
the analogous results hold for the spaces $\widehat{Q}$ and $Q$
after replacing $\pi^+$ by $\pi^-$
and $-\infty$ by $+\infty$.
\end{corollary}
\textit{Proof.} 
For each $n\geq 1$ we let $y_n\in\widehat{\mathcal{P}}(\tau_n)$ and $z_n\in\widehat{\mathcal{P}}(\tau_0)$ denote the left prolongations of $\phi_n$ and $\psi_n$ respectively. Moreover, we write $x_n(t)=y_n(t+\tau_n-\tau_0)-z_n(t)$ for $t\leq \tau_0+\rmax$. Then $x_n$ satisfies the inhomogeneous version of (\ref{ditishetprobleem}) given by
\begin{equation}\label{eq:inhomg}
    \begin{array}{lcl}
         \dot{x}_n(t)&=&\sum\limits_{j=-\infty}^\infty A_j(t)x_n(t+r_j)+\int\limits_{\R}\mathcal{K}(\xi;t)x_n(t+\xi)d\xi+h_n(t),
    \end{array}
\end{equation}
in which $h_n$ is defined by
\begin{equation}
    \begin{array}{lcl}
         h_n(t)&=&\sum\limits_{j=-\infty}^\infty \big(A_j(t+\tau_n-\tau_0)-A_j(t)\big)y_n(t+r_j+\tau_n-\tau_0)\\[0.2cm]
         &&\qquad +\int\limits_{\R}\big(\mathcal{K}(\xi;t+\tau_n-\tau_0)-\mathcal{K}(\xi;t)\big)y_n(t+\xi+\tau_n-\tau_0)d\xi.
    \end{array}
\end{equation}
Since $x_n$ satisfies the inhomogeneous equation (\ref{eq:inhomg}), since $\sum\limits_{|j|=N}^\infty \nrm{A_j}_\infty\rightarrow 0$ as $N\rightarrow\infty$, since $\sup\limits_{t\in\R}\nrm{\mathcal{K}(\cdot;t)}_{\tilde{\eta}}<\infty$, and since both $y_n$ and $z_n$ enjoy the uniform exponential estimates in Theorem \ref{thm:expdecay}, we see that the sequence $\{x_n\}_{n\geq 1}$ is uniformly bounded and equicontinuous. Hence we can apply the Ascoli-Arzela theorem to pass over to a subsequence for which $x_n\rightarrow x$ uniformly on compact subsets of $(-\infty,\tau_0]$. Moreover, $x$ is bounded and 
the convergence
$x_n\rightarrow x$ is uniform on $D_X^++\tau_0$
since $\{\pi^+(\phi_n-\psi_n)\}_{n\geq 1}$ converges in $X^+$.
However, in contrast to \cite{MPVL} we cannot conclude that this convergence is uniform on $D_X$, since this interval is not necessarily compact.\\

We see that $h_n\rightarrow 0$ in $L^1(I)$ for any bounded interval $I\subset (-\infty,\tau_0]$, again using the limit $\sum\limits_{|j|=N}^\infty \nrm{A_j}_\infty\rightarrow 0$ as $N\rightarrow\infty$, the bound $\sup\limits_{t\in\R}\nrm{\mathcal{K}(\cdot;t)}_{\tilde{\eta}}<\infty$ and the fact that the sequence $\{y_n\}_{n\geq 1}$ is bounded uniformly on $D_0^\ominus$. Similarly to the proof of Lemma \ref{thm:expdecaybewijs1bewijs2}, we obtain that $x:D_{\tau_0}^\ominus\rightarrow \C^M$ is a bounded solution of (\ref{ditishetprobleem}) on $(-\infty,\tau_0]$,
which yields
$x\in\mathcal{P}(\tau_0)$. Finally, for every $w\in\mathcal{B}$ we obtain
\begin{equation}\label{eq:4.27equiv}
    \begin{array}{lcl}
         0&=& \int\limits_{-\infty}^{\tau_n+\rmax}w(t)^\dagger y_n(t)dt-\int\limits_{-\infty}^{\tau_0+\rmax}w(t)^\dagger z_n(t)dt\\[0.2cm]
         &=&\int\limits_{-\infty}^{\tau_0+\rmax}w(t)^\dagger x_n(t)dt-\int\limits_{-\infty}^{\tau_n+\rmax}\big(w(t)-w(t-\tau_n-\tau_0)\big)^\dagger y_n(t)dt\\[0.2cm]
         &\rightarrow &\int\limits_{-\infty}^{\tau_0+\rmax}w(t)^\dagger x(t)dt,
    \end{array}
\end{equation}
since $w$ decays exponentially on account of Proposition \ref{prop:fredholmproperties}. Therefore we must have $x\in\widehat{\mathcal{P}}(\tau_0)$ and thus $\phi:=x_{\tau_0}\in\widehat{P}(\tau_0)$.\\

The result for $P(\tau_n)$ where $\tau_n\rightarrow-\infty$ follows a similar proof. We now use the estimate (\ref{eq:4.8equiv:p}), which is valid for sufficiently small $\tau$. Naturally, the integral computation (\ref{eq:4.27equiv}) is not needed in this proof.
The remaining results follow by symmetry. \qed\\

\subsection{Fundamental properties of the Hale inner product}\label{subsec:fundhale}
We now shift our focus towards the Hale inner product, which plays an important role throughout the remainder of the paper.
In particular, we establish the identity (\ref{eq:fundHale}), which requires
special care on account of the infinite sums. In addition, we study the limiting behaviour of the Hale inner product and establish a uniform estimate
that holds
for exponentially decaying functions.

\begin{lemma}\label{lemma:diffhale} Assume that (\asref{aannamesconstanten}{\text{H}}), (\asref{aannamesconvoluties}{\text{H}}) and (\asref{aannameslimit}{\text{H}}) are satisfied and fix two functions $x,y\in C_b(\R)$. 
Suppose furthermore that $x$ and $y$ are both differentiable at some time $t \in \R$. Then we have the identity
\begin{equation}\label{eq:diffhale}
    \begin{array}{lcl}
         \frac{d}{dt}\ip{y^t,x_t}_t&=&y^\dagger (t)[\Lambda x](t)+[\Lambda^* y](t)^\dagger x(t).
    \end{array}
\end{equation}
In particular, if either $x\in\mathcal{P}(\tau)$ and $y\in\mathcal{B}^*$ or if $x\in\mathcal{B}$ and $y\in\mathcal{Q}(\tau)$ for some $\tau\in\R$, then $\frac{d}{dt}\ip{y^t,x_t}_t=0$ for all $t\leq \tau$ or all $t\geq\tau$ respectively.
\end{lemma}
\textit{Proof.} For any $t\in\R$ we can rewrite the Hale inner product
in the form
\begin{equation}\label{eq:diffhale1}
    \begin{array}{lcl}
         \ip{y^t,x_t}_t&=&y^t(0)^\dagger x_t(0)-\sum\limits_{j=-\infty}^\infty \int\limits_0^{r_j}y^t(s-r_j)^\dagger A_j(t+s-r_j)x_t(s)ds\\[0.2cm]
         &&\qquad -\int\limits_{\R}\int\limits_0^r y^t(s-r)^\dagger \mathcal{K}(r;t+s-r)x_t(s)dsdr\\[0.2cm]
         &=&y(t)^\dagger x(t)-\sum\limits_{j=-\infty}^\infty \int\limits_t^{t+r_j}y(s-r_j)^\dagger A_j(s-r_j)x(s)ds\\[0.2cm]
         &&\qquad -\int\limits_{\R}\int\limits_t^{t+r} y(s-r)^\dagger \mathcal{K}(r;s-r)x(s)dsdr.
    \end{array}
\end{equation}
We aim to compute the derivative $\frac{d}{dt}\ip{y^t,x_t}_t$, so the main difficulty compared to \cite{MPVL} is that we need to interchange a derivative and an infinite sum as well as a derivative and an integral instead of a derivative and a finite sum. Since we can estimate
\begin{equation}\label{eq:diffhale2}
    \begin{array}{lcl}
        \sum\limits_{j=-\infty}^\infty  \Big| \frac{d}{dt}\int\limits_t^{t+r_j}y(s-r_j)^\dagger A_j(s-r_j)x(s)ds\Big|
         &=&\sum\limits_{j=-\infty}^\infty \Big|y(t)^\dagger A_j(t)x(t+r_j)-y(t-r_j)^\dagger A_j(t-r_j)x(t)\Big|\\[0.2cm]
         &\leq &2\nrm{x}_\infty\nrm{y}_\infty\sum\limits_{j=-\infty}^\infty \nrm{A_j(\cdot)}_\infty,
    \end{array}
\end{equation}
we see that this series converges uniformly. In a similar fashion
we can estimate
\begin{equation}\label{eq:diffhale2.5}
\begin{array}{lcl}
\int\limits_{\R}\Big|\frac{d}{dt}\int\limits_t^{t+r} y(s-r)^\dagger \mathcal{K}(r;s-r)x(s)ds\Big|dr&= &\int\limits_{\R}\big| y(t)^\dagger \mathcal{K}(r;t)x(t+r)-y(t-r)^\dagger \mathcal{K}(r;t-r)x(t)\big|dr\\[0.2cm]
&\leq &\nrm{x}_\infty\nrm{y}_\infty \big[\sup\limits_{t\in\R}\nrm{\mathcal{K}(\cdot;t)}_{\tilde{\eta}}+\sup\limits_{t\in\R}\nrm{\mathcal{K}(\cdot;t-\cdot)}_{\tilde{\eta}}\big].
\end{array}
\end{equation}
We can hence freely exchange
a time derivative with the integral
and sum in \eqref{eq:diffhale2}
to obtain
\begin{equation}\label{eq:diffhale3}
    \begin{array}{lcl}
         \frac{d}{dt}\ip{y^t,x_t}_t&=&\dot{y}(t)^\dagger x(t)+y(t)^\dagger \dot{x}(t)-\sum\limits_{j=-\infty}^\infty \big[y(t)^\dagger A_j(t)x(t+r_j)-y(t-r_j)^\dagger A_j(t-r_j)x(t)\big]\\[0.2cm]
         &&\qquad -\big[\int\limits_{\R} y(t)^\dagger \mathcal{K}(r;t)x(t+r)dr-\int\limits_{\R}y(t-r)^\dagger \mathcal{K}(r;t-r)x(t)dr\big]\\[0.2cm]
         &=&y^\dagger (t)[\Lambda x](t)+[\Lambda^* y](t)^\dagger x(t).
    \end{array}
\end{equation}
The final statement follows trivially from (\ref{eq:diffhale}).\qed\\ 

\begin{lemma}\label{lemma:diffhale2} Assume that (\asref{aannamesconstanten}{\text{H}}), (\asref{aannamesconvoluties}{\text{H}}) and (\asref{aannameslimit}{\text{H}}) are satisfied and fix two functions $x,y\in C_b(\R)$. 
Suppose furthermore that $y(t)$ decays exponentially as $t\rightarrow\infty$. Then we have the limit
\begin{equation}
    \begin{array}{lcl}
         \lim\limits_{t\rightarrow\infty}\ip{y^t,x_t}_t&=&0.
    \end{array}
\end{equation}
The corresponding estimate holds for $t\rightarrow-\infty$ if $y(t)$ decays exponentially as $t\rightarrow -\infty$.
\end{lemma}
\textit{Proof.} Pick $0<\beta<\tilde{\eta}$ and $K>0$ in such a way that $|y(t)|\leq K e^{-\beta t}$ for $t\geq 0$. Upon choosing a small $\epsilon>0$, we first pick $N\in\Z_{\geq 1}$ in such a way that the bound

\begin{equation}
    \begin{array}{lcl}
         \sum\limits_{|j|=N+1}^\infty |r_jA_j(s-r_j)|\nrm{x}_\infty\nrm{y}_\infty+\int\limits_{(-\infty,-N]\cup [N,\infty)} |r\mathcal{K}(r;s-r)|\nrm{x}_\infty\nrm{y}_\infty dr&\leq &\frac{\epsilon}{6}
    \end{array}
\end{equation}
holds for all $s\in\R$. 
We pick $T > N$ in such a way that
also $T > \max \{|r_j|:-N\leq j\leq N\}$ and that we have the estimates
\begin{equation}
    \begin{array}{lcl}
         |y(t)|\nrm{x}_\infty&\leq &\frac{\epsilon}{3},\\[0.2cm]
         
         Ke^{-\beta t}\sum\limits_{j=-N}^N e^{\beta |r_j|}|r_j A_j(s)|\nrm{x}_\infty&\leq &\frac{\epsilon}{6},\\[0.2cm]
         Ke^{-\beta t}\int_{-N}^N e^{\beta |r|}|r\mathcal{K}(r;s-r)|\nrm{x}_\infty dr&\leq &\frac{\epsilon}{6}
    \end{array}
\end{equation}
for all $t\geq T$ and all $s\in\R$. In particular, we can estimate 
\begin{equation}
    \begin{array}{lcl}
         \sum\limits_{j=-\infty}^\infty\big| \int\limits_t^{t+r_j}y(s-r_j)^\dagger A_j(s-r_j)x(s)ds\big|&=&\sum\limits_{j=-N}^N \big|\int\limits_t^{t+r_j}y(s-r_j)^\dagger A_j(s-r_j)x(s)ds\big| \\[0.2cm]
         &&\qquad +\sum\limits_{|j|=N+1}^\infty \big|\int\limits_t^{t+r_j}y(s-r_j)^\dagger A_j(s-r_j)x(s)ds\big|\\[0.2cm]
         &\leq &\sup\limits_{s\in\R}\sum\limits_{j=-N}^NK\max\{e^{-\beta t},e^{-\beta(t-r_j)}\}|r_jA_j(s-r_j)|\nrm{x}_\infty\\[0.2cm]
         &&\qquad +\sup\limits_{s\in\R}\sum\limits_{|j|=N+1}^\infty |r_jA_j(s-r_j)|\nrm{x}_\infty\nrm{y}_\infty\\[0.2cm]
         &\leq &Ke^{-\beta t}\sup\limits_{s\in\R}\sum\limits_{j=-N}^Ne^{\beta|r_j|}|r_jA_j(s-r_j)|\nrm{x}_\infty+\frac{\epsilon}{4}\\[0.2cm]
         &\leq &\frac{\epsilon}{3}
    \end{array}
\end{equation}
for any $t\geq T$. In a similar fashion, we obtain the estimate
\begin{equation}
    \begin{array}{lcl}
         \int\limits_{\R}\big|\int\limits_t^{t+r} y(s-r)^\dagger \mathcal{K}(r;s-r)x(s)ds\big|dr&=&\int\limits_{-N}^{N}\big|\int\limits_t^{t+r} y(s-r)^\dagger \mathcal{K}(r;s-r)x(s)ds\big|dr \\[0.2cm]
         &&\qquad +\int\limits_{(-\infty,-N]\cup[N,\infty)}\big|\int\limits_t^{t+r} y(s-r)^\dagger \mathcal{K}(r;s-r)x(s)ds\big|dr\\[0.2cm]
         &\leq &\frac{\epsilon}{3}
    \end{array}
\end{equation}
for $t\geq T$. The representation (\ref{eq:diffhale1}) now allows us to estimate
\begin{equation}
    \begin{array}{lcl}
         |\ip{y^t,x_t}_t|&\leq &|y(t)| |x(t)|+\sum\limits_{j=-\infty}^\infty \big|\int\limits_t^{t+r_j}\big|y(s-r_j)^\dagger A_j(s-r_j)x(s)\big|ds\big|\\[0.2cm]
         &&\qquad +\int\limits_{\R}\big|\int\limits_t^{t+r} \big|y(s-r)^\dagger \mathcal{K}(r;s-r)x(s)\big|ds\big|dr\\[0.2cm]
         &\leq &\epsilon,
    \end{array}
\end{equation}
for any $t\geq T$, as desired.\qed\\

\begin{lemma}\label{lemma:lin:4.5:bewijs1}
Assume that (\asref{aannamesconstanten}{\text{H}}), (\asref{aannamesconvoluties}{\text{H}}) and (\asref{aannameslimit}{\text{H}}) and (\asref{aannames0ophalflijn}{\text{H}}).
Suppose furthermore that $\rmin=-\infty$ and consider
a pair of constants $(K_0,\alpha_0)\in\R_{> 0}^2$. 
Then there exists a positive constant $B>0$ so that
the estimate
\begin{equation}\label{eq:lin:Haleafsch4.5}
    \begin{array}{lcl}
         |\ip{\psi,y_0}_0|&\leq & B \nrm{\psi}_\infty\nrm{y_\tau}_\infty e^{-\alpha\tau}
    \end{array}
\end{equation}
holds for any $\tau\geq 0$,
 any  $\psi\in Y$ 
 and
 any $y\in C_b(D_{\tau}^\ominus)$ that satisfies the exponential bound
\begin{equation}\label{eq:expdecayy}
    \begin{array}{lcl}
         |y(t)|&\leq &K_0 e^{-\alpha_0(\tau-t)}\nrm{y_\tau}_\infty,
         \qquad \qquad
   t\leq \tau.        
    \end{array}
\end{equation}
\end{lemma}
\textit{Proof.} 
Recall the constants $ (K_{\mathrm{exp}},\alpha,p)\in\R_{>0}^3$ from Lemma \ref{thm:expdecaybewijs0}. By lowering $\alpha$ and increasing $ K_{\mathrm{exp}}$ if necessary, we may assume that $\alpha\leq \alpha_0$ and $K_{\mathrm{exp}}\geq K_0$. 
A first crude estimate yields
\begin{equation}
    \begin{array}{lcl}
         \Big|\sum\limits_{j=-\infty}^\infty \int\limits_{0}^{r_j}x(s-r_j)A_j(t+s-r_j)y(s)ds\Big|&\leq &\nrm{x}_\infty\sum\limits_{j=-\infty}^\infty  \nrm{A_j(\cdot)}_\infty \big|\int\limits_0^{r_j} y(s) ds \big|.
    \end{array}
\end{equation}
Splitting this sum into two parts and using the decay (\ref{eq:expdecayy}), we obtain the bound
\begin{equation}
    \begin{array}{lcl}
         \sum\limits_{j=-\infty}^\infty  \nrm{A_j(\cdot)}_\infty \big|\int\limits_0^{r_j} y(s)ds\big|&=&\sum\limits_{ r_j\leq \tau}  \nrm{A_j(\cdot)}_\infty \big|\int\limits_0^{r_j} y(s)ds\big|+\sum\limits_{ r_j> \tau}  \nrm{A_j(\cdot)}_\infty \int\limits_0^{r_j} |y(s)| ds\\[0.2cm]
         &\leq &\sum\limits_{ r_j\leq \tau}  \nrm{A_j(\cdot)}_\infty
         K_{\mathrm{exp}}
         \nrm{y_\tau}_\infty
         \big|\int\limits_0^{r_j}  e^{-\alpha(\tau-s)} ds\big|
         \\[0.2cm]
         & & \qquad
         +\sum\limits_{ r_j> \tau}  \nrm{A_j(\cdot)}_\infty r_j\nrm{y_0}_\infty\\[0.2cm]
         &\leq &\sum\limits_{ r_j\leq \tau}  \nrm{A_j(\cdot)}_\infty  K_{\mathrm{exp}}\nrm{y_\tau}_\infty\frac{1}{\alpha}\big|e^{\alpha(r_j-\tau)}-e^{-\alpha\tau}\big|
         \\[0.2cm]
         & & \qquad + K_{\mathrm{exp}} e^{-2\alpha \tau}\nrm{y_0}_\infty\\[0.2cm]
         &\leq & K_{\mathrm{exp}}\nrm{y_\tau}_\infty\frac{1}{\alpha}e^{-\alpha\tau}\sum\limits_{j=-\infty}^\infty  \nrm{A_j(\cdot)}_\infty + K_{\mathrm{exp}} e^{-2\alpha \tau}\nrm{y_\tau}_\infty,
    \end{array}
\end{equation}
where we used $\rmin=-\infty$ to conclude $\nrm{y_0}_\infty\leq \nrm{y_\tau}_\infty$. A similar computation for the convolution term 
yields the desired bound (\ref{eq:lin:Haleafsch4.5}).\qed\\

\subsection{Exponential splitting of the state space \texorpdfstring{$X$}{X}}
In the remainder of this section, we set out to establish Proposition \ref{thm:4.3equivbewijs3}
and complete the proof of
Theorem \ref{thm:4.3equiv}. In particular, the main technical goal is to establish the identity (\ref{eq:4.11}). 
We start by considering the inclusion $X^\perp(\tau) \subset S(\tau) $,
which will follow from Proposition
\ref{thm:4.3equivbewijs3} and the closedness of $S(\tau)$. In particular, we show that
$C^1(D_X) \cap X^\perp(\tau)$ is contained in $S(\tau)$. \\

Again, the main complication
is that the derivatives of functions $x$ in this subset need not be bounded, which hence also holds for $\Lambda x$. However, we do know that $\dot{x} - \Lambda x$ is bounded, which allows us to use a technical splitting of $x$ to achieve the desired result. In order to establish the consequences of this splitting, we will need to exploit the fundamental properties of the Hale inner product from \S \ref{subsec:fundhale}.

\begin{lemma}\label{thm:4.3equivbewijs2.5} Assume that (\asref{aannamesconstanten}{\text{H}}), (\asref{aannamesconvoluties}{\text{H}}) and (\asref{aannameslimit}{\text{H}}) are satisfied. Fix $\tau\in\R$ 
and pick
a differentiable function $x\in C_b(\R)\cap C^1(\R)$ with $\phi:=x_\tau\in X^\perp(\tau)$. Recall the operator $\Lambda$ from (\ref{def:Lambda}), write $h=\Lambda x$ and consider the functions
$h_-$ and $h_+$ given by
\begin{equation}
    \begin{array}{lclcl}
         h_-(t)&=&\begin{cases} h(t),\enskip &t\leq \tau,\\ 0,\enskip &t>\tau,\end{cases}\qquad h_+(t)&=&\begin{cases} 0,\enskip &t\leq \tau,\\ h(t),\enskip &t>\tau.\end{cases}
    \end{array}
\end{equation}
Then there exists a decomposition
$x = x_- + x_+$ with $x_-,x_+\in C_b(\R)\cap C^1(\R)$
for which we have the inclusions
\begin{equation}
\label{eq:exp:incl:h:pm}
    \begin{array}{lclcl}
         h_--\Lambda x_-&\in &\Range(\Lambda),\qquad  h_+-\Lambda x_+ &\in &\Range(\Lambda).
    \end{array}
\end{equation}
\end{lemma}
\textit{Proof.}
We choose the decomposition $x_-+x_+=x$ with $x_\pm \in C_b(\R) \cap C^1(\R)$ in such a way that $x_-=0$ on $[\tau+1,\infty)$, while $x_+=0$ on $(-\infty,\tau-1]$.
Although the derivative
of $x_-$ need not be bounded on $(-\infty,\tau]$, while the derivative of $x_+$ need not be bounded on $[\tau,\infty)$, we do claim that
\begin{equation}
    \begin{array}{lclcl}
         h_--\Lambda x_-&\in &L^\infty(\R),\qquad  h_+-\Lambda x_+&\in &L^\infty(\R).
    \end{array}
\end{equation}
To see this, we note that by construction
$\Lambda x^-$ is bounded on $[\tau, \infty)$, while $\Lambda x^+$ is bounded on $(-\infty, \tau]$. In particular, $h_+ - \Lambda x^+$
is automatically bounded on $(-\infty, \tau]$. On the other hand,
for $t \ge \tau$
we may compute
\begin{equation}
\label{eq:exp:id:for:h:plus:expr}
    \begin{array}{lclcl}
         h_+(t)-[\Lambda x_+](t) 
         &=& h(t)-
         [\Lambda x](t)
        +[\Lambda x_-](t)
         &=&[\Lambda x_-](t),
    \end{array}
\end{equation}
which shows that $h_+ - \Lambda x_+$
is also bounded on $[\tau, \infty)$. 
The claim for $x_-$ follows by symmetry.\\

We now set out to show
that $h_+-\Lambda x_+ \in\Range(\Lambda)$ by exploiting
the characterization
(\ref{eq:RangeLambda}).
In particular, pick any
$u\in\mathcal{B}^*$
and consider the integral
\begin{equation}
    \begin{array}{lcl}
         \mathcal{I}&:=&\int\limits_{-\infty}^{\infty} u(t)^\dagger \big[h_+(t)-(\Lambda x_+)(t)\big]dt\\[0.2cm]
         &=&\int\limits_{\tau}^\infty u(t)^\dagger \big[h_+(t)-(\Lambda x_+)(t)\big]dt+\int\limits_{-\infty}^\tau u(t)^\dagger \big[h_+(t)-(\Lambda x_+)(t)\big]dt\\[0.2cm]
         &=:&\mathcal{I}^++\mathcal{I}^-.
    \end{array}
\end{equation}
Exploiting \eqref{eq:exp:id:for:h:plus:expr}
we obtain
\begin{equation}
    \begin{array}{lcl}
         \mathcal{I}^+&=&
         \int\limits_{\tau}^\infty u(t)^\dagger \big[\dot{x}_-(t)-\sum\limits_{j=-\infty}^\infty A_j(t)x_-(t+r_j)-\int\limits_{\R}\mathcal{K}(\xi;t)x_-(t+\xi)d\xi\big]dt\\[0.2cm]
         &=&\int\limits_{\tau}^\infty u(t)^\dagger [\Lambda x_-](t) dt.
    \end{array}
\end{equation}
Since $\Lambda^* u=0$ we can immediately exploit the fundamental property of the Hale inner product from Lemma \ref{lemma:diffhale} to obtain
\begin{equation}
    \begin{array}{lcl}
         \mathcal{I}^+
         &=&\int\limits_{\tau}^\infty u(t)^\dagger [\Lambda x_-] (t)dt\\[0.2cm]
         &=&\int\limits_{\tau}^\infty \Big( u(t)^\dagger [\Lambda x_-](t)+[\Lambda^* u](t)x_-(t) \Big) \, dt\\[0.2cm]
         &=&\int\limits_{\tau}^\infty \frac{d}{dt}\ip{u^t,(x_-)_t}_t \, dt\\[0.2cm]
         &=&\lim\limits_{t\rightarrow\infty}\ip{u^t,(x_-)_t}_t-\ip{u^\tau,(x_-)_\tau}_\tau\\[0.2cm]
         &=&-\ip{u^\tau,(x_-)_\tau}_\tau.
    \end{array}
\end{equation}
The final equality follows in consideration of Lemma \ref{lemma:diffhale2}, since the function $u\in\mathcal{B}^*$ decays exponentially on account of Proposition \ref{prop:fredholmproperties}. In a similar fashion, we obtain
\begin{equation}
    \begin{array}{lcl}
         \mathcal{I}^-
         &=&-\ip{u^\tau,(x_+)_\tau}_\tau.
    \end{array}
\end{equation}
As such, we can
use $\phi\in X^\perp(\tau)$
to compute
\begin{equation}
    \begin{array}{lclclclcl}
         \mathcal{I}
         &=&
         \mathcal{I}^++\mathcal{I}^-
         &=&
         -\ip{u^\tau,x_\tau}_\tau
         &=&-\ip{u^\tau,\phi}_\tau
         &=&0.
    \end{array}
\end{equation}
The identity (\ref{eq:RangeLambda}) now yields the desired conclusion.\qed\\


\textit{Proof of Proposition \ref{thm:4.3equivbewijs3}.} Inspecting
the definition
of the Hale inner product (\ref{eq:defHale}),
we readily see that
the map $\phi\mapsto \ip{\psi,\phi}_\tau$ is continuous for any $\tau\in\R$ and any $\psi\in B^*(\tau)$.
In particular,
the space $X^\perp(\tau)$ is closed and has finite codimension in $X$. We now write
\begin{equation}
    \begin{array}{lcl}
         \mathcal{E}&=&C^1(D_X)\cap X^\perp(\tau)
    \end{array}
\end{equation}
and note that 
$\mathcal{E}$ is indeed dense in $X^\perp(\tau)$
by 
\cite[Lem. 4.14]{MPVL}.
Pick any $\phi \in \mathcal{E}$ and extend it arbitrarily to a bounded $C^1$ function $x:\R\rightarrow \C^M$
that has $x_\tau=\phi$. 
Recalling the functions $h_\pm$ and $x_\pm$ from Lemma \ref{thm:4.3equivbewijs2.5},
we use this result to find
a function 
$\tilde{y}\in W^{1,\infty}$ that has $\Lambda \tilde{y}=h_+-\Lambda x_+$.
Writing 
$y=\tilde{y}+x_+ \in C_b(\R)$, 
we see that
\begin{equation}
    \begin{array}{lcl}
        [\Lambda y](t) &=& h_+(t)  
    \end{array}
\end{equation}
which vanishes for $t \le \tau$.
In particular, we have $y \in \mathcal{P}(\tau)$. In a similar fashion, we can find a function $z \in \mathcal{Q}(\tau)$ with $\Lambda z = h_-$.\\

Writing $w = x - y - z \in C_b(\R) $, we readily compute
\begin{equation}
    \begin{array}{lclcl}
        \Lambda w &=&h-h_+-h_- &=&0,
    \end{array}
\end{equation}
which implies $w \in \mathcal{B}$
and hence
\begin{equation}
    \begin{array}{lclclcl}
        \phi &=& x_{\tau} &=& y_{\tau} + z_{\tau} + w_{\tau} 
        \in P(\tau) + Q(\tau) + B(\tau)
         &=& S(\tau),
    \end{array}
\end{equation}
as desired.
\qed\\

We now turn to the remaining inclusion $S(\tau)\subset X^\perp(\tau)$. As before, we exploit the fundamental identity (\ref{eq:fundHale}).
Combined with the exponential decay
of functions in $P(\tau)$ and $Q(\tau)$,
this will allow us to show that both spaces are contained in $X^\perp(\tau)$. \\

\begin{lemma}\label{thm:4.3equivbewijs2} Assume that (\asref{aannamesconstanten}{\text{H}}), (\asref{aannamesconvoluties}{\text{H}}) and (\asref{aannameslimit}{\text{H}}) are satisfied. 
Then for each $\tau \in \R$ we have 
the inclusion $S(\tau)\subset X^\perp(\tau)$.
\end{lemma}
\textit{Proof.} By symmetry, it suffices to show that $P(\tau)\subset X^\perp(\tau)$.
To this end, we pick
$x \in \mathcal{P}(\tau)$
and $y \in \mathcal{B}^*$ and note
that 
\begin{equation}
    \begin{array}{lcl}
         \frac{d}{dt}\ip{y^t,x_t}_t&=&0
    \end{array}
\end{equation}
for all $t \le \tau$ by Lemma \ref{lemma:diffhale}.
Since 
$x(t)$ is bounded as $t\rightarrow-\infty$ while $y(t)\rightarrow 0$ at an exponential rate, 
we may use 
Lemma \ref{lemma:diffhale2}
to obtain
\begin{equation}\label{eq:diffhale4}
    \begin{array}{lclcl}
         \ip{y^{\tau},x_\tau}_\tau
         &=&
         \lim\limits_{t\rightarrow-\infty}\ip{y^t,x_t}_t
         &=&
         0,
    \end{array}
\end{equation}
as desired.
\qed\\

The remainder of the proof of 
Theorem \ref{thm:4.3equiv} uses arguments that are very similar
to those in \cite{MPVL}. The main point is that the compactness properties obtained in {\S}\ref{sec:ex:restr} allow us to show that $\widehat{P}(\tau)$ and $\widehat{Q}(\tau)$ are closed, which allows the computations above to be leveraged.

\begin{lemma}\label{thm:4.3equivbewijs0.5}
Assume that (\asref{aannamesconstanten}{\text{H}}), (\asref{aannamesconvoluties}{\text{H}}) and (\asref{aannameslimit}{\text{H}}) are satisfied. Then for each $\tau\in\R$, the spaces $P(\tau)$, $Q(\tau)$, $\widehat{P}(\tau)$ and $\widehat{Q}(\tau)$
are all  closed subspaces of $X$. 
\end{lemma}
\textit{Proof.} Let $\{\phi_n\}_{n\geq 1}$ be a sequence in $\widehat{P}(\tau)$ that converges in $X$ to some $\phi\in X$. Picking $\tau_n=\tau$ and $\psi_n=0$ in Corollary \ref{cor:4.11equiv} then immediately implies that $\{\phi_n\}_{n\geq 1}$ converges uniformly on compact sets to some $\hat{\phi}\in \widehat{P}(\tau)$.
By necessity we hence have
$\phi =\hat{\phi}$, which means that
$\widehat{P}(\tau)$
and by symmetry $\widehat{Q}(\tau)$
are both closed. This subsequently must also hold for the finite dimensional extensions $P(\tau)$ and $Q(\tau)$.
\qed\\

\begin{lemma}
\label{thm:4.3equivbewijs1} Assume that (\asref{aannamesconstanten}{\text{H}}), (\asref{aannamesconvoluties}{\text{H}}) and (\asref{aannameslimit}{\text{H}}) are satisfied. Then for each $\tau\in\R$ the spaces $S(\tau)$ and $\widehat{S}(\tau)$ are closed subspaces of $X$. Moreover, 
the decompositions \eqref{eq:4.10} hold.
\end{lemma}

\textit{Proof.} 
In view of Proposition \ref{compactoperator} and Lemma \ref{thm:4.3equivbewijs0.5},
the result can be obtained
by following the proof of \cite[Prop. 4.12 \& Prop. 4.13]{MPVL}, together with the first part of the proof of \cite[Thm. 4.3]{MPVL}.
\qed\\

\textit{Proof of Theorem \ref{thm:4.3equiv}.} Every statement except the identity (\ref{eq:4.11}) follows from Lemma \ref{thm:4.3equivbewijs0.5} and Proposition \ref{thm:4.3equivbewijs1}. In addition, Lemma \ref{thm:4.3equivbewijs2} yields the inclusion $S(\tau)\subset X^\perp(\tau)$, while Proposition \ref{thm:4.3equivbewijs3} yields the inclusion $D\subset S(\tau)$ for some dense set $D\subset X^\perp(\tau)$. Since $S(\tau)$ is closed, we immediately obtain (\ref{eq:4.11}). \qed\\

\section{Fredholm properties of the projections \texorpdfstring{$\Pi_{\widehat{P}}$}{PiP} and \texorpdfstring{$\Pi_{\widehat{Q}}$}{PiQ}}\label{section:fredholm}

The goal of this section is to understand the projection operators $\Pi_{\widehat{P}}$ and $\Pi_{\widehat{Q}}$ associated to the decomposition
(\ref{eq:4.14}). 
In contrast to the previous section,
we can follow the approach from \cite{MPVL} relatively smoothly here. The main difficulty is that the arguments in \cite{MPVL}
often use Corollary \ref{cor:4.11equiv} to conclude
that certain subsequences converge uniformly,
while we can only conclude that this convergence takes place on compact subsets.
The primary way in which we circumvent this issue is by appealing to the exponential estimates in Theorem \ref{thm:expdecay}.\\

As a bonus, we also obtain information on the Fredholm properties of the restriction operators $\pi^\pm$ introduced in {\S}\ref{sec:ex:restr}. In particular,
besides proving Theorem \ref{theorem4.6}, we also establish the following two results. \\

\begin{proposition}[{cf. \cite[Thm. 4.5]{MPVL}}]\label{theorem4.5} Assume that (\asref{aannamesconstanten}{\text{H}}), (\asref{aannamesconvoluties}{\text{H}}) and (\asref{aannameslimit}{\text{H}}) are satisfied. Then the operators $\pi^+_{P(\tau)}$, 
$\pi^-_{Q(\tau)}$, 
$\pi^+_{\widehat{P}(\tau)}$ and $\pi^-_{\widehat{Q}(\tau)}$ are all Fredholm  for every $\tau\in\R$. Recalling the function $\beta(\tau)$ defined in (\ref{eq:4.12}), the Fredholm indices
satisfy the identities
\begin{equation}\label{eq:4.13}
    \begin{array}{lcl}
         \ind(\pi^+_{P(\tau)})+\ind(\pi^-_{Q(\tau)})&=&-M+\dim B(\tau)-\beta(\tau),\\[0.2cm]
         \ind(\pi^+_{\widehat{P}(\tau)})+\ind(\pi^-_{\widehat{Q}(\tau)})&=&-\big(M+\dim B(\tau)+\beta(\tau)\big).
    \end{array}
\end{equation}
\end{proposition}

\begin{proposition}[{cf. \cite[Thm. 4.6]{MPVL}}]\label{theorem4.6:indices} Assume that (\asref{aannamesconstanten}{\text{H}}), (\asref{aannamesconvoluties}{\text{H}}) and (\asref{aannameslimit}{\text{H}}) are satisfied. Fix $\tau_0\in\R$
and consider the projections
$\Pi_{\widehat{P}}$ and $\Pi_{\widehat{Q}}$ 
associated to the decomposition
\eqref{eq:4.14}. Then we have the identities
\begin{equation}\label{eq:4.18}
    \begin{array}{lcl}
        \ind(\pi_{\widehat{P}(\tau)}^+)&=&\ind(\pi^+_{\widehat{P}(\tau_0)})-\codim_{\widehat{P}(\tau_0)}\Pi_{\widehat{P}}\big(\widehat{P}(\tau)\big),\\[0.2cm]
        \ind(\pi_{\widehat{Q}(\tau)}^-)&=&\ind(\pi^-_{\widehat{Q}(\tau_0)})-\codim_{\widehat{Q}(\tau_0)}\Pi_{\widehat{Q}}\big(\widehat{Q}(\tau)\big).
    \end{array}
\end{equation}
Moreover, the quantities $\ind(\pi_{\widehat{P}(\tau)}^+)$ and $\ind(\pi_{\widehat{Q}(\tau)}^-)$ vary upper semicontinuously with $\tau$. In addition, we have the identities
\begin{equation}\label{eq:4.20:indices}
    \begin{array}{lclcllcl}
         \ind(\pi_{\widehat{P}(\tau)}^+)+\dim B(\tau)&=&\ind(\pi_{P(\tau)}^+)&=&\ind(\pi_{P(-\infty)}^+),\\[0.2cm]
         \ind(\pi_{\widehat{Q}(\tau)}^-)+\dim B(\tau)&=&\ind(\pi_{Q(\tau)}^-)&=&\ind(\pi_{Q(\infty)}^-),
    \end{array}
\end{equation}
for sufficiently negative values of $\tau$ in the first line of (\ref{eq:4.20:indices}) and for sufficiently positive values of $\tau$ in the second line of (\ref{eq:4.20:indices}).

\end{proposition}

We first need to study the projection operators $\pi^+$ and $\pi^-$ from (\ref{eq:defpi}) in more detail. We proceed largely along the lines of \cite{MPVL}, taking a small detour in order to establish that the ranges are closed.

\begin{lemma}
\label{lemma3.8.1} Assume that (\asref{aannamesconstanten}{\text{H}}), (\asref{aannamesconvoluties}{\text{H}}) and (\asref{aannameslimit}{\text{H}}) are satisfied. Then the operators $\pi^+_{\widehat{P}(\tau)}$ and $\pi^-_{\widehat{Q}(\tau)}$ have finite dimensional kernels for each $\tau\in\R$.
\end{lemma}
\textit{Proof.} This can be established by repeating the first
half of the proof of \cite[Lem. 3.8]{MPVL}.
\qed\\

\begin{lemma}[{cf. \cite[Lem. 3.8]{MPVL}}]\label{lemma3.8.2} Assume that (\asref{aannamesconstanten}{\text{H}}), (\asref{aannamesconvoluties}{\text{H}}) and (\asref{aannameslimit}{\text{H}}) are satisfied. Then the operators $\pi^+_{\widehat{P}(\tau)}$ and $\pi^-_{\widehat{Q}(\tau)}$ have closed ranges for each $\tau\in\R$.
\end{lemma}
\textit{Proof.} By symmetry, we pick $\tau\in\R$ and restrict attention to the operator $\pi^+_{\widehat{P}(\tau)}$. We fix a closed complement $C\subset \widehat{P}(\tau)$ for the finite dimensional space $\ker(\pi^+_{\widehat{P}(\tau)})$, so that $\widehat{P}(\tau)=\ker(\pi^+_{\widehat{P}(\tau)})\oplus C$. 
We now consider a sequence
$\{\phi_n\}_{n \ge 1} \subset C$ 
and suppose that the restrictions
$\psi_n=\pi^+_{\widehat{P}(\tau)}\phi_n$
satisfy the uniform convergence $\psi_n\rightarrow \psi$ on $D_X^+$.
If the sequence $\{\phi_n\}_{n \ge 1}$ is bounded, then an application of Corollary \ref{cor:4.11equiv} immediately yields that $\phi_n\rightarrow \phi\in\widehat{P}(\tau)$ uniformly on compacta, after passing to a subsequence. This implies 
that $\psi=\pi^+_{\widehat{P}(\tau)}\phi$ and thus $\psi\in\Range(\pi^+_{\widehat{P}(\tau)})$, as desired.\\

Let us assume therefore that 
$\nrm{\phi_n}_\infty\uparrow \infty$ and consider the rescaled
sequence
$\tilde{\phi}_n=\nrm{\phi_n}_\infty^{-1}\phi_n$,
which satisfies
\begin{equation}
    \begin{array}{lcl}
    \pi^+_{\widehat{P}(\tau)}\tilde{\phi}_n&=&\nrm{\phi_n}_{\infty}^{-1}\psi_n\rightarrow 0
    \end{array}
\end{equation}
uniformly on $D_X^+$.
We may again apply Corollary \ref{cor:4.11equiv} to obtain $\tilde{\phi}_n\rightarrow \tilde{\phi} \in \widehat{P}(\tau)$ uniformly on compacta,
with $\pi^+ \tilde{\phi} = 0$.
In contrast to the setting
of \cite[Lem. 3.8]{MPVL}, this convergence is not immediately
uniform on the 
(possibly unbounded) interval $D_X^-$.
On account of Proposition \ref{compactoperator}, the operator $\pi^-_{\widehat{P}(\tau)}$ is compact, so we can pass to yet another subsequence to obtain the limit $\tilde{\phi}_n\rightarrow \tilde{\phi}$ uniformly on $D_X^-$. As such, $\tilde{\phi}_n\rightarrow\tilde{\phi}$ uniformly both on $D_X^-$ and on $D_X^+$, so the convergence is uniform on $D_X$.  Moreover, $\pi^+_{\widehat{P}(\tau)}\tilde{\phi}=0$, so $\tilde{\phi}\in\ker(\pi^+_{\widehat{P}(\tau)})$. Since the convergence $\tilde{\phi}_n\rightarrow\tilde{\phi}$ is uniform on $D_X$ we get $\nrm{\tilde{\phi}}_\infty=1$, as $\nrm{\tilde{\phi}_n}_\infty=1$ for each $n$. However, $C$ is closed and $\tilde{\phi}_n\in C$ for each $n$, so $\tilde{\phi}\in C$. Therefore, $\tilde{\phi}$ is a non-trivial element of $\ker(\pi^+_{\widehat{P}(\tau)})\cap C$, which yields a contradiction.\qed\\


\textit{Proof of Proposition \ref{theorem4.5}.} The proof is identical to that of \cite[Prop. 4.12]{MPVL} and, as such, will be omitted. It uses Theorem \ref{thm:4.3equiv}, together with Lemmas \ref{lemma3.8.1} and \ref{lemma3.8.2}.\qed\\


\begin{lemma}\label{theorem4.6bewijs1} Assume that (\asref{aannamesconstanten}{\text{H}}), (\asref{aannamesconvoluties}{\text{H}}) and (\asref{aannameslimit}{\text{H}}) are satisfied. 
Fix $\tau_0\in\R$
and consider the projections
$\Pi_{\widehat{P}}$ and $\Pi_{\widehat{Q}}$ 
associated to the decomposition
\eqref{eq:4.14}. Then for $\tau$ sufficiently close to $\tau_0$, the restrictions
\begin{equation}\label{eq:4.16bewijs}
    \begin{array}{lcl}
         \Pi_{\widehat{P}}:\widehat{P}(\tau)&\rightarrow &\Pi_{\widehat{P}}\big(\widehat{P}(\tau)\big)\subset \widehat{P}(\tau_0),\\[0.2cm]
         \Pi_{\widehat{Q}}:\widehat{Q}(\tau)&\rightarrow &\Pi_{\widehat{Q}}\big(\widehat{Q}(\tau)\big)\subset \widehat{Q}(\tau_0)
    \end{array}
\end{equation}
to the subspaces $\widehat{P}(\tau)$ and $\widehat{Q}(\tau)$ are isomorphisms onto their ranges, which are closed. Moreover, we have the limits
\begin{equation}\label{eq:4.17bewijs}
    \begin{array}{lclcl}
         \lim\limits_{\tau\rightarrow\tau_0}\nrm{I-\Pi_{\widehat{P}}|_{\widehat{P}(\tau)}}&=&0,
         \qquad \qquad \lim\limits_{\tau\rightarrow\tau_0}\nrm{I-\Pi_{\widehat{Q}}|_{\widehat{Q}(\tau)}}&=&0,
    \end{array}
\end{equation}
in which $I$ denotes the inclusion of $\widehat{P}(\tau)$ or $\widehat{Q}(\tau)$ into $X$.\end{lemma}
\textit{Proof.} By symmetry, we only have to consider the projection $\Pi_{\widehat{P}}$. 
In order to establish the limit
(\ref{eq:4.17bewijs}),
we pick an arbitrary bounded sequence $\{\phi_n\}_{n\geq 1}$
that has $\phi_n\in \widehat{P}(\tau_n)$ and $\tau_n\rightarrow \tau_0$. Using the decomposition (\ref{eq:4.14}), 
we write
\begin{equation}
    \begin{array}{lcl}
         \phi_n&=&\rho_n+\psi_n+\sigma_n,
    \end{array}
\end{equation}
with $\rho_n\in \widehat{P}(\tau_0)$, $\psi_n\in \widehat{Q}(\tau_0)$ and $\sigma_n\in \Gamma$ for each $n$. Then each of the sequences $\{\rho_n\}_{n\geq 1}$, $\{\psi_n\}_{n\geq 1}$ and $\{\sigma_n\}_{n\geq 1}$ is bounded. It is sufficient to show that
$\phi_n-\rho_n\rightarrow 0$
for some subsequence. Note that this also establishes the claim that the restriction in (\ref{eq:4.16bewijs}) is an isomorphism with closed range. \\

By Proposition \ref{compactoperator} and the finite dimensionality of $\Gamma$, we can pass over to a subsequence for which both $\{\pi^+_{\widehat{Q}(\tau_0)}\psi_n\}_{n\geq 1}$ and $\{\sigma_n\}_{n\geq 1}$ converge. As such, $\{\pi^+(\phi_n-\rho_n)\}_{n\geq 1}$ converges, so Corollary \ref{cor:4.11equiv} implies that $\phi_n-\rho_n\rightarrow \phi\in\widehat{P}(\tau_0)$ uniformly on compact subsets of $D_X$ after passing to a further subsequence. In particular, we obtain the convergence $\psi_n+\sigma_n\rightarrow \phi$, uniformly on $D_X^+$ and uniformly on compact subsets of $D_X^-$.\\

If $\rmin\neq-\infty$ then the convergence $\psi_n+\sigma_n\rightarrow \phi$ is in fact uniform on $D_X$, allowing us to follow
the approach in \cite{MPVL}.
In particular, we obtain
$\phi\in \big(\widehat{Q}_{\tau_0}\oplus \Gamma \big)\cap \widehat{P}_{\tau_0}$
and hence $\phi=0$ as desired.
Assuming therefore that $\rmin=-\infty$, we use
the convergence of $\{\sigma_n\}_{n\geq 1}$ 
to conclude that
 $\psi_n\rightarrow \psi$ uniformly on $D_X^+$ and uniformly on compact subsets of $(-\infty,0]$. For any $n\geq 1$ we write $y_n\in \widehat{\mathcal{Q}}(\tau_0)$
 for the right-extension of $\psi_n$, i.e.,
 $\psi_n=(y_n)_{\tau_0}$.
 Using the uniform estimates in Theorem \ref{thm:expdecay} for large positive $t$, we can use the Ascoli-Arzela theorem to pass to a subsequence that has $y_n\rightarrow y$, uniformly on compact subsets of $\R$. Necessarily we have 
 \begin{equation}
    \begin{array}{lcl}
         y_{\tau_0}(t)&=&\psi(t),
         \qquad \qquad
         t \in \overline{(-1,\rmax)}.
    \end{array}
 \end{equation} Since $\psi_n\rightarrow \psi$ uniformly on $D_X^+$, it follows that $y_n\rightarrow y$ uniformly on $\tau_0+D_X^+$. We can hence follow the proof of Lemma \ref{thm:expdecaybewijs1bewijs2} to see that $y$ is a solution of (\ref{ditishetprobleem}) on $[\tau_0,\infty)$ and therefore $\psi\in Q(\tau_0)$. Similarly to the proof of Corollary \ref{cor:4.11equiv} we even get $\psi\in \widehat{Q}(\tau_0)$. This yields $\phi\in \big(\widehat{Q}(\tau_0)\oplus \Gamma\big)\cap \widehat{P}(\tau_0)$ and therefore $\phi=0$.\qed \\

\textit{Proof of Theorem \ref{theorem4.6}.} 
The first statement and (\ref{eq:4.17}) follow from Lemma \ref{theorem4.6bewijs1}, while the lower semicontinuity of $\dim B(\tau)$ and $\beta(\tau)$ and the limit in (\ref{eq:4.19}) can be established in a fashion similar to the proof of \cite[Thm. 4.6]{MPVL}.\\

It remains to show that (\ref{eq:4.20}) holds. Following \cite{MPVL}, it suffices to find a bounded function $y:\R\rightarrow\C^M$ that satisfies the inhomogeneous system
\begin{equation}\label{eq:4.36}
    \begin{array}{lcl}
         \dot{y}(t)&=&\sum\limits_{j=-\infty}^\infty A_j^{(\tau)} (t)y(t+r_j)+\int\limits_{\R}\mathcal{K}^{(\tau)}(\xi;t)y(t+\xi)d\xi+h^{(\tau)}(t),
    \end{array}
\end{equation}
in which we have introduced
the coefficients 
\begin{equation}
    \begin{array}{lclcl}
         A_j^{(\tau)}(t)&=&\begin{cases} A_j(t+\tau),\\ A_j(-\infty), \end{cases}\qquad \begin{array}{l}t<0,\\ t\geq 0,\end{array}
          \qquad \qquad
         \mathcal{K}^{(\tau)}(\xi;t)
         &=&\begin{cases} \mathcal{K}(\xi;t+\tau),\\ \mathcal{K}(\xi;-\infty),\end{cases}\qquad \begin{array}{l}t<0,\\ t\geq 0,\end{array}
    \end{array}
\end{equation}
together with the inhomogeneity
\begin{equation}
    \begin{array}{lcll}
         h^{(\tau)}(t)&=&\sum\limits_{j=-\infty}^\infty \big(A_j^\tau(t)-A_j(-\infty)\big)x(t+r_j) +\int\limits_{\R}\big(\mathcal{K}^\tau(\xi;t)-\mathcal{K}(\xi;-\infty)\big)x(t+\xi)d\xi.
    \end{array}
\end{equation}
This can be achieved by following the same steps as in \cite{MPVL}, 
but now using the proof of \cite[Lem. 3.1 (step 3)]{Faye2014}\footnote{The matrices $A_j^\tau$ and $\mathcal{K}^\tau$ need not be continuous, while in \cite{Faye2014} the coefficients are assumed to be continuous. However, the continuity is not used in the parts of the proof that are relevant for us.} instead of the results in \cite{MPA}. \qed\\

\textit{Proof of Proposition \ref{theorem4.6:indices}.} The proof is identical to that of \cite[Thm. 4.6]{MPVL} and, as such, will be omitted. It uses Theorems \ref{thm:4.3equiv} and \ref{theorem4.6}.
\qed\\ 

\section{Exponential Dichotomies on half-lines}\label{section:halfline}

In this section, we adapt the approach 
of \cite{HJHLIN} to obtain exponential splittings for (\ref{ditishetprobleem})
on the half-line $[0, \infty)$. The main idea
is to explicitly construct  suitable finite-dimensional enlargements of $\mathcal{P}(\tau)$ for $\tau \ge 0$.
The extra functions $\{y_{(\tau)}\}_{\tau \ge 0}$ satisfy
(\ref{ditishetprobleem}) on $[0, \tau]$,
but not on $(-\infty, \tau]$. In fact,
we will exploit the fundamental identity
\eqref{eq:fundHale} to guarantee that the
segments $\{(y_{(\tau)})_{\tau} \}$
are not contained in $S(\tau)$.\\

In order to achieve this, we need to construct inverses for the Fredholm operator $\Lambda$ restricted to half-lines. In the ODE case one can write down explicit variation-of-constants formula's to achieve this, but such constructions are problematic at best in the current setting. Instead, we follow \cite{HJHLIN} and solve $\Lambda x = h$ by appropriately modifying $h$ \textit{outside} the half-line of interest in order to satisfy $\langle y , h \rangle_{L^2} = 0$
for all $y \in \mathcal{B}^*$. 
In order to ensure that such a modification is not precluded by degeneracy issues,
we need to assume that (\asref{aannames0ophalflijn}{\text{H}}) holds.\\

The main complication in the setting $r_{\min} = -\infty$ is that  this modification of $h$
is visible directly in the equation satisfied by $x$, rather than only indirectly via the Fredholm properties of $\Lambda$ as in \cite{HJHLIN}. This raises issues when using a standard bootstrapping procedure to obtain estimates on $\dot{x}$. Naturally, the unbounded shifts also cause technical problems similar to those encountered in {\S}\ref{section:existence}-\ref{section:fredholm}, but fortunately the same tricks also work here.

\begin{remark}
In fact, in this section it suffices to assume a weaker version of the non-triviality condition (\asref{aannames0ophalflijn}{\text{H}}). In particular, we do not need the condition that each nonzero $d\in\mathcal{B}$ vanishes on the intervals $(-\infty,\tau]$ for $\tau<0$ or $[\tau,\infty)$ for $\tau>0$. This is because the formulation of Theorem \ref{thm:lin:4.1}  references the specific half-line $\R^+$,
rather than arbitrary half-lines.
\end{remark}

\subsection{Strategy}

In order to ensure that the spaces we construct are invariant with respect to $\tau$, we need to slightly modify the $\tau$-dependent normalization condition
used in \eqref{defhatmathcalspaces}. Indeed,
upon writing
\begin{equation}
    \begin{array}{lcl}
         \tilde{\mathcal{P}}(\tau)&=&\{x\in \mathcal{P}(\tau)\text{ }|\text{ }\int_{-\infty}^{\min(\tau+\rmax,0)}y(t)^\dagger x(t)dt=0\text{ for every }y\in\mathcal{B}\},\\[0.2cm]
         \tilde{\mathcal{Q}}(\tau)&=&\{x\in \mathcal{Q}(\tau)\text{ }|\text{ }\int^{\max(\tau+\rmin,0)}_\infty y(t)^\dagger x(t)dt=0\text{ for every }y\in\mathcal{B}\},\\[0.2cm]
         \tilde{P}(\tau)&=&\{\phi \in X\text{ }|\text{ }\phi=x_\tau\text{ for some }x\in\tilde{\mathcal{P}}(\tau)\},\\[0.2cm]
         \tilde{Q}(\tau)&=&\{\phi \in X\text{ }|\text{ }\phi=x_\tau\text{ for some }x\in\tilde{\mathcal{Q}}(\tau)\},
    \end{array}
\end{equation}
we see that the upper bounds for the defining integrals 
are now constant for $\tau\geq 0$ and $\tau\leq 0$ respectively.
In view of the
non-triviality assumption (\asref{aannames0ophalflijn}{\text{H}}),
all the conclusions from the previous sections 
remain valid for these new spaces.
In particular, 
we have the following result.
\begin{corollary}[{cf. \cite[Prop. 4.2]{HJHLIN}}]
Assume that (\asref{aannamesconstanten}{\text{H}}), (\asref{aannamesconvoluties}{\text{H}}), (\asref{aannameslimit}{\text{H}}) and (\asref{aannames0ophalflijn}{\text{H}}) are satisfied. Recall the spaces $S(\tau)$ from Theorem \ref{thm:4.3equiv}. Then we have the direct sum decomposition
\begin{equation}
    \begin{array}{lcl}
         S(\tau)&=&\tilde{P}(\tau)\oplus\tilde{Q}(\tau)\oplus B(\tau)
    \end{array}
\end{equation}
for any $\tau\in\R$.
\end{corollary}

Our first goal is to find an explicit complement for the space $S(\tau)$ in $X$. In view of the identification $S(\tau) = X^\perp(\tau)$,
we actually build a duality basis for $B^*(\tau)$ 
with respect to the Hale inner product; see \eqref{eq:lin:4.9}. 
\begin{proposition}[{cf. \cite[Lem. 4.3]{HJHLIN}}]\label{lemma:lin:4.3}
Assume that (\asref{aannamesconstanten}{\text{H}}), (\asref{aannamesconvoluties}{\text{H}}), (\asref{aannameslimit}{\text{H}}) and (\asref{aannames0ophalflijn}{\text{H}}) are satisfied. Write $n_d=\dim(\mathcal{B}^*)$ and choose a basis $\{d^i:1\leq i\leq n_d\}$ for $\mathcal{B}^*$. Then there exists a constant $r_0>0$, together with a family of functions $y_{(\tau)}^i\in C_b(D_{\tau}^\ominus)$, defined for every $\tau\geq 0$ and every integer $1\leq i\leq n_d$, that satisfies the following properties.
\begin{enumerate}[label=(\roman*)]
    \item For any $\tau\geq 0$ and any integer $1\leq i\leq n_d$ we have $\big[\Lambda y^i_{(\tau)}\big](t)=0$ for every $t\in(\infty,-r_0]\cup[0,\tau]$.
    \item For any $0\leq t\leq \tau$ and any pair $1\leq i,j\leq n_d$ we have the identity
    \begin{equation}\label{eq:lin:4.9}
        \begin{array}{lcl}
             \ip{(d^i)^t,(y^j_{(\tau)})_t}_t&=&\delta_{ij}.
        \end{array}
    \end{equation}
    \item For any fixed constant $t\geq 0$ and fixed integer $1\leq i\leq n_d$, the map $\tau\mapsto (y^i_{(\tau)})_t$ is continuous from the interval $[t,\infty)$ into the state space $X$.
    \item For any triplet $0\leq t\leq \tau_1\leq \tau_2$ and any integer $1\leq i\leq n_d$, we have the inclusion
    \begin{equation}
        \begin{array}{lcl}
             \big[y^i_{(\tau_1)}-y^i_{(\tau_2)}\big]_t\in \tilde{P}(t).
        \end{array}
    \end{equation}
    \item For any $\tau\geq 0$ and any integer $1\leq i\leq n_d$, the integral condition
    \begin{equation}\label{eq:lin:4.11}
        \begin{array}{lcl}
             \int_{-\infty}^0 b(t)^\dagger y_{(\tau)}(t)dt&=&0
        \end{array}
    \end{equation}
    holds for all $b\in\mathcal{B}$.
\end{enumerate}
\end{proposition}
Upon using 
the functions in Proposition \ref{lemma:lin:4.3} to introduce the 
finite-dimensional spans
\begin{equation}\label{eq:lin:Ytau}
    \begin{array}{lclcl}
         \mathcal{Y}(\tau)
         &=&\span\big\{y^i_{(\tau)}\big\}_{i=1}^{n_d},
         \qquad \qquad
         Y(\tau)
         &=&
         \span\big\{(y^i_{(\tau)})_\tau\big\}_{i=1}^{n_d},
    \end{array}
\end{equation}
we can now define the spaces $R(\tau)$ and $\mathcal{R}(\tau)$
that appear in Theorem \ref{thm:lin:4.1} by writing
\begin{equation}\label{eq:lin:4.14}
    \begin{array}{lclcl}
         \mathcal{R}(\tau)
         &=&
         
         \tilde{\mathcal{P}}(\xi)\oplus\mathcal{Y}(\tau),
         \qquad \qquad
         R(\tau)
         &=&
         
         \tilde{P}(\xi)\oplus Y(\tau).
    \end{array}
\end{equation}
The identities in (\ref{eq:lin:4.9}) show that the dimension of the space $Y(\tau)$ is precisely $n_d$. Moreover, in combination with Theorem \ref{thm:4.3equiv} they yield 
\begin{equation}
    \begin{array}{lcl}
         S(\tau)\cap Y(\tau)&=&\{0\},
    \end{array}
\end{equation}
which means that we have the direct sum decomposition
\begin{equation}\label{eq:lin:4.33}
    \begin{array}{lcl}
         X&=&\tilde{P}(\tau)\oplus Y(\tau)\oplus Q(\tau)
    \end{array}
\end{equation}
for any $\tau\geq 0$. \\

Our final main result here generalizes the exponential decay estimates
contained in Theorem \ref{thm:expdecay} to the half-line setting. The main 
obstacle here is that it is more involved to control the derivative of functions in $\mathcal{Y}(\tau)$, preventing a direct application of the Ascoli-Arzela theorem. Indeed, these functions have
a non-zero right-hand side on the interval $[-r_0, 0]$
when substituted into (\ref{ditishetprobleem}).


\begin{proposition}[{cf. \cite[Prop. 4.4]{HJHLIN}}]\label{prop:lin:4.4}
Assume that (\asref{aannamesconstanten}{\text{H}}), (\asref{aannamesconvoluties}{\text{H}}), (\asref{aannameslimit}{\text{H}}) and (\asref{aannames0ophalflijn}{\text{H}}) are satisfied. Then for any $\tau \ge 0$,
every function $x\in \mathcal{R}(\tau)$ is $C^1$-smooth on $(-\infty,\tau]$.
In addition,
there exist constants $K_{\mathrm{dec}}>0$ and $\alpha>0$ in such a way that for all $\tau\geq 0$ and all $t\leq \tau$ we have the pointwise estimate
\begin{equation}
    \begin{array}{lcl}
         |x(t)|+|\dot{x}(t)|&\leq &K_{\mathrm{dec}}e^{-\alpha (\tau-t)}\nrm{x_\tau}_\infty
    \end{array}
\end{equation}
for every $x\in \mathcal{R}(\tau)$. 
\end{proposition}

\subsection{Construction of \texorpdfstring{$\mathcal{Y}(\tau)$}{Y(tau)}}

In order to construct the functions $\{y^i_{(\tau)}\}$ from Proposition \ref{lemma:lin:4.3}, we will use the freedom we still have to choose complements for the range and the kernel of the operator $\Lambda$.
\begin{lemma}[{cf. \cite[Lem. 3.4]{HJHLIN}}]\label{lin:Lemma3.4} Assume that (\asref{aannamesconstanten}{\text{H}}), (\asref{aannamesconvoluties}{\text{H}}), (\asref{aannameslimit}{\text{H}}) and (\asref{aannames0ophalflijn}{\text{H}}) are satisfied and fix $\tau\in\R$. Write $n_d=\dim(\mathcal{B}^* )$ and choose a basis $\{d^i:1\leq i\leq n_d\}$ for $\mathcal{B}^*$. Then there exists a constant $r_0^{(\tau)}>0$, together with functions 
\begin{equation}
\{\phi^i_{(\tau)}:1\leq i\leq n_d\}\subset C_b[\tau,\tau+r_0^{(\tau)}],
\qquad \qquad
\{\psi^i_{(\tau)}:1\leq i\leq n_d\}\subset C_b[\tau-r_0^{(\tau)},\tau]
\end{equation}
that satisfy the identities
\begin{equation}
\begin{array}{lclcl}
\int\limits_{\tau}^{\tau+r_0^{(\tau)}}d^i(t)^\dagger \phi^j_{(\tau)}(t)dt&=&\delta_{i,j},\\[0.2cm]
\int\limits_{\tau-r_0^{(\tau)}}^{\tau}d^i(t)^\dagger \psi^j_{(\tau)}(t)dt&=&\delta_{i,j}\\[0.2cm]
\end{array}
\end{equation}
for any $1\leq i,j\leq n_d$, together with
\begin{equation}\label{eq:normalizationpsij}
\begin{array}{rclcl}
0&=&\phi^j_{(\tau)}(\tau)&=&\phi^j_{(\tau)}(\tau+r_0^{(\tau)}),\\[0.2cm]
0&=&\psi^j_{(\tau)}(\tau-r_0^{(\tau)})&=&\psi^j_{(\tau)}(\tau),
\end{array}
\end{equation}
for any $1\leq j\leq n_d$.
\end{lemma}
\textit{Proof.} By symmetry, we only consider the construction of the functions $\{\psi^i_{(\tau)}:1\leq i\leq n_d\}$. We first note that the restriction operator 
\begin{equation}
\mathcal{B}^*\rightarrow C_b[-r_0^{(\tau)}+\tau,\tau],
\qquad d\mapsto d|_{[-r_0^{(\tau)}+\tau,\tau]}  \end{equation}
is injective for some $r_0^{(\tau)}>0$. This follows trivially from (\asref{aannames0ophalflijn}{\text{H}}) and the fact that $\mathcal{B}^*$ is finite dimensional.\\

Let us denote $[\cdot,\cdot]_\tau$ for the integral product
\begin{equation}
    \begin{array}{lcl}
         [\psi,\phi]_\tau&=&\int\limits_{\tau-r_0^{(\tau)}}^{\tau}\psi(t)^\dagger \phi(t)dt.
    \end{array}
\end{equation}
Consider any set of functions $\{\tilde{\psi}^i:1\leq i\leq n_d\}\subset C_b[\tau-r_0^{(\tau)},\tau]$ with 
\begin{equation}
\begin{array}{rclcl}
0&=&\tilde{\psi}^j(\tau-r_0^{(\tau)})&=&\tilde{\psi}^j(\tau),
\end{array}
\end{equation}
for which the $n_d\times n_d$-matrix $Z$ with entries $Z_{ij}=[d^i|_{[-r_0^{(\tau)}+\tau,\tau]},\tilde{\psi}^j]_\tau$ is invertible. This is possible on account of the linear independence of the sequence $\{d^i|_{[-r_0^{(\tau)}+\tau,\tau]}:1\leq i\leq n_d\}$. For any integer $1\leq j\leq n_d$ we can now choose
\begin{equation}
    \begin{array}{lcl}
         \psi^j_{(\tau)}&=&\sum\limits_{k=1}^{n_d}Z^{-1}_{kj} \tilde{\psi}^k.
    \end{array}
\end{equation}
By construction, we have $\psi^j_{(\tau)}(\tau-r_0^{(\tau)})=\psi^j_{(\tau)}(\tau)=0$ and we can compute
\begin{equation}
    \begin{array}{lcl}
         \int\limits_{\tau-r_0^{(\tau)}}^{\tau}d^i(t)^\dagger \psi^j_{(\tau)}(t)dt&=&\Big[d^i|_{[-r_0^{(\tau)}+\tau,\tau]},\psi^j\Big]_{\tau}\\[0.2cm]
         &=&\sum\limits_{k=1}^{n_d}Z^{-1}_{kj}\Big[d^i|_{[-r_0^{(\tau)}+\tau,\tau]},\tilde{\psi}^k \Big]_{\tau}\\[0.2cm]
         &=&\sum\limits_{k=1}^{n_d}Z^{-1}_{kj}Z_{ik}\\[0.2cm]
         &=&\delta_{i,j}
    \end{array}
\end{equation}
for any $1\leq i,j\leq n_d$, as desired.
\qed\\

\begin{lemma}[cf. {\cite[Pg. 13]{HJHLIN}}]\label{lemma:lin:halflineinverse}
Assume that (\asref{aannamesconstanten}{\text{H}}), (\asref{aannamesconvoluties}{\text{H}}), (\asref{aannameslimit}{\text{H}}) and (\asref{aannames0ophalflijn}{\text{H}}) are satisfied and fix  $\tau\in\R$. Then there exist bounded linear operators
\begin{equation}
    \begin{array}{lcl}
         \Lambda_{+;\tau}^{-1}:L^\infty\big([\tau,\infty);\C^M\big)&\rightarrow & W^{1,\infty}(D_\tau^\oplus;\C^M),\\[0.2cm]
         \Lambda_{-;\tau}^{-1}:L^\infty\big((-\infty,\tau];\C^M\big)&\rightarrow & W^{1,\infty}(D_\tau^\ominus;\C^M)
    \end{array}
\end{equation}
with the property that  the
identities
\begin{equation}
\begin{array}{lclcl}
[\Lambda\Lambda_{+;\tau}^{-1}f](t)
&=&
f(t), & \qquad & t \ge \tau,
\\[0.2cm]
[\Lambda\Lambda_{-;\tau}^{-1}g](s)
&=&g(s) & \qquad & t \le \tau
\end{array}
\end{equation}
hold for $f\in L^\infty\big([\tau,\infty);\C^M\big)$ 
and $g\in L^\infty\big((-\infty,\tau];\C^M\big)$.
\end{lemma}
\textit{Proof.} By symmetry, we will only construct the operator $\Lambda_{+;\tau}^{-1}$. We write $R=\Range(\Lambda)$ and $K=\mathcal{B}$. Let $R^{\perp}$ and $K^\perp$ be arbitrary complements of $R$ and $K$ respectively, so that we have
\begin{equation}
    \begin{array}{lclcl}
         W^{1,\infty}(\R;\C^M)&=&K\oplus K^\perp,
         \qquad \qquad
         L^{\infty}(\R;\C^M)&=&R\oplus R^\perp.
    \end{array}
\end{equation}
Let $\pi_{R}$ and $\pi_{R^\perp}$ denote the projections corresponding to this splitting. Then $\Lambda:K^\perp\rightarrow R$ is invertible, with a bounded inverse $\Lambda^{-1}\in\mathcal{L}(R,K^\perp)$.\\

We let $r_0^{(\tau)}>0$ and $\{\psi^i_{(\tau)}:1\leq i\leq n_d\}$ be the constant and the functions from Lemma \ref{lin:Lemma3.4} for this value of $\tau$. For $1\leq i\leq n_d$ we write $g^i_{(\tau)}\in L^\infty(\R;\C^M)$ for the function that has $g^i_{(\tau)}=\psi^i_{(\tau)}$ on $[-r_0^{(\tau)}+\tau,\tau]$, while $g^i_{(\tau)}=0$ on $(-\infty,-r_0^{(\tau)}+\tau)\cup(\tau,\infty)$. Since  we have $g^i_{(\tau)}\notin R$
for $1\leq i\leq n_d$ by Proposition \ref{prop:fredholmproperties} and these functions are linearly independent, we can explicitely choose the projection $\pi_{R^{\perp}}$ to be given by
\begin{equation}
\begin{array}{lcl}
\pi_{R^{\perp}}f&=&\sum\limits_{i=1}^{n_d}\big[\int_{-\infty}^\infty d^i(t)^\dagger f(t)dt\big]g^i_{(\tau)}.
\end{array}
\end{equation}

Upon writing $\1_{[\tau,\infty)}$ for the indicator function on $[\tau,\infty)$, we can define the inverse of $\Lambda$ on the positive half-line $[\tau,\infty)$ by
\begin{equation}
\begin{array}{lcl}
\Lambda_{+;\tau}^{-1}f&=&\Lambda^{-1}\pi_{\mathcal{R}}\1_{[\tau,\infty)}f.
\end{array}
\end{equation}
By construction, we have $g^i_{(\tau)}(t)=0$ for all $t\geq \tau$ and all $1\leq i\leq n_d$. As such, we have $[\pi_{R^{\perp}}\1_{[\tau,\infty)}f](t)=0$ for any $t\geq \tau$ and any $f\in L^\infty([\tau,\infty);\C^M)$. Hence a short computation shows that  we have 
\begin{equation}
\begin{array}{lclcl}
[\Lambda\Lambda_{+;\tau}^{-1}f](t)
&=&

[\1_{[\tau,\infty)}f](t)-\big[\pi_{\R^{\perp}}\1_{[\tau,\infty)}f\big](t)
&=&
f(t)
\end{array}
\end{equation}
for $t\geq \tau$ and $f\in L^\infty([\tau,\infty);\C^M)$,
as desired.\qed\\

For notational convenience, we write 
\begin{equation}\label{eq:def:r0}
    \begin{array}{lclcl}
         r_0&:=&%
          
         r_0^{(0)},
         \qquad \qquad
         \psi^i &:=&
         
         \psi^i_{(0)}
    \end{array}
\end{equation}
for the constant and 
functions obtained in Lemma \ref{lin:Lemma3.4} for $\tau=0$. 
As in the proof of Lemma \ref{lemma:lin:halflineinverse}, we
also write $g^i\in L^\infty(\R;\C^M)$ for the function 
\begin{equation}\label{eq:def:gi}
    \begin{array}{lcl}
         g^i(t)&=&\begin{cases}
         \psi^i(t),\enskip &t\in[-r_0,0]\\[0.2cm]
         0,\enskip&t\in(\infty,-r_0]\cup[0,\infty).\end{cases}
    \end{array}
\end{equation}
On account of the identity (\ref{eq:normalizationpsij}), we note that the function $g^i$ is continuous.\\

\textit{Proof of Proposition \ref{lemma:lin:4.3}.} 
For any $k\in\Z_{\geq 1}$ we write $\Lambda_{-;k}^{-1}$ for the inverse operator constructed in Lemma \ref{lemma:lin:halflineinverse} 
for the half-line $(-\infty,k]$,
together with 
$y^i_{(k)}=\Lambda_{-;k}^{-1}g^i$. Assumption (\asref{aannames0ophalflijn}{\text{H}}) implies that any basis of $\mathcal{B}$ remains linearly independent when restricted to the interval $(-\infty,0]$. As such, we can add an appropriate element of $\mathcal{B}$ to $y_{(k)}$ to ensure that the integral condition (\ref{eq:lin:4.11}) is satisfied. For any integer $1\leq j\leq n_d$, Lemma \ref{lemma:diffhale}
and the exponential decay of the function $d^j$ 
allow us to compute 
\begin{equation}
    \begin{array}{lcl}
         \ip{(d^j)^t,(y^i_{(k)})_t}_t&=&\int_{-\infty}^{t}d^j(s)^\dagger \big[\Lambda y^i_{(k)}\big](s)dt\\[0.2cm]
         &=&\int_{-r_0}^{0}d^j(s)^\dagger g^i(s)ds\\[0.2cm]
         &=&\delta_{ij}
    \end{array}
\end{equation}
for any $0\leq t\leq k$. We now pick a continuous function $\chi:[0,\infty)\rightarrow [0,1]$ that is zero near even integers and one near odd integers. Upon defining
\begin{equation}
    \begin{array}{lcl}
         y^i_{(\tau)}&=&\chi(2\tau)y^i_{(\lceil \tau\rceil)}+\big[1-\chi(2\tau)\big]y^i_{(\lceil \tau+\frac{1}{2}\rceil)},
    \end{array}
\end{equation}
in which $\lceil \tau\rceil$ denotes the closest integer larger or equal to $\tau$, it is easy to see that 
properties (i) through (v) are all satisfied.\qed\\

\subsection{Exponential decay}

We now focus on the exponential decay of functions in $\mathcal{Y}(\tau)$, noting that Theorem \ref{thm:expdecay}
already captures the corresponding
behaviour for functions in $\tilde{\mathcal{P}}(\tau)$.
The technical issues that we encountered during the proof of Theorem \ref{thm:expdecay} persist in this half-line setting. In particular, we need to control the behaviour of functions in $\mathcal{Y}(\tau)$ on a left half-line and a right half-line at the same time. \\

In addition, in the proof of the corresponding result in \cite{HJHLIN}, the authors 
were explicitly able
to avoid the region where $\Lambda y^i_{(\tau)}$ is non-zero
when considering the states $(y^i_{(\tau)})_\tau$.
This is of course no longer possible in our setting when $|\rmin|$ is infinite.
As such, we need to control the value of $\Lambda y^i_{(\tau)}$ in a more rigorous fashion. \\

Our first result can be see as the analogue of Lemma \ref{thm:expdecaybewijs1bewijs1},
but now the goal is to obtain estimates on
 $\Lambda y_n$ 
 for bounded sequences
 $\{y_n  \in \mathcal{Y}(\tau_n) \}$.
As a preparation, we recall
from the proof of Proposition \ref{lemma:lin:4.3}
that the identity
\begin{equation}\label{eq:lin:4.19}
    \begin{array}{lcl}
         \Lambda y&=&\sum\limits_{i=1}^{n_d}g^i\ip{(d^i)^0,y_0}_0
    \end{array}
\end{equation}
holds for $y\in\mathcal{Y}(\tau)$.
In addition, we recall the constants
$p>0,  K_{\mathrm{exp}}>0$ and $\alpha>0$ introduced in Lemma \ref{thm:expdecaybewijs0}.

\begin{lemma}\label{lemma:lin:unifLambda}
Assume that (\asref{aannamesconstanten}{\text{H}}), (\asref{aannamesconvoluties}{\text{H}}), (\asref{aannameslimit}{\text{H}}) and (\asref{aannames0ophalflijn}{\text{H}}) are satisfied and let $\{\sigma_n\}_{n\geq 1}$, $\{y_n\}_{n\geq 1}$ and $\{\tau_n\}_{n\geq 1}$ be sequences with the following properties.
\begin{enumerate}[label=(\alph*)]
    \item We have $\sigma_n > 0$ 
    for each $n$, together with $\sigma_n\uparrow \infty$.
    \item We have $y_n\in\mathcal{Y}(\tau_n)$ and $\tau_n\geq 0$ for each $n$.
    
    \item For each $n \ge 1$ we have the bound
    \begin{equation}\label{eq:lin:4.18}
    \begin{array}{lclclcl}
         |y_n(-\sigma_n+\tau_n)|&\geq &\frac{1}{2},\end{array}\end{equation}
together with the normalization
\begin{equation}
    \begin{array}{lclclcl}
    \sup\limits_{s\in(-\infty,\tau_n+p]}|y_n(s)|&=&1.
    \end{array}
\end{equation}
\item If $\rmax=\infty$, then
we have the additional bound
\begin{equation}\label{eq:lin:4.18.2}
    \begin{array}{lclclcl}
         |y_n(-\sigma_n+\tau_n)|&\geq & K_{\mathrm{exp}}e^{\alpha(-\sigma_n+\tau_n)}\sup\limits_{s\in[p+\tau_n,\infty)}e^{-\alpha s}|y_n(s)|.\end{array}\end{equation}
         \item The limit $-\sigma_n+\tau_n\rightarrow\beta_0$ holds for some $\beta_0\in\R$.
\end{enumerate}
Then the set of scalars $\{\ip{(d^i)^0,(y_n)_0}_0\}$ is bounded uniformly for $n\geq 1$ and $1\leq i\leq n_d$.
\end{lemma}
\textit{Proof.} Suppose first that $\rmax=\infty$. Fixing $n\in\Z_{\geq 1}$ and $1\leq i\leq n_d$,
we can use the bounds (\ref{eq:lin:4.18}) and (\ref{eq:lin:4.18.2}) to estimate
\begin{equation}\label{eq:lin:Haleafsch1}
    \begin{array}{lcl}
         |\ip{(d^i)^0,(y_n)_0}_0|&\leq &|d^i(0)^\dagger y_n(0)|+\Big|\sum\limits_{j=-\infty}^\infty \int_0^{r_j}d^i(s-r_j)^\dagger A_j(s-r_j)y_n(s)ds\Big|\\[0.2cm]
         &&\qquad +\Big|\int_\R\int_{0}^r d^i(s-r)^\dagger \mathcal{K}(r;s-r)y(s)dsdr\Big|\\[0.2cm]
         &\le&|d^i(0)|+S_1(i,n)+S_2(i,n)+S_3(i,n)+I_1(i,n)+I_2(i,n)+I_3(i,n),
    \end{array}
\end{equation}
in which we have defined
\begin{equation}
\begin{array}{lcl}
    S_1(i,n) & = & 
    \sum\limits_{r_j\leq p+\tau_n}\big| \int_0^{r_j}d^i(s-r_j)^\dagger A_j(s-r_j)ds \big|,
    \\[0.2cm]
    S_2(i,n) & = & 
    \sum\limits_{ r_j> p+\tau_n} \big|\int_0^{p+\tau_n}d^i(s-r_j)^\dagger A_j(s-r_j)ds\big| ,
    \\[0.2cm]
    S_3(i,n) & = & 
    \sum\limits_{ r_j> p+\tau_n} \big|\int_{p+\tau_n}^{r_j}d^i(s-r_j)^\dagger A_j(s-r_j) K_{\mathrm{exp}}^{-1}e^{\alpha(\sigma_n-\tau_n)}e^{\alpha s}ds\big| ,
\end{array}
\end{equation}
together with the corresponding
expressions $I_1(i,n)$, $I_2(i,n)$
and $I_3(i,n)$ related to the integrals involving $\mathcal{K}$.\\

We easily obtain the 
bounds 
\begin{equation}\label{eq:lin:Haleafsch2}
    \begin{array}{lcl}
         |S_1(i,n)|&\leq &\max\limits_{1\leq k\leq n_d}\sum\limits_{j=-\infty}^\infty \big|\int_0^{r_j}|d^k(s-r_j)^\dagger A_j(s-r_j)|ds\big|,\\[0.2cm]
         |S_2(i,n)|&\leq &\max\limits_{1\leq k\leq n_d}\sum\limits_{ r_j> p} \int_0^{\infty}|d^k(s-r_j)^\dagger A_j(s-r_j)|ds,\\[0.2cm]
         |I_1(i,n)|&\leq &\max\limits_{1\leq k\leq n_d}\int_{\R}\big|\int_{0}^r |d^k(s-r)^\dagger \mathcal{K}(r;s-r)|ds\big|dr,\\[0.2cm]
         |I_2(i,n)|&\leq &\max\limits_{1\leq k\leq n_d}\int_{p}^\infty\int_{0}^{\infty} |d^k(s-r)^\dagger \mathcal{K}(r;s-r)|dsdr,
    \end{array}
\end{equation}
which are uniform in $i$ and $n$.
Turning to the remaining expressions,
we pick a small $\epsilon > 0$ and assume that $n$ is large enough to
have $|\beta_0+\sigma_n-\tau_n|<\epsilon$. 
This allows us to estimate 
\begin{equation}\label{eq:lin:Haleafsch3}
    \begin{array}{lcl}
         |S_3(i,n)|&\leq& K_{\mathrm{exp}}^{-1}e^{\alpha(\sigma_n-\tau_n)}\sum\limits_{ r_j> p+\tau_n}\nrm{A_j}_\infty |r_j|e^{\alpha r_j}\nrm{d^i}_{\infty}\\[0.2cm]
         &\leq &  K_{\mathrm{exp}}^{-1}e^{\alpha(\sigma_n-\tau_n)} K_{\mathrm{exp}}e^{-2\alpha(p+\tau_n)}\nrm{d^i}_{\infty}\\[0.2cm]
         &\leq &e^{-2\alpha(p+\tau_n)}e^{\alpha\beta_0+\alpha\epsilon}\max\limits_{1\leq k\leq n_d}\nrm{d^k}_{\infty},
    \end{array}
\end{equation} 
with a corresponding bound for $I_3$.
In particular, 
both $S_3(i,n)$ and $I_3(i,n)$ converge to 0 as $n\rightarrow \infty$, so 
they can be bounded from above uniformly in $i$ and $n$.\\

In the case where $\rmax<\infty$, we can repeat this procedure with $p=\rmax$. The quantities $S_3(i,n)$ and $I_3(i,n)$ 
are identically zero in this case. \qed\\

\begin{lemma}\label{prop:lin:4.4bewijs1}
Assume that (\asref{aannamesconstanten}{\text{H}}), (\asref{aannamesconvoluties}{\text{H}}), (\asref{aannameslimit}{\text{H}}) and (\asref{aannames0ophalflijn}{\text{H}}) are satisfied and suppose that $\rmax=\infty$. Then for each $\tau\geq 0$ and each $y\in \mathcal{Y}(\tau)$ we have the bound
\begin{equation}\label{eq:lin:4.16equiv}
    \begin{array}{lcl}
         |y(t)|&\leq &\max\Big\{\frac{1}{2}\sup\limits_{s\in(-\infty,\tau+p]}|y(s)|, K_{\mathrm{exp}}\sup\limits_{s\in[p+\tau,\infty)}e^{-\alpha(s-t)}|y(s)|\Big\},\qquad t\leq -\sigma +\tau.
    \end{array}
\end{equation}
If $\rmax<\infty$, then the same statements hold with (\ref{eq:lin:4.16equiv}) replaced by
\begin{equation}\label{eq:lin:4.16equivalt}
    \begin{array}{lcl}
         |x(t)|&\leq &\frac{1}{2}\sup\limits_{s\in(-\infty,\tau+\rmax]}|x(s)|,\qquad t\leq -\sigma +\tau.
    \end{array}
\end{equation}
\end{lemma}
\textit{Proof.} Arguing by contradiction, let us consider sequences $\{\sigma_n\}_{n\geq 1}$, $\{\tau_n\}_{n\geq 1}$ and $\{y_n\}_{n\geq 1}$ 
that satisfy properties (a)-(d) in Lemma \ref{lemma:lin:unifLambda}.
If the sequence $\{-\sigma_n+\tau_n\}_{n\geq 1}$ is unbounded then we can follow the proof of Proposition \ref{thm:expdecaybewijs1} to arrive at a contradiction, since the interval $[-r_0,0]$ on which $\Lambda y_n$ might be nonzero gets `pushed out' towards $\pm\infty$.\\

Suppose therefore that
$-\sigma_n+\tau_n\rightarrow \beta_0 \in \R$,
possibly after passing to a subsequence. Combining
Lemma \ref{lemma:lin:unifLambda}
with \eqref{eq:lin:4.19}
shows that $\{\Lambda y_n\}_{n\geq 1}$ is uniformly bounded,
which allows us to apply the Ascoli-Arzela theorem to conclude that $y_n\rightarrow y_*$ uniformly on compact subsets of $\R$.
A computation similar to the proof of Lemma \ref{thm:expdecaybewijs1bewijs2} shows that $[\Lambda y_*](t)=0$ for every $t\geq 0$, since the functions $g^i$ 
vanish for these values of $t$. In particular, we must have $(y_*)_0\in Q(0)$. On account of Theorem \ref{thm:4.3equiv}, we obtain $(y_*)_0\in X^\perp(0)$, which yields 
\begin{equation}
    \begin{array}{lclcl}
         \ip{(d^i)^0,(y_*)_0}_0&=&0
    \end{array}
\end{equation}
for each $1\leq i\leq n_d$.
In view of (\ref{eq:lin:4.19}) this means that $\Lambda y_n\rightarrow 0$ uniformly on every compact subset of the real line. In particular, we must have $\Lambda y_*=0$ on the entire real line, which implies that $y_*\in\mathcal{B}$. However, this 
contradicts the integral condition (\ref{eq:lin:4.11}).
\qed\\

\begin{lemma}\label{prop:lin:4.4bewijs2}
Assume that (\asref{aannamesconstanten}{\text{H}}), (\asref{aannamesconvoluties}{\text{H}}), (\asref{aannameslimit}{\text{H}}) and (\asref{aannames0ophalflijn}{\text{H}}) are satisfied. Then there exists $C>0$ so that for all $\tau\geq 0$ and all $y\in \mathcal{Y}(\tau)$ we have the bound
\begin{equation}\label{eq:lin:4.17equiv}
    \begin{array}{lcl}
         \nrm{y}_{C_b(D_{\tau}^\ominus)}&\leq & C\nrm{y_{\tau}}_\infty.
    \end{array}
\end{equation}
\end{lemma}
\textit{Proof.} The bound (\ref{eq:lin:4.17equiv}) is in fact an equality with $C=1$ if $\rmin=-\infty$. Hence we assume that $\rmin>-\infty$. Arguing by contradiction, we can pick sequences $\{y_n\}_{n\geq 1},\{\tau_n\}_{n\geq 1}$ and $\{C_n\}_{n\geq 1}$ with $C_n\rightarrow\infty$ with $\tau_n\geq 0$ and $y_n\in \mathcal{Y}(\tau_n)$ for each $n$ in such a way that we have the identity
\begin{equation}\label{eq:lin:4.20}
    \begin{array}{lclcl}
         \nrm{y_n}_{C_b(D_{\tau_n}^\ominus)}&=&C_n\nrm{(y_n)_{\tau_n}}_\infty&=&1.
    \end{array}
\end{equation}
If the sequence $\{\tau_n\}_{n\geq 1}$ is unbounded we can follow the first half of the proof of Lemma \ref{thm:expdecaybewijs2} to arrive at a contradiction.\\

Hence we suppose that, after passing to a subsequence, we have $\tau_n\rightarrow \tau_*\geq 0$. Since the bounds on the functions $\{y_n\}_{n\geq 1}$ are stronger than those in (\ref{eq:lin:4.18}) or (\ref{eq:lin:4.18.2}), we can repeat the procedure from Lemma \ref{lemma:lin:unifLambda} to conclude that $y_n\rightarrow y_*$ uniform on compact subsets of $(-\infty,\tau_*]$. For each $n\geq 1$ we pick $s_n$ in such a way that  $|y_n(-s_n+\tau_n)|=1$. On account of Lemma \ref{prop:lin:4.4bewijs1}, the set $\{s_n\}_{n\geq 1}$ is bounded. Hence, we obtain that
\begin{equation}
    \begin{array}{lcl}
         y_*(t)&\neq &0,\enskip\text{for some }t\in \overline{(\rmin+\tau_*,\sigma+\tau_*)}.
    \end{array}
\end{equation}
In addition, we have $(y_n)_{\tau_n}\rightarrow 0$ uniformly as $n\rightarrow\infty$, so we even obtain that $y_n\rightarrow y_*$ uniformly on $D_X^++\tau_*$. If $\rmax<\infty$, we set $y_*=0$ on $(\rmax,\infty)$. In particular, we have $(y_*)_{\tau_*}=0$ and thus $[\Lambda y_*](t)=0$ for any $t\in[\tau_*,\infty)$. Moreover, we have $[\Lambda y_*](t)=0$ for any $t\in[0,\tau_*]$, since $\Lambda y_n$ is zero for these values of $t$ for each $n\in\Z_{\geq 1}$. This means that $y_*\in \mathcal{Q}(0)$ and, as before, this yields a contradiction.\qed\\

\textit{Proof of Proposition \ref{prop:lin:4.4}.} 
%
Using Lemmas \ref{prop:lin:4.4bewijs1} and \ref{prop:lin:4.4bewijs2},
we can extend the proof of
Theorem \ref{thm:expdecay}
to also include functions in $\mathcal{Y}(\tau)$. As such, for all
$\tau\geq 0$  and 
$x\in \mathcal{R}(\tau)$
we have the pointwise estimate
\begin{equation}\label{eq:unifdecRtau}
    \begin{array}{lcl}
         |x(t)|&\leq &K_{\mathrm{dec}}e^{-\alpha (\tau-t)}\nrm{x_\tau}_\infty,
         \qquad \qquad
         t \le \tau.
    \end{array}
\end{equation}
The exponential decay of $\dot{x}$ for $x\in \tilde{P}(\tau)$ follows directly from Theorem \ref{thm:expdecay}. 
Let us therefore consider an arbitrary $y \in \mathcal{Y}(\tau)$,
which satisfies the exponential bound \eqref{eq:expdecayy}.
Recalling the constant $B >0$
from Lemma \ref{lemma:lin:4.5:bewijs1},
we write
\begin{equation}
    \begin{array}{lcl}
         C&=&\sum\limits_{j=-\infty}^\infty \nrm{A_j(\cdot)}_\infty e^{\alpha |r_j|}+\sup\limits_{t\in\R}\nrm{\mathcal{K}(\cdot;t)}_{\alpha},\\[0.2cm]
         \tilde{B}&=&Be^{\alpha r_0}\sum\limits_{i=1}^{ n_d}\nrm{(d^i)^0}_\infty \nrm{g^i}_\infty .
    \end{array}
\end{equation}
Recalling the bound (\ref{eq:lin:Haleafsch4.5}) and the identity (\ref{eq:lin:4.19}), we obtain that
\begin{equation}\label{eq:lambdayafsch}
    \begin{array}{lcl}
         |\Lambda y|(t)&=& \big|\sum\limits_{i=1}^{n_d}g^i\ip{(d^i)^0,y_0}_0\big|\\[0.2cm]
         &\leq &\sum\limits_{i=1}^{n_d}\nrm{g^i}_\infty \big|\ip{(d^i)^0,y_0}_0\big|\\[0.2cm]
         &\leq &B\sum\limits_{i=1}^{n_d}\nrm{g^i}_\infty\nrm{(d^i)^0}_\infty\nrm{y_\tau}_\infty e^{-\alpha\tau}\\[0.2cm]
         &\leq &e^{-\alpha(\tau-t)}\tilde{B}\nrm{y_\tau}_\infty
    \end{array}
\end{equation}
for any $-r_0\leq t\leq 0$. Since $g^i(t)=0$ for $t\geq 0$ and $t\leq -r_0$ and since $g^i$ is continuous, we see that (\ref{eq:lambdayafsch}) is, in fact, valid for any $t\leq \tau$. As such, we immediately obtain
\begin{equation}
    \begin{array}{lcl}
         |\dot{y}(t)|&\leq & K_{\mathrm{dec}}e^{-\alpha(\tau-t)}\Big(C+\tilde{B}\Big)\nrm{y_\tau}_\infty,
         \qquad \qquad 
         t \le \tau
    \end{array}
\end{equation}
for any $\tau\geq 0$ and any $y\in\mathcal{Y}(\tau)$.\\

For the final statement we first recall the identity (\ref{eq:lin:4.19}). Since the coefficients $A_j(t)$ and $\mathcal{K}(\cdot;t)$ depend continuously on $t$ and since the functions $g^i$ are continuous, the identity (\ref{eq:lin:4.19}) yields that $\dot{x}$ is continuous on $(-\infty,\tau]$ for any $x\in\mathcal{R}(\tau)$ and any $\tau\geq 0$.
%
\qed\\

\subsection{
Projection operators}

In order to complete the proof of Theorem \ref{thm:lin:4.1},
we need to consider the behaviour of several projection operators.
In particular, we recall
the splitting
\begin{equation}\label{eq:autdecomp}
    \begin{array}{lcl}
         X&=&P(\infty)\oplus Q(\infty)
    \end{array}
\end{equation}
corresponding to the hyperbolic limiting system (\ref{autlimitsystems}) at $+\infty$, together with the notation $\overrightarrow{\Pi}_P$ and $\overrightarrow{\Pi}_Q$ for the projections onto the factors $P(\infty)$ and $Q(\infty)$.
In addition, we recall the decompositions
\begin{equation}
\label{eq:hlf:decomp:x:r:q}
\begin{array}{lclcl}
    X &= &R(\tau) \oplus Q(\tau) &= &\tilde{P}(
    \tau) \oplus Y(\tau) \oplus Q(\tau),
    \qquad 
    \qquad 
    \tau \ge 0
    \end{array}
\end{equation}
obtained above in this section
and write 
$\Pi_{\tilde{P}(\tau)}$, $\Pi_{Y(\tau)}$ and $\Pi_{Q(\tau)}$
for the corresponding projections.\\

Our first result can be seen as a supplement for the bound
(\ref{eq:4.19}) in Theorem \ref{theorem4.6}. Indeed, together
these bounds allow
the full structure of the two decompositions above
to be compared with each other
for $\tau \gg 1$.
\begin{lemma}[{cf. \cite[Lem. 4.5]{HJHLIN}}]\label{lemma:lin:4.5}
Assume that (\asref{aannamesconstanten}{\text{H}}), (\asref{aannamesconvoluties}{\text{H}}), (\asref{aannameslimit}{\text{H}}) and (\asref{aannames0ophalflijn}{\text{H}}) are satisfied. Then we have the limit
\begin{equation}
    \begin{array}{lcl}
         \lim\limits_{\tau\rightarrow\infty}\nrm{I-\overrightarrow{\Pi}_{P}|_{R(\tau)}}&=&0.
    \end{array}
\end{equation}
\end{lemma}
\textit{Proof.} If $\rmin>-\infty$ then we can follow the proof of \cite[Lem. 4.5]{HJHLIN} to obtain the desired result, so we assume that $\rmin=-\infty$. Recalling the positive constants $ K_{\mathrm{dec}}$ and $\alpha$ from Proposition \ref{prop:lin:4.4},
we write
\begin{equation}
    \begin{array}{lcl}
         C&=&\sum\limits_{j=-\infty}^\infty \nrm{A_j(\cdot)}_\infty e^{\alpha |r_j|}+\sup\limits_{t\in\R}\nrm{\mathcal{K}(\cdot;t)}_{\alpha}+\sum\limits_{j=-\infty}^\infty |A_j(\infty)| e^{\alpha |r_j|}+\nrm{\mathcal{K}(\cdot;\infty)}_\alpha.
         
    \end{array}
\end{equation}
Fix an arbitrary $\epsilon>0$ and pick $\tau_0 \gg 1$ in such a way that
the bounds
\begin{equation}
    \begin{array}{lcl}
         4K_{\mathrm{dec}}\big(1+C\big)e^{-\alpha \tau_0}&<&\frac{\epsilon}{2},\\[0.2cm]
         \sum\limits_{j=-\infty}^\infty \big|A_j(t)-A_j^+(\infty)\big|+\nrm{\mathcal{K}(\cdot;t)-\mathcal{K}(\cdot;\infty)}_{\alpha}+\nrm{\mathcal{K}(\cdot;t-\cdot)-\mathcal{K}(\cdot;\infty)}_{\alpha}&<&\frac{\epsilon}{2}
    \end{array}
\end{equation}
hold for all $t\geq \tau_0$. Recall the constant $r_0$ from (\ref{eq:def:r0}) and fix any $\tau\geq 2\tau_0+p+r_0$. 
\\

First we pick any $y\in\mathcal{R}(\tau)$ and write $\phi=y_{\tau}\in R(\tau)$. We now set out to show that 
\begin{equation}
    \begin{array}{lcl}
         \nrm{\overrightarrow{\Pi}_Q \phi}_\infty&\leq &\epsilon C'\nrm{\phi}_\infty,
    \end{array}
\end{equation}
for some constant $C'>0$. Indeed, this upper bound implies that
\begin{equation}
    \begin{array}{lclcl}
         \nrm{I-\overrightarrow{\Pi}_P|_{Y(\tau)}}&=&\nrm{\overrightarrow{\Pi}_Q|_{Y(\tau)}}
         &\leq &
          C'\epsilon,
    \end{array}
\end{equation}
which yields the desired result.\\

On account of Proposition \ref{prop:lin:4.4} we note that $y$ is continuously differentiable on $(-\infty,\tau]$, which yields that $\phi$ is continuously differentiable on $(-\infty,0]$. In addition, Proposition \ref{prop:lin:4.4} implies that
both $\phi$ and $\dot{\phi}$ decay exponentially for $t\rightarrow -\infty$.
which means that $\phi|_{(-\infty,0]}\in C^1_b\big((-\infty,0]\big)$. 
We can hence
approximate $\phi$ by functions $\{\phi_k\}_{k\geq 1}$ in $ C_b^1(D_X)$ which have $\phi_k(t)=\phi(t)$ for any $t\in (-\infty,0]$.
 These functions can be extended to $C^1$-smooth functions $\{y_k\}_{\geq 1}$, defined on $\R$, which have $(y_k)_{\tau}=\phi_k$. As such, they have $y_k(t)=y(t)$  for any $ t\leq \tau$. Due to the uniform bound on both $y$ and $\dot{y}$ from Proposition \ref{prop:lin:4.4} we can pick the functions $\{y_k\}_{k\geq 1}$ in such a way that the bound
\begin{equation}\label{eq:lin:4.27}
    \begin{array}{lcl}
         |\dot{y}_k(t)|+|y_k(t)|&\leq & 4K_{\mathrm{dec}}e^{-\alpha(\tau-t)}\Big(1+C\Big)\nrm{y_\tau}_\infty
    \end{array}
\end{equation}
holds for any $t\leq 0$ and any $k\in\Z_{\geq 1}$.\\

We now introduce the Heaviside function $H_\tau$ that acts as $H_\tau(t)=I$ if $t\geq \tau$ and zero otherwise, together with the operator
\begin{equation}
    \begin{array}{lcl}
         [\Lambda_\infty x](t)&=&\dot{x}(t)-\sum\limits_{j=-\infty}^\infty A_j(\infty)x(t+r_j)-\int_{\R}\mathcal{K}(s;\infty)x(t+s)ds .
    \end{array}
\end{equation}
Recalling the splitting (\ref{eq:autdecomp}), we observe that for any function $x\in C_b^1(\R)$ we have
\begin{equation}
    \begin{array}{lclcl}
         \big(\Lambda_\infty^{-1}H_\tau\Lambda_\infty x\big)_\tau  &\in &P(\infty),\qquad \big(\Lambda_\infty^{-1}[I-H_\tau]\Lambda_\infty x\big)_\tau &\in &Q(\infty),
    \end{array}
\end{equation}
together with
\begin{equation}
    \begin{array}{lcl}
         x_\tau&=& \big(\Lambda_\infty^{-1}H_\tau\Lambda_\infty x\big)_\tau+\big(\Lambda_\infty^{-1}[I-H_\tau]\Lambda_\infty x\big)_\tau.
    \end{array}
\end{equation}
As such, we have the representation
\begin{equation}
    \begin{array}{lcl}
         \overrightarrow{\Pi}_Qx_\tau&=&\big(\Lambda_\infty^{-1}[I-H_\tau]\Lambda_\infty x\big)_\tau
    \end{array}
\end{equation}
for any $C^1$-smooth function $x$. For any $t\in\R$ and any $k\in\Z_{\geq 1}$, we observe that
\begin{equation}
    \begin{array}{lcl}
         [\Lambda_\infty y_k](t)&=&[\Lambda y_k](t)+\sum\limits_{j=-\infty}^\infty \big[A_j(t)-A_j(\infty)\big]y_k(t+r_j)+\int_{\R}\big(\mathcal{K}(s;t)-\mathcal{K}(s;\infty\big)y_k(t+s)ds.
    \end{array}
\end{equation}
Since $[\Lambda y_k](t)=[\Lambda y](t)=0$ for $\tau_0\leq t\leq \tau$ and any $k\in\Z_{\geq 1}$,
we may hence estimate
\begin{equation}
    \begin{array}{lcl}
         \nrm{[I-H_\tau]\Lambda_\infty y_k}_\infty&\leq &\sup\limits_{t \leq \tau_0}\big[|\dot{y}_k(t)|+C\nrm{y_t}_\infty\big]+\sup\limits_{\tau_0\leq t\leq \tau}\frac{\epsilon}{2}\nrm{(y_k)_t}_\infty\\[0.2cm]
         &\leq &4K_{\mathrm{dec}}\big(1+C\big)e^{-\alpha(\tau-\tau_0)}\nrm{\phi_k}_\infty+\frac{\epsilon}{2}\nrm{\phi_k}_\infty\\[0.2cm]
         &\leq &4K_{\mathrm{dec}}\big(1+C\big)e^{-\alpha \tau_0}\nrm{\phi_k}_\infty+\frac{\epsilon}{2}\nrm{\phi_k}_\infty\\[0.2cm]
         &\leq &\epsilon\nrm{\phi_k}_\infty.
    \end{array}
\end{equation}
By the boundedness of the operator $\Lambda_\infty^{-1}$, we find that there exists a constant $C'>0$ that allows us to write
\begin{equation}
    \begin{array}{lclcl}
         \nrm{\overrightarrow{\Pi}_Q \phi_k}_\infty&=&\nrm{\big(\Lambda_\infty^{-1}[I-H_\tau]\Lambda_\infty y_k\big)_\tau}_\infty
         &\leq &\epsilon C'\nrm{\phi_k}_\infty.
    \end{array}
\end{equation}
The operator $\overrightarrow{\Pi}_Q$ is continuous, so we can take the limit $k\rightarrow \infty$ to obtain
\begin{equation}
    \begin{array}{lcl}
         \nrm{\overrightarrow{\Pi}_Q \phi}_\infty&\leq &\epsilon C'\nrm{\phi}_\infty.
    \end{array}
\end{equation}
This yields the desired bound
\begin{equation}
    \begin{array}{lclcl}
         \nrm{I-\overrightarrow{\Pi}_P|_{Y(\tau)}}&=&\nrm{\overrightarrow{\Pi}_Q|_{Y(\tau)}}
         &\leq & C'\epsilon.
    \end{array}
\end{equation}
\qed\\

\begin{lemma}[{cf. \cite[Lem. 4.6]{HJHLIN}}]\label{lemma:lin:4.6}
Assume that (\asref{aannamesconstanten}{\text{H}}), (\asref{aannamesconvoluties}{\text{H}}), (\asref{aannameslimit}{\text{H}}) and (\asref{aannames0ophalflijn}{\text{H}}) are satisfied and fix $\tau_0 \ge 0$. 
Then we have the limits
\begin{equation}
    \begin{array}{lcll}
         \nrm{[I-\Pi_{\tilde{P}(\tau_0)}]|_{\tilde{P}(\tau)}} &\rightarrow &0\enskip &\text{as }\tau\rightarrow\tau_0,\\[0.4cm]
         \nrm{[I-\Pi_{Y(\tau_0)}]|_{Y(\tau)}} &\rightarrow &0\enskip &\text{as }\tau\rightarrow\tau_0,\\[0.4cm]
         \nrm{[I-\Pi_{Q(\tau_0)}]|_{Q(\tau)}} &\rightarrow &0\enskip &\text{as }\tau\rightarrow\tau_0.
    \end{array}
\end{equation}
\end{lemma}
\textit{Proof.} The first and the third limit follow from Theorem \ref{theorem4.6}. The second limit follows from the finite dimensionality of the spaces $Y$ and from item (iii) of Proposition \ref{lemma:lin:4.3}.\qed\\

\begin{lemma}[{cf. \cite[Lem. 4.7]{HJHLIN}}]\label{lemma:lin:4.7}
Assume that (\asref{aannamesconstanten}{\text{H}}), (\asref{aannamesconvoluties}{\text{H}}), (\asref{aannameslimit}{\text{H}}) and (\asref{aannames0ophalflijn}{\text{H}}) are satisfied. Then the projections $\Pi_{Q(\tau)}$ from Lemma \ref{lemma:lin:4.6} can be uniformly bounded for all $\tau\geq 0$. 
\end{lemma}
\textit{Proof.} The proof is identical to that of \cite[Lem. 4.7]{HJHLIN} and, as such, will be omitted. It uses Proposition \ref{compactoperator}, together with Lemmas \ref{lemma:lin:4.5} and \ref{lemma:lin:4.6}.\qed\\

\begin{corollary}[{cf. \cite[Cor. 4.8]{HJHLIN}}]\label{cor:lin:4.8}
Assume that (\asref{aannamesconstanten}{\text{H}}), (\asref{aannamesconvoluties}{\text{H}}), (\asref{aannameslimit}{\text{H}}) and (\asref{aannames0ophalflijn}{\text{H}}) are satisfied. 
Then the 
projections $\Pi_{R(\tau)}$ and $\Pi_{Q(\tau)}$ 
corresponding to the first
splitting in \eqref{eq:hlf:decomp:x:r:q}
depend continuously on $\tau\in\R_{\geq 0}$. In addition, we have the limits
\begin{equation}
    \begin{array}{lclcl}
         \lim\limits_{\tau\rightarrow\infty}\nrm{\Pi_{Q(\tau)}-\overrightarrow{\Pi}_Q}&=&0,\qquad \qquad \lim\limits_{\tau\rightarrow\infty}\nrm{\Pi_{R(\tau)}-\overrightarrow{\Pi}_P}&=&0.
    \end{array}
\end{equation}
\end{corollary}
\textit{Proof.} The proof is identical to that of \cite[Cor. 4.8]{HJHLIN} and, as such, will be omitted. It uses Lemmas \ref{lemma:lin:4.5} and \ref{lemma:lin:4.7}.\qed\\

\textit{Proof of Theorem \ref{thm:lin:4.1}.} Upon defining the space $R(\tau)$ by (\ref{eq:lin:4.14}), the exponential decay rates follow from Theorem \ref{thm:expdecay} and Proposition \ref{prop:lin:4.4}. The continuity of the projections follows from Corollary \ref{cor:lin:4.8},
while the uniform bounds on the projections follow from Lemma \ref{lemma:lin:4.7}.\qed\\

\section{Degeneracies and their avoidance}\label{section:Hale}

In this section we set out to prove Corollaries \ref{cor:4.7:a} and \ref{cor:4.7:b}.
In fact, our main result below formulates
alternative conditions that can be used instead of (\asref{aannames0ophalflijn}{\text{H}}) to obtain the same conclusions.
These alternatives involve the Hale inner product, which we require to be (partially) nondegenerate
in the following sense.

\begin{definition} Let $F\subset Y$ be a subset with $0\in F$ and fix $\tau\in\R$. We say that the Hale inner product is left-nondegenerate at $\tau$ for functions in $F$ if $\psi=0$ is the only function $\psi\in F$ for which $\ip{\psi,\phi}_\tau=0$ 
holds for every $\phi\in X$. 
\end{definition}

\begin{definition}
Let $E\subset X$ be a subset with $0\in E$ and fix $\tau \in \R$. We say that the Hale inner product is right-nondegenerate at $\tau$ for functions in $E$ if $\phi=0$ is the only function $\phi\in E$ for which $\ip{\psi,\phi}_\tau=0$ holds for every $\psi\in Y$.
\end{definition}

\begin{proposition}[{cf. \cite[Cor. 4.7]{MPVL}}]\label{corollary4.7:indices} Assume that (\asref{aannamesconstanten}{\text{H}}), (\asref{aannamesconvoluties}{\text{H}}) and (\asref{aannameslimit}{\text{H}}) are satisfied. Suppose furthermore that at least one of the following three conditions is satisfied.
\begin{enumerate}[label=(\alph*)]
    \item The non-triviality condition (\asref{aannames0ophalflijn}{\text{H}}) holds.
    \item\label{enum:hale} We have $|\rmin|=\rmax=\infty$ and the Hale inner product is left-nondegenerate for functions in $B^*(\tau)$ at each $\tau\in\R$.
    \item\label{enum:dubbelhale} We have $\rmin < 0 < r_{\max}$ and for each $\tau\in\R$ the Hale inner product at $\tau$ is both left-nondegenerate 
    for functions in $B^*(\tau)$
    and right-nondegenerate
    for functions in $B(\tau)$.
\end{enumerate}
Then the identities
\begin{equation}\label{eq:4.21alt}
    \begin{array}{lclclcl}
         \dim B(\tau)&=&\dim \mathcal{B},
         \qquad \qquad
         \beta(\tau)&=&\dim B^*(\tau)&=&\dim\mathcal{B}^*
    \end{array}
\end{equation}
hold for every $\tau\in\R$. Moreover, the four Fredholm indices appearing in (\ref{eq:4.13}) are independent of $\tau$
and given by (\ref{eq:4.21alt}).
In addition, the first equation in (\ref{eq:4.13}) becomes
\begin{equation}\label{eq:4.22}
    \begin{array}{lcl}
         \ind(\pi^+_{P(\tau)})+\ind(\pi^-_{Q(\tau)})&=&-M+\ind(\Lambda)
    \end{array}
\end{equation}
with $\Lambda$ as in (\ref{def:Lambda}). Finally, the spaces $P(\tau)$, $Q(\tau)$, 
$\widehat{P}(\tau)$ and $\widehat{Q}(\tau)$ 
all vary continuously with respect to $\tau$.
\end{proposition}

In {\S}\ref{subsec:assump}
we provide various structural conditions
on the system \eqref{ditishetprobleem}
that allow the conditions (a)-(c)
above to be verified. 
They turn out to be closely related, as illustrated
by the examples that we provide
in {\S}\ref{subsec:examples}. We establish our main result in {\S}\ref{sec:ver:codim:cnt}, where we 
also describe how partial results can be obtained under weaker conditions.

\subsection{Structural conditions}
\label{subsec:assump}
In order to use Proposition \ref{corollary4.7:indices}
to compute the codimension of the space $S(\tau)$ in $X$, we either need to establish the nondegeneracy of the Hale inner product or show that the nontriviality condition (\asref{aannames0ophalflijn}{\text{H}}) is satisfied. However, it is by no means clear how this can be achieved in practice
for concrete systems. Our goal
here is to describe several more-or-less explicit criteria that can be used to verify these nondegeneracy and nontriviality conditions.\\

Some of these criteria reference the adjoint of the system (\ref{ditishetprobleem}),
which is closely related to
the operator $\Lambda^*$ defined in (\ref{def:Lambda^*}). This system is given by
\begin{equation}\label{ditishetadjointprobleem}
    \begin{array}{lcl}
         \dot{y}(t)&=&-\sum\limits_{j=-\infty}^\infty A_j(t-r_j)^\dagger y(t-r_j)-\int\limits_{\R}\mathcal{K}(\xi;t-\xi)^\dagger y(t-\xi)d\xi.
    \end{array}
\end{equation}
Most of our conditions impose the following basic structural condition,
which demands that the coefficients
corresponding to large shifts are autonomous. This is valid
for many common reaction-diffusion systems such as those studied in
\cite{BatesInfRange,HJHFHNINFRANGE}.
Indeed, the large shifts usually arise from discretizations of the diffusion, which is typically linear. The nonlinear reaction terms are typically localized in space.

\begin{assumption}{B}{\text{h}}\label{aannamesatomicbasis} There exists  a constant $K_{\mathrm{const}} \in \mathbb{Z}_{\ge 1}$ together with families of diagonal matrices
\begin{equation}
\label{eq:deg:coeffs:tilde:a:k}
    \big\{ \tilde{A}_j : j \in \Z \hbox{ with } \abs{j} \ge K_{\mathrm{const}} \big\} \subset \C^{M \times M},
    \qquad
        \big\{ \tilde{\mathcal{K}}(\xi) : \xi \in \R \hbox{ with } \abs{\xi} \ge K_{\mathrm{const}} 
        \big\} \subset \C^{M \times M},
\end{equation}
so that the following structural conditions are satisfied.
\begin{itemize}
    \item[(a)]{
      We have $r_j=j$ for $j\in\Z$, which implies $r_{\min} = - \infty$ and $\rmax=\infty$.
    }
    \item[(b)]{
      We have $A_j(t) = \tilde{A}_j$
      for all $t \in \R$ whenever $|j|\geq K_{\mathrm{const}}$.
    }
    \item[(c)]{
      We have $\mathcal{K}(\xi;t)=
      \tilde{\mathcal{K}}(\xi)$
      for all $t \in \R$
      whenever $|\xi|\geq K_{\mathrm{const}}$.
    }
\end{itemize}
\end{assumption}

\begin{remark}
The assumption (\asref{aannamesatomicbasis}{\text{h}}) can be relaxed by assuming that there exists a basis for $\C^M$ on which the matrices $\tilde{A}_{j}$ for $j\leq - K_{\mathrm{const}}$ are diagonal,
together with a separate basis on which the matrices $\tilde{A}_{j}$ for $j\geq  K_{\mathrm{const}}$ are diagonal. However, for notational simplicity, we do not pursue such an approach.
\end{remark}

\begin{remark}
The condition (\asref{aannamesatomicbasis}{\text{h}}) can be relaxed to include shifts $r_j$ with $|r_j|< K_{\mathrm{const}}$ that are not equidistant. In addition, there does not need to be any limit on the number of these small shifts. However, for notational simplicity we do not pursue such a level of generality.
\end{remark}

We divide our discussion into several scenarios for the unbounded coefficients that we each discuss in turn. Our general results are formulated at the end of this subsection.\\

\subsubsection{Bounded shifts and compact support}

The methods from \cite{MPVL} can be applied almost directly when the 
nonlocal terms all have finite range,
except that we need to take care of accumulation points of the shifts. In any case,
it is straightforward to formulate the appropriate atomic condition at a point $\tau\in\R$.
\begin{assumption}{Fin}{\text{h}}\label{aannamesatomicbdd} We have 
 $\abs{\rmin } + \rmax < \infty$ and
 there is a small $\delta>0$ so that the convolution kernel $\mathcal{K}(\cdot;t)$ is supported in the interval $[\rmin+\delta,\rmax-\delta]$ for each $t\in\R$.
 In addition, neither $\rmin$ nor $\rmax$ is an accumulation point of the set of shifts $\mathcal{R}$ and
  there are unique integers $j_{\min}, j_{\max}$ that satisfy
 \begin{equation}
    \begin{array}{lclcl}
        r_{\min} &=& r_{j_{\min}},
        \qquad \qquad 
        r_{\max} = r_{j_{\max}}.
    \end{array}
 \end{equation}
 Finally, we have $\det\big(A_{j_{\mathrm{min}}}(t)\big)\neq 0$ for a dense set of $t\in[\tau+\rmin,\tau-\rmin]$,
 together with $\det\big(A_{j_{\mathrm{max}}}(t)\big)\neq 0$
 for a dense set of $t\in[\tau-\rmax,\tau+\rmax]$.
\end{assumption}

\subsubsection{Unbounded shifts and compact support}

\label{sec:deg:unb:cmp}

We here consider the case where the discrete shifts are unbounded, but the convolution kernels all have finite support. For convenience, we formulate this as an assumption.

\begin{assumption}{Sh1}{\text{h}}\label{aannamesatomic2:prep} Assumption (\asref{aannamesatomicbasis}{\text{h}}) is satisfied. In addition, 
$\mathcal{K}(\cdot;t)$ is supported in the interval $[- K_{\mathrm{const}}, K_{\mathrm{const}}]$ for each $t\in\R$. 
\end{assumption}

Our approach here exploits the functional analytic framework of cyclic vectors for the backward shift operator on $\ell^2$, which was first described in \cite{DOUGSHAPSHIELDS}. 
This framework allows us to find sufficient conditions under which
the non-triviality condition (\asref{aannames0ophalflijn}{\text{H}})
holds and the Hale inner product is nondegenerate for exponentially decaying functions.
Reversely, we also provide a condition that guarantees the Hale inner product to be degenerate, even for exponentially decaying functions; see Proposition \ref{prop:4.16:infrangereverse} below. \\

Let us first collect the necessary terminology. We consider the backward shift operator $S$ on the sequence space $\ell^2(\N_0;\C)$, defined by
\begin{equation}
\label{eq:deg:def:bck:shift}
    \begin{array}{lclcl}
         S:\ell^2(\N_0;\C)&\rightarrow &\ell^2(\N_0;\C),
         \qquad (a_n)_{n\geq 0}&\mapsto &(a_n)_{n\geq 1}.
    \end{array}
\end{equation}
We call a vector $a=(a_n)_{n\geq 0}\in \ell^2(\N_0;\C)$ \textit{cyclic} if the 
span of the set $\{S^N a : N\geq 0\}$ is dense in $\ell^2(\N_0;\C)$. Our main condition here demands that 
the diagonal elements of the
matrices $\tilde{A}_j$ can be used to form such cyclic sequences. Our first result shows that this is in fact essential
for the nondegeneracy of the Hale inner product.

\begin{assumption}{Sh2}{\text{h}}\label{aannamesatomic2}  Upon writing $j_n = K_{\mathrm{const}} + n$ together with
\begin{equation}
\label{eq:deg:def:alpha:beta:k}
\begin{array}{lclcl}
    \alpha^{(k)}&=&
    \big(\tilde{A}_{-j_n}^{(k,k)}
    \big)_{n \ge 0} 
    \subset \ell^2(\mathbb{N}_0 ; \mathbb{C})
    \qquad \qquad
    \beta^{(k)}&=&
    \big(\tilde{A}_{j_n}^{(k,k)}\big)_{n \ge 0}
    \subset \ell^2(\mathbb{N}_0 ; \mathbb{C}),
\end{array}
\end{equation}
    the sequences $\alpha^{(k)}$ and $\beta^{(k)}$ are cyclic for the backwards shift operator for any $1\leq k\leq M$. 
\end{assumption}

\begin{proposition}[{see \S \ref{subsec:hale}}]\label{prop:4.16:infrangereverse}
Assume that (\asref{aannamesconstanten}{\text{H}}), (\asref{aannamesconvoluties}{\text{H}}) and (\asref{aannameslimit}{\text{H}})
and (\asref{aannamesatomic2:prep}{\text{h}})
are all satisfied. If the cyclicity condition (\asref{aannamesatomic2}{\text{h}}) is not satisfied, then there exists a non-zero function $\psi\in Y$ that decays exponentially
and satisfies $\ip{\psi,\phi}_\tau=0$ for every $\phi\in X$ and each $\tau\in\R$.
\end{proposition}

For the backward shift operator on $\ell^2(\N_0;\C)$, the criterion for an exponentially decaying sequence to be cyclic can be made explicit; see
{\S}\ref{sec:ver:cyc}.
This allows us to formulate two results that can be used to verify (\asref{aannamesatomic2}{\text{h}}).\\

\begin{lemma}[{see \S \ref{sec:ver:cyc}}]\label{lem:cyclicshifts1} 
Assume that (\asref{aannamesconstanten}{\text{H}}), (\asref{aannamesconvoluties}{\text{H}}) and (\asref{aannameslimit}{\text{H}})
and (\asref{aannamesatomic2:prep}{\text{h}})
are all satisfied. Consider the functions $f^{(k)}$ and $g^{(k)}$ that are defined on their natural domain by
\begin{equation}
    \begin{array}{lclcl}
         f^{(k)}(z)&=&\sum\limits_{j=K_{\mathrm{const}}}^\infty \tilde{A}_{-j}^{(k,k)}z^j,\qquad g^{(k)}(z)&=&\sum\limits_{j=K_{\mathrm{const}}}^\infty \tilde{A}_{j}^{(k,k)}z^j.
    \end{array}
\end{equation}
Then the cyclicity condition (\asref{aannamesatomic2}{\text{h}}) is satisfied if and only if the functions $f^{(k)}$ and $g^{(k)}$ are not rational functions for any $1\leq k\leq M$.
\end{lemma}

\begin{lemma}[{see \S \ref{sec:ver:cyc}}]\label{lem:cyclicshifts2} 
Assume that (\asref{aannamesconstanten}{\text{H}}), (\asref{aannamesconvoluties}{\text{H}}) and (\asref{aannameslimit}{\text{H}})
and (\asref{aannamesatomic2:prep}{\text{h}})
are all satisfied
and consider the sequences $\alpha^{(k)}$ and $\beta^{(k)}$ defined in
\eqref{eq:deg:def:alpha:beta:k}.
Then the sets
$\{S^N\alpha^{(k)}:N\geq 0\}$ and $\{S^N\beta^{(k)}:N\geq 0\}$ are both infinite dimensional for each $1\leq k\leq M$ if and only if the cyclicity condition (\asref{aannamesatomic2}{\text{h}}) is satisfied.
\end{lemma}

\subsubsection{Bounded shifts, unbounded support}
\label{sec:deg:bnd:unbnd:kernel}

We now consider the reverse of the setting discussed in {\S}\ref{sec:deg:unb:cmp}. In particular,
we assume that the discrete shifts are bounded.
\begin{assumption}{Cyc1}{\text{h}}\label{aannamesatomicconv2:prep} Assumption (\asref{aannamesatomicbasis}{\text{h}}) is satisfied,
with 
$\tilde{A}_j =0 $ whenever $\abs{j} \ge K_{\mathrm{const}}$.
\end{assumption}

In this case, one is interested in the translation semigroup $\{S_t\}_{t\geq 0}$ on the space $L^1$, which acts as
\begin{equation}\label{def:semigroup}
    \begin{array}{lcl}
         (S_tf)(s)=f(s+t)
    \end{array}
\end{equation}
for $f\in L^1\big([0,\infty);\C\big)$. A function $f\in L^1\big([0,\infty);\C\big)$ is said to be \textit{cyclic} for the translation semigroup if $\span\{S_t f:t\geq 0\}$ is dense in $ L^1\big([0,\infty);\C\big)$. We  impose the following counterpart to (\asref{aannamesatomic2}{\text{h}}),
which will allow us to establish
(\asref{aannames0ophalflijn}{\text{H}})
together with the
 nondegeneracy of the Hale inner product for bounded functions. 

\begin{assumption}{Cyc2}{\text{h}}\label{aannamesatomicconv2} 
For any $1\leq k\leq M$,
the functions
\begin{equation}
\begin{array}{lclcl}
f^{(k)}(s)&=&\tilde{\mathcal{K}}( K_{\mathrm{const}}+s)^{(k,k)},
\qquad \qquad
g^{(k)}(s)&=&\tilde{\mathcal{K}}(- K_{\mathrm{const}}-s)^{(k,k)}
\end{array}
\end{equation}
are cyclic for the translation semigroup on $ L^1\big([0,\infty);\C\big)$.
\end{assumption}


It is well-known that there exist kernels that satisfy (\asref{aannamesatomicconv2}{\text{h}})
and (\asref{aannamesconvoluties}{\text{H}}), see Lemma \ref{prop:oldresultcyclictranslation} below.  In addition, translates of such kernels remain cyclic.
However, we are unaware of any criterion to explicitely characterize them. This prevents us from 
formulating a result analogous to Lemma \ref{lem:cyclicshifts1}.
\\

\subsubsection{Positive-definite coefficients}

Our final scenario requires information
on the sign of the coefficient functions \eqref{eq:deg:coeffs:tilde:a:k} and the kernel elements in $\mathcal{B}^*$. Such information can typically be obtained by applying Krein-Rutman type arguments, see for example \cite{MPB,HJHNEGDIF,CHENGUOWU2008}. In each of these examples the kernels $\mathcal{B}$ and $\mathcal{B}^*$ are at most one-dimensional. Notice that our main condition here is weaker than the requirements formulated in  Proposition \ref{prop:4.16:positiveweaker}. For convenience we split the conditions on the coefficients and the kernels into separate assumptions.

\begin{assumption}{Pos1}{\text{h}}\label{aannamesatomicpositivestronger}  Assumption (\asref{aannamesatomicbasis}{\text{h}}) is satisfied and the matrices
\eqref{eq:deg:coeffs:tilde:a:k} are all positive semidefinite. Finally, at least one of the following
two conditions holds.
\begin{enumerate}[label=(\roman*)]
\item[(a)] For each $m\geq  K_{\mathrm{const}}$ there exist $i\geq m$ and $j\leq -m$ for which the matrices $\tilde{A}_{i}$ and $\tilde{A}_{j}$ are positive definite.
\item[(b)] The map $s\mapsto\tilde{\mathcal{K}}(s)$ is continuous on $(-\infty,- K_{\mathrm{const}}]\cup[ K_{\mathrm{const}},\infty)$. In addition, for each $m\geq  K_{\mathrm{const}}$ there exists $s\geq m$ and $r\leq -m$ for which the matrices $\tilde{\mathcal{K}}(s)$ and $\tilde{\mathcal{K}}(r)$ are positive definite.
\end{enumerate}
\end{assumption}

\begin{assumption}{Pos2}{\text{h}}\label{aannamesatomicpositivestronger2}  The adjoint kernel satisfies $\mathcal{B}^* = \{0 \}$
or $\mathcal{B}^* = \mathrm{span}\{ b \}$
for some nonnegative function $b$.  
\end{assumption}

In Proposition \ref{assumptionimplications} below, we show that the non-triviality condition (\asref{aannames0ophalflijn}{\text{H}}) is satisfied if 
(\asref{aannamesatomicpositivestronger}{\text{h}}) holds, while (\asref{aannamesatomicpositivestronger2}{\text{h}}) holds both for  the system (\ref{ditishetprobleem}) 
as well as
its adjoint
(\ref{ditishetadjointprobleem}). 
On the other hand, the left-nondegeneracy of the Hale inner product follows from the positivity condition (\asref{aannamesatomicpositivestronger}{\text{h}}) without any additional assumptions on $\mathcal{B}$ or $\mathcal{B}^*$.

\subsubsection{Summary of results}

Our main results for this subsection can now be formulated as follows.
\begin{proposition}[{see \S \ref{subsec:hale}}]\label{prop:4.16volledig} Assume that (\asref{aannamesconstanten}{\text{H}}), (\asref{aannamesconvoluties}{\text{H}}) and (\asref{aannameslimit}{\text{H}}) are satisfied. Then we have the following implications.
\begin{enumerate}[label=(\roman*)]
    \item\label{enum:halebdd} If the atomic condition (\asref{aannamesatomicbdd}{\text{h}}) is satisfied at some point $\tau\in\R$, then the Hale inner product $\ip{\cdot,\cdot}_\tau$ is left-nondegenerate at $\tau$ for functions in $Y$ and right-nondegenerate at $\tau$ for functions in $X$.
    \item\label{enum:halecycsh} If the cyclicity conditions
    (\asref{aannamesatomic2:prep}{\text{h}})
    and (\asref{aannamesatomic2}{\text{h}}) are satisfied, then 
    at each $\tau \in \R$ the Hale inner product $\ip{\cdot,\cdot}_\tau$ 
    is left-nondegenerate and right-nondegenerate for exponentially decaying functions.
    \item\label{enum:halecycconv} If the cyclicity conditions
    (\asref{aannamesatomicconv2:prep}{\text{h}}) and
    (\asref{aannamesatomicconv2}{\text{h}}) are satisfied, then at each $\tau \in \R$ the Hale inner product $\ip{\cdot,\cdot}_\tau$ 
    is left-nondegenerate for functions in $Y$ and right-nondegenerate for functions in $X$.
    \item\label{enum:halepos} If the positivity condition (\asref{aannamesatomicpositivestronger}{\text{h}}) is satisfied, then 
    at each $\tau \in \R$ the Hale inner product $\ip{\cdot,\cdot}_\tau$ 
    is left-nondegenerate and right-nondegenerate for nonnegative functions.
\end{enumerate}
In each of the cases \ref{enum:halebdd}-\ref{enum:halecycconv},
the quantity in (\ref{eq:4.12}) satisfies $\beta(\tau)=\dim B^*(\tau)$.
This also holds for case \ref{enum:halepos} provided that  the positivity condition (\asref{aannamesatomicpositivestronger2}{\text{h}}) is satisfied.
\end{proposition}
\begin{proposition}[{see \S \ref{subsec:HAker}}]\label{assumptionimplications} Assume that (\asref{aannamesconstanten}{\text{H}}), (\asref{aannamesconvoluties}{\text{H}}) and (\asref{aannameslimit}{\text{H}}) are satisfied. Then we have the following implications.
\begin{enumerate}[label=(\roman*)]
    \item If the atomic condition (\asref{aannamesatomicbdd}{\text{h}}) is satisfied at each $\tau\in\R$,
    then the non-triviality condition (\asref{aannames0ophalflijn}{\text{H}}) is satisfied for the system (\ref{ditishetprobleem}).
    \item\label{enum:cyclic} If the cyclicity conditions
    (\asref{aannamesatomic2:prep}{\text{h}})
    and (\asref{aannamesatomic2}{\text{h}}) are satisfied,
    then the non-triviality condition (\asref{aannames0ophalflijn}{\text{H}}) is satisfied for the system (\ref{ditishetprobleem}).
    \item If the cyclicity conditions
    (\asref{aannamesatomicconv2:prep}{\text{h}}) and
    (\asref{aannamesatomicconv2}{\text{h}}) are satisfied, then the non-triviality condition (\asref{aannames0ophalflijn}{\text{H}}) is satisfied for the system (\ref{ditishetprobleem}).
    \item If the positivity condition (\asref{aannamesatomicpositivestronger}{\text{h}}) is satisfied and 
    (\asref{aannamesatomicpositivestronger2}{\text{h}}) 
    holds both for 
    (\ref{ditishetprobleem}) and its adjoint (\ref{ditishetadjointprobleem}), then the non-triviality condition (\asref{aannames0ophalflijn}{\text{H}}) is satisfied for the system (\ref{ditishetprobleem}).
\end{enumerate}
\end{proposition}

Note that the non-triviality condition (\asref{aannames0ophalflijn}{\text{H}}) does not directly imply that the the Hale inner product is nondegenerate in some form. Instead, it enables us construct an explicit complement to the space $S(\tau)$. In particular, the nondegeneracy of the Hale inner product is useful, but not necessary to compute the codimension $\beta(\tau)$.\\

\subsection{Examples}\label{subsec:examples}

In order to illustrate the results above, we consider
the infinite-range nonlinear MFDE
\begin{equation}\label{eq:nagumomfde}
    \begin{array}{lcl}
         \dot{u}(t)&=&\sum\limits_{k=1}^\infty \gamma_k[u(t+k)+u(t-k)-2u(t)]+\int\limits_{0}^\infty \theta(\xi) [u(t+\xi) + u(t - \xi)  - 2 u(t)] d\xi \\[0.2cm]
          & & \qquad  
            +g\big(u(t);a\big),
    \end{array}
\end{equation}
in which the nonlinearity $g$ is given by the cubic nonlinearity
\begin{equation}
    g(u;a)=u(1-u)(u-a), \qquad \qquad a \in (0,1),
\end{equation}
while the sequence $\gamma$ and the function $\theta$ decay exponentially.
This MFDE can be interpreted as the travelling wave equation for a nonlocal version of the Nagumo PDE. One is typically interested in the front solutions, which satisfy the limits
\begin{equation}
\label{eq:deg:limits:front}
    \begin{array}{lclcl}
         \lim\limits_{t\rightarrow-\infty}\overline{u}(t)&=&0,\qquad \lim\limits_{t\rightarrow\infty}\overline{u}(t)&=&1.
    \end{array}
\end{equation}
Results concerning the existence of such these solutions in a variety of settings can be found in \cite{VL22,BatesInfRange,VL28,MPB}. For our purposes here, we will simply assume such a solution exists
and consider the associated linearization
of (\ref{eq:nagumomfde}), which is given by
\begin{equation}\label{eq:nagumomfdelinearization}
    \begin{array}{lcl}
         \dot{u}(t)&=&\sum\limits_{k=1}^\infty \gamma_k[u(t+k)+u(t-k)-2u(t)]+\int\limits_{0}^\infty \theta(\xi)\big[ u(t+\xi) + u(t -  \xi) - 2 u(t) \big] d\xi
         \\[0.2cm]
         & & \qquad
         +g_u(\overline{u}(t);a)u(t).
    \end{array}
\end{equation}
We remark that a simple differentiation automatically yields $\frac{d}{dt}\overline{u} \in \mathcal{B}$.\\ 

In this setting, the Hale inner product
is given by 
\begin{equation}
    \begin{array}{lcl}
          \ip{\psi,\phi}_\tau&=&\overline{\psi(0)} \phi(0)+\sum\limits_{k=1}^\infty \int\limits_{-k}^0 \overline{\psi(s+k)} \gamma_{k}\phi(s)ds-\sum\limits_{k=1}^\infty \int\limits_0^{k} \overline{\psi(s-k)} \gamma_{|k|}\phi(s)ds\\[0.2cm]
         &&\qquad -\int\limits_{\R}\int\limits_0^r \overline{\psi(s-r)} \theta(|r|)\phi(s)dsdr,
    \end{array}
\end{equation}
which is independent of $\tau$ and the function $\overline{u}$.
With the exception of (\asref{aannamesatomicpositivestronger2}{\text{h}}),
we can hence investigate
the
validity of our assumptions
and the nondegeneracy of the Hale inner product without any knowledge regarding
the wave $\overline{u}$ besides the limits \eqref{eq:deg:limits:front}.\\

For example,
we note that (\asref{aannamesatomicbasis}{\text{h}}) is automatically satisfied with $K_{\mathrm{const}} = 1$ and
\begin{equation}
\begin{array}{lclcl}
    \tilde{A}_j &=& \gamma_{\abs{j}}
    \qquad \qquad
    \tilde{\mathcal{K}}(\xi)& =& \theta(\abs{\xi})
\end{array}
\end{equation}
for $\abs{j} \ge 1$ and $\xi \neq 0$. In addition, we have
\begin{equation}
\begin{array}{lclcl}
    A_0(t)& =& -2 \sum\limits_{k=1}^\infty \gamma_k - 2 \int_0^\infty \theta(\xi) d \xi + g_u(\overline{u}(t) ; a ).
\end{array}
\end{equation}
In particular, it is clear that (\asref{aannamesconstanten}{\text{H}}) and (\asref{aannamesconvoluties}{\text{H}}) hold. However, one needs
additional information on the coefficients
in order to verify the hyperbolicity assumption (\asref{aannameslimit}{\text{H}}).\\

We consider various choices for $\gamma$ and $\theta$ in our discussion below. In each case we are able to distinguish whether or not the Hale inner product is degenerate.  For each of the two degenerate cases, 
we construct an explicit non-trivial function $ \psi\in Y$ for which $\ip{\psi,\phi}_\tau=0$ for all $\phi\in X$ and all $\tau\in\R$. However, we emphasize again that this does not prevent us from showing that (\asref{aannames0ophalflijn}{\text{H}}) holds.

\subsubsection{Positive coefficients}\label{subsubsec:posnagumo}

Consider the system (\ref{eq:nagumomfdelinearization}) and suppose that the coefficients $\{\gamma_k\}_{k\geq 1}$ and the convolution kernels $\theta(\xi)$ are positive. The bistablity of the nonlinearity $g$ then allows us to conclude that
the hyperbolicity condition (\asref{aannameslimit}{\text{H}}) is satisfied. In addition,
(\asref{aannamesatomicpositivestronger}{\text{h}})
holds and hence
the Hale inner product is nondegenerate for nonnegative functions.\\

These positivity conditions
imply that a comparison principle holds 
for (\ref{eq:nagumomfdelinearization}). In such a setting, one can typically derive that the kernels $\mathcal{B}$ and $\mathcal{B}^*$ are both 
one-dimensional and spanned by a strictly positive function. For example,
the wave $\overline{u}$ is typically monotonically increasing and the associated
derivative $\frac{d}{dt} \overline{u}$ spans $\mathcal{B}$ and is strictly positive.
Results of this type have been proven in various settings, see for example \cite{bates1997traveling,BatesDiscConv,CHEN1997}. In each case, the system (\ref{eq:nagumomfdelinearization}) together with its adjoint \eqref{ditishetadjointprobleem} satisfy 
(\asref{aannamesatomicpositivestronger2}{\text{h}}). In particular, the nontriviality condition (\asref{aannames0ophalflijn}{\text{H}}) holds.

\subsubsection{Non-cyclic shift coefficients}
Consider the system (\ref{eq:nagumomfdelinearization}) with $\theta(\xi)=0$ for each $t\in\R$ and $\gamma_k=e^{-k}$ for $k\geq 1$. This system satisfies (\asref{aannamesatomic2:prep}{\text{h}}). Since the coefficients $\{\gamma_k\}_{k\geq 1}$ are positive, the results from \S \ref{subsubsec:posnagumo} show that (\asref{aannameslimit}{\text{H}}) is satisfied and that the Hale inner product for the system (\ref{eq:nagumomfdelinearization})
is
nondegenerate for nonnegative functions. \\

However, it is easy to see that
\begin{equation}
\begin{array}{lcl}
    \sum\limits_{k\geq 1} \gamma_{k}z^k&=&\frac{z}{e-z},
\end{array}
\end{equation}
 which is a rational function. Hence, this system does not satisfy (\asref{aannamesatomic2}{\text{h}}) on account of Lemma \ref{lem:cyclicshifts1}. Alternatively, letting $S$ denote the backwards shift operator on $\ell^2(\N_0;\C)$, the sequence $\alpha=(\gamma_k)_{k\geq 1}$ satisfies $S^N\alpha=e^{-N}\alpha$ for any $N\geq 0$. In particular, the set $\span\{S^N\alpha:N\geq 0\}$ is one-dimensional, which in view of
 Lemma \ref{lem:cyclicshifts2}
 again shows that (\asref{aannamesatomic2}{\text{h}})
 is not satisfied.
 In particular, Proposition \ref{prop:4.16:infrangereverse} implies that the Hale inner product  
 is not nondegenerate
 for all exponentially decaying functions.
  \\

To make this more explicit, 
we consider the continuous, bounded function $\psi:\R\rightarrow \R$
that has
\begin{equation}
    \psi(s)=0 \hbox{ for } s\leq 1,
    \qquad
    \psi\Big(\frac{3}{2}\Big)=1,
    \qquad
    \psi\Big(\frac{5}{2}\Big)=-e,
    \qquad
    \psi(s)=0 \hbox{ for } s\geq 3
\end{equation}
and is linear in the missing segments.
This choice 
is motivated by the fact that
\begin{equation}
\beta=(1,-e,0,0,...)\in\ell^2(\mathbb{N}_0;\C)
\end{equation}
is perpendicular to the set $\span\{S^N\alpha:N\geq 0\}$
and ensures that  
\begin{equation}\label{eq:chardefpsispecific}
    \begin{array}{lclcl}
         \sum\limits_{k=m}^\infty \psi\big(\tilde{s}+k+ 1-m\big) \gamma_k&=&\tilde{s} \big(e^{-m}-e\cdot e^{-(m+1)}\big)
         &=&0, 
    \end{array}
\end{equation}
for any $m\in\Z_{\geq 1}$ and any $\tilde{s}\in[0,1)$. For an arbitrary $s \le 0$
we make the decomposition
\begin{equation}
\begin{array}{lcl}
    s &=& \tilde{s} + 1 - m
\end{array}
\end{equation}
for some integer $m \ge 1$ and $\tilde{s} \in [0,1)$. Applying \eqref{eq:chardefpsispecific},
we now compute
\begin{equation}
\begin{array}{lclcl}
    \sum\limits_{k \ge 1 - s} \psi(s +k) \gamma_k
    &=& \sum\limits_{k \ge m - \tilde{s}} \psi(\tilde{s} + k + 1 - m ) \gamma_k
    &=& 0
\end{array}
\end{equation}
since the final sum in fact ranges over $k \ge m$.\\

Since $\psi(s)=0$ for $s\leq 1$, the Hale inner product reduces to
\begin{equation}\label{eq:chardefpsi:expl}
    \begin{array}{lcl}
          \ip{\psi,\phi}_\tau&=&\sum\limits_{k=1}^\infty \int\limits_{-k}^0 \psi(s+k) \gamma_{k}\phi(s)ds\\[0.2cm]
          &=&\sum\limits_{k=1}^\infty \int\limits_{-k+1}^0 \psi(s+k) \gamma_{k}\phi(s)ds,
         
    \end{array}
\end{equation}
for $\phi\in C_b(\R)$ and $\tau\in\R$. The dominated convergence theorem allows us to interchange the sum and the infinite integral, which yields
\begin{equation}
    \begin{array}{lclcl}
         \ip{\psi,\phi}_{\tau}
         &=&\int\limits_{-\infty}^0  \sum\limits_{ k\geq 1-s}  \psi(s+k) \gamma_k\phi(s)ds
         &=&0,
    \end{array}
\end{equation}
for any $\phi\in C_b(\R)$ and $\tau\in\R$.
Since $\gamma_k > 0$ for any $k\in\Z_{\geq 1}$, this example shows that
a naive generalization of the atomic condition (\asref{aannamesatomicbdd}{\text{h}}) is not sufficient
to establish the nondegeneracy of the Hale inner product.

\subsubsection{Non-cyclic convolution kernel}

Consider the system (\ref{eq:nagumomfdelinearization}) with $\theta(\xi)=\exp(-\xi)$ and $\gamma_k=0$ for $k\geq 1$. This system satisfies (\asref{aannamesatomicconv2:prep}{\text{h}}). Since the kernel $\theta$ is positive, the results from \S \ref{subsubsec:posnagumo} again show that (\asref{aannameslimit}{\text{H}}) is satisfied and that the Hale inner product for the system (\ref{eq:nagumomfdelinearization})
is nondegenerate for nonnegative functions. \\

However, the identity
\begin{equation}
\begin{array}{lcl}
    \theta(t+\xi)&=&\exp(-t)\theta(\xi),
    \qquad \qquad
    (\xi, t) \in R_{\ge 0}^2
\end{array}
\end{equation}
directly implies that
 $\span\{\theta(\cdot+t):t\geq 0\}$ is one dimensional in $L^1\big([0,\infty);\C\big)$.
 In particular, 
 the cyclicity condition (\asref{aannamesatomicconv2}{\text{h}}) fails to be satisfied. While we cannot appeal to a general result here, we can show by hand that the Hale inner product is degenerate for an exponentially decaying function.\\

To this end, we consider the bounded, continuous function $\psi:\R\rightarrow \R$
that has
\begin{equation}
    \psi(s)=0 \hbox{ for } s\leq 0,
    \qquad
    \psi(1)=-1,
    \qquad
    \psi(2)=0,
    \qquad
    \psi(3)=e^2,
    \qquad
    \psi(s)=0 \hbox{ for } s\geq 4
\end{equation}
and is linear in the missing segments.
By construction, the identity
\begin{equation}\label{eq:chardefpsispecific2}
\begin{array}{lclcl}
\int\limits_{0}^\infty \psi(r)\theta(r-s)dr&=&\exp(s)\int\limits_{0}^\infty \psi(r)\theta(r)dr
&=&0
\end{array}
\end{equation}
holds for any $s\leq 0$. 
For any $\phi\in C_b(\R)$ we can again use the dominated convergence theorem to compute 
\begin{equation}
    \begin{array}{lcl}
         \ip{\psi,\phi}_{\tau}&=&\psi(0)\phi(0)-\int\limits_{\R} \int\limits_0^r \psi(s-r)^\dagger \theta(|r|)\phi(s)dsdr\\[0.2cm]
         &=&-\int\limits_{0}^{\infty} \int\limits_{-r}^0 \psi(s+r)^\dagger \theta(r)\phi(s)dsdr\\[0.2cm]
         &=&-\int\limits_{-\infty}^0 \int\limits_0^{\infty} \1_{\{s\in[-r,0]\}} \psi(s+r)^\dagger \theta(r)\phi(s)drds\\[0.2cm]
         &=&-\int\limits_{-\infty}^0 \int\limits_{-s}^{\infty}  \psi(s+r)^\dagger \theta(r)\phi(s)drds\\[0.2cm]
         &=&-\int\limits_{-\infty}^0 \int\limits_0^{\infty}  \psi(r)^\dagger \theta(r-s)\phi(s)drds\\[0.2cm]
         &=&0
    \end{array}
\end{equation}
for any $\tau\in\R$.

\subsubsection{Cyclic shifts with mixed coefficients}

For our final example, we choose $\theta =0 $ and consider a sequence $\gamma$ that admits Gaussian decay. In particular, we write
\begin{equation}
\begin{array}{lcl}
    \gamma_k&=&\frac{1}{h^2}c_k \exp(-k^2),
    \qquad \qquad h > 0,
\end{array}
\end{equation}
for some bounded sequence
$\{c_k\}_{k\geq 1}$ that can have both positive and negative elements, but must be uniformly bounded away from zero.
In particular,
(\asref{aannamesatomic2:prep}{\text{h}}) is satisfied,
but this may not hold for the positivity condition (\asref{aannamesatomicpositivestronger}{\text{h}}). In order to verify the hyperbolicity condition
(\asref{aannameslimit}{\text{H}}), it suffices to impose the restriction
\begin{equation}\begin{array}{lcl}\sum\limits_{k>0}c_k\exp(-k^2)\Big(1-\cos(kz)\Big)&>& 0, \qquad \qquad z\in(0,2\pi);
\end{array}\end{equation}
see \cite[Lem. 5.6]{HJHFHNINFRANGE}.
This can be interpreted as the statement
that the sum in  \eqref{eq:nagumomfdelinearization}
is spectrally similar to the Laplacian. \\

We now set out to establish 
the cyclicity condition (\asref{aannamesatomic2}{\text{h}}) by appealing to 
Lemma \ref{lem:cyclicshifts2}.
Recalling the backward shift operator
\eqref{eq:deg:def:bck:shift},
we consider the vector
$e = (e_n)_{n \ge 0} \in \ell^2(\N_0;\C)$ given by 
\begin{equation}
    \begin{array}{lcl}
         e_{k-1}&=& c_k\exp(-k^2),
         \qquad \qquad k \ge 1
    \end{array}
\end{equation}
and set out to show that the set
\begin{equation}
    \begin{array}{lcl}
         \mathcal{A}&:=&\span\{S^N e:N\geq 0\}
    \end{array}
\end{equation}
is an infinite dimensional subspace of $\ell^2(\N_0;\C)$. \\

Arguing by induction, we pick $\ell \geq 1$ and assume that the vectors $e,S e,...,S^{\ell-1}e$ are linearly independent. Suppose now that we have a nonzero multiplet $(\lambda_0,...,\lambda_\ell)\in \C^{\ell+1}$ for which
\begin{equation}
\label{eq:id:for:sum:lambda:S}
    \begin{array}{lcl}
         \sum\limits_{i=0}^{\ell}\lambda_i S^i e&=&0. 
    \end{array}
\end{equation}
Let $0\leq i_*< \ell$ be the smallest integer with $\lambda_{i_*} \neq 0$. Our assumption
on $c$ implies that the sequence
$\{|\frac{c_{k}}{c_{k+1}}| \}_{k \ge 1}$
is uniformly bounded away from zero, 
which implies that the quotient
\begin{equation}
\begin{array}{lcl}
    \abs{\frac{e_{k-1}}{e_{k}}}
    &= &
    \abs{\frac{c_{k}}{c_{k+1}} } \exp(2k+1)
\end{array}
\end{equation}
grows to infinity as $k\rightarrow \infty$.
In particular, by picking a sufficiently large index $K \gg 1$ we obtain the bound
\begin{equation}
    \begin{array}{lcl}
         |\sum\limits_{i=i_*+1}^\ell \lambda_i (S^i e)_K|
         &\leq &\sum\limits_{i=i_*+1}^\ell |\lambda_i e_{K+i}|\\[0.2cm]
         &<&|\lambda_{i_*} e_{K+i_*}|\\[0.2cm]
         &=&|\lambda_{i_*} (S^{i_*} e)_K|,
    \end{array}
\end{equation}
which contradicts 
the $K$-th component
of the identity
\eqref{eq:id:for:sum:lambda:S}.
In particular, Proposition \ref{prop:4.16volledig} yields that there is no exponentially decaying $\psi\in Y$ that has $\ip{\psi,\phi}_\tau=0$ for each $\phi \in X$.\\

If $h > 0$ is sufficiently small,
then  the existence of a travelling front solution for (\ref{eq:nagumomfde}) is guaranteed by \cite[Thm. 1]{BatesInfRange}.
One can subsequently use
Proposition \ref{assumptionimplications} to conclude that the non-triviality condition (\asref{aannames0ophalflijn}{\text{H}}) is satisfied.

\subsection{(Co)-dimension counting}
\label{sec:ver:codim:cnt}

The main goal of this subsection is to 
establish the identities
(\ref{eq:4.21alt}) concerning the dimensions of $B(\tau)$ and $B^*(\tau)$
and the codimension of $S(\tau)$. The remainder of the statements
in Proposition \ref{corollary4.7:indices}
follow readily from these computations,
using the main results in {\S}\ref{section:mainresults}.
We aim to use as little information as possible, providing partial results under weaker conditions.

\begin{lemma}\label{cor4.7:bewijs1} Assume that (\asref{aannamesconstanten}{\text{H}}), (\asref{aannamesconvoluties}{\text{H}}) and (\asref{aannameslimit}{\text{H}}) are satisfied. Fix $\tau\in\R$ and suppose first that the Hale inner product is left-nondegenerate at $\tau$ for functions in $B^*(\tau)$. Then the identity
\begin{equation}\label{eq:4.21.2}
    \begin{array}{lclclcl}
         \beta(\tau)&=&\dim B^*(\tau)
    \end{array}
\end{equation}
holds.
Alternatively, if the non-triviality condition (\asref{aannames0ophalflijn}{\text{H}}) is satisfied, then the identity (\ref{eq:4.21.2})
is valid for all $\tau\in\R$.
\end{lemma}
\textit{Proof.} In the first case, this follows directly from the characterisation of $S(\tau)$ given by (\ref{eq:4.11}). In the second case, the statement
for $\tau\geq 0$
follows from the direct sum decomposition (\ref{eq:lin:4.33})
and the identities in (\ref{eq:lin:4.9}). Using symmetry arguments this can be extended to $\tau<0$. \qed\\

\begin{lemma}\label{cor4.7:bewijs2} Assume that (\asref{aannamesconstanten}{\text{H}}), (\asref{aannamesconvoluties}{\text{H}}) and (\asref{aannameslimit}{\text{H}}) are satisfied. Fix $\tau\in\R$ and suppose first that any nonzero $d\in \mathcal{B}\cup\mathcal{B}^*$ does not vanish on $(-\infty,\tau]$ and does not vanish on $[\tau,\infty)$. Then we have the identities
\begin{equation}\label{eq:4.21.1}
    \begin{array}{lclclcl}
         \dim B(\tau)&=&\dim \mathcal{B},
         \qquad \qquad
         \dim B^*(\tau)&=&\dim\mathcal{B}^*.
    \end{array}
\end{equation}
In particular, if the non-triviality condition (\asref{aannames0ophalflijn}{\text{H}}) holds then (\ref{eq:4.21.1}) is valid for each $\tau\in\R$. 
\end{lemma}
\textit{Proof.} 
Since the statements
hold trivially if $|\rmin|=\rmax=\infty$ on account of Lemma \ref{lemma:elementaryproperties},
we will use symmetry to assume without loss that $r_{\max} < \infty$.
Arguing by contradiction
to establish the first identity, let us consider a non-trivial kernel element $x\in\mathcal{B}$ that has $x_\tau=0$.
If $\rmin=-\infty$, this means that $x$ vanishes identically on $D_\tau^\ominus$ and hence $(-\infty, \tau]$, violating our assumption. On the other hand, if $r_{\min} > - \infty$ we can assume without loss that $x$ does not vanish on $(r_{\max}, \infty)$. Upon introducing the new function
\begin{equation}
    \begin{array}{lcl}
       \tilde{x}(t)&=&\begin{cases}
         x(t),\qquad &t\geq \tau+\rmin,\\[0.2cm]
         0,\qquad &t<\tau+\rmin,\end{cases}
    \end{array}
\end{equation}
we see that $\tilde{x}$ is a non-trivial element of $\mathcal{B}$ that vanishes on $D_\tau^\ominus$, again violating our assumption.
The second identity in \eqref{eq:4.21.1} can be obtained in a similar fashion.
\qed\\

\begin{lemma}\label{cor4.7:bewijs3} Assume that (\asref{aannamesconstanten}{\text{H}}), (\asref{aannamesconvoluties}{\text{H}}) and (\asref{aannameslimit}{\text{H}}) are satisfied and that $\rmin<0<\rmax$. Suppose that for each $\tau\in\R$ the Hale inner product is left-nondegenerate for functions in $B^*(\tau)$. Then we have the identity
\begin{equation}\label{dim:counting1}
    \begin{array}{lcl}
         \dim B^*(\tau)&=&\dim \mathcal{B}^* 
    \end{array}
\end{equation}
for any $\tau\in\R$. Similarly, if for each $\tau\in \R$ the Hale inner product is right-nondegenerate for functions in $B(\tau)$, then the identity
\begin{equation}
    \begin{array}{lcl}
         \dim B(\tau)&=&\dim \mathcal{B}
    \end{array}
\end{equation}
holds for each $\tau\in\R$.
\end{lemma}
\textit{Proof.} Both identities follow trivially from Lemma \ref{lemma:elementaryproperties} if $|\rmin|=\rmax=\infty$. By symmetry we only consider the identity (\ref{dim:counting1}). Suppose that (\ref{dim:counting1}) fails, allowing us to pick
a non-zero $y\in \mathcal{B}^*$
that has $y^\tau=0$ for some $\tau\in\R$.
Possibly after increasing $\tau$, we may assume by symmetry that $\rmin>-\infty$ and that there exists a small $0<\epsilon<|\rmin|$ so that
\begin{equation}
    \begin{array}{lcl}
         y(\tau-\rmin+\delta)&\neq &0
    \end{array}
\end{equation}
holds for each $\delta\in (0,\epsilon)$. In particular, $0\neq y^{\tau+\epsilon} \in B^*(\tau+\epsilon)$, so by the left-nondegeneracy of the Hale inner product at $\tau+\epsilon$, we can pick $\phi\in X$ with 
\begin{equation}\label{eq:contradictionHale}
    \begin{array}{lcl}
         \ip{y^{\tau+\epsilon},\phi}_{\tau+\epsilon}&\neq &0.
    \end{array}
\end{equation}
Without loss, we can assume that $\phi$ is differentiable, allowing us to pick a differentiable function $x\in C_b(\R)$ that has $\phi=x_{\tau+\epsilon}$. On account of Lemma \ref{lemma:diffhale} we can compute
\begin{equation}
    \begin{array}{lclcl}
        \frac{d}{dt}\ip{y^t,x_t}_t&=&y^*(t)[\Lambda x](t)+[\Lambda^*y](t)x(t)&=&0
    \end{array}
\end{equation}
for any $t\in (\tau-\rmax,\tau-\rmin)$, since $y^\tau=0$ and since $y\in \mathcal{B}^*$. As such, $\ip{y^t,x_t}_t$ is constant on $(\tau-\rmax,\tau-\rmin]$. Since $y^\tau=0$, it follows that $\ip{y^\tau,x_\tau}_\tau=0$. However, this yields the identity
\begin{equation}
    \begin{array}{lclcl}
         0&=&\ip{y^{\tau+\epsilon},x_{\tau+\epsilon}}_{\tau+\epsilon}
         &=&\ip{y^{\tau+\epsilon},\phi}_{\tau+\epsilon},
    \end{array}
\end{equation}
which contradicts (\ref{eq:contradictionHale}). \qed\\

\textit{Proof of Proposition \ref{corollary4.7:indices}.}
We first aim to establish (\ref{eq:4.21alt}). If the non-triviality condition (\asref{aannames0ophalflijn}{\text{H}}) holds, this follows by combining Lemmas \ref{cor4.7:bewijs1} and Lemma \ref{cor4.7:bewijs2}. 
Alternatively, if \ref{enum:hale} holds, then (\ref{eq:4.21alt}) follows by combining 
Proposition \ref{prop:4.16volledig}
with Lemmas \ref{lemma:elementaryproperties}
and \ref{cor4.7:bewijs1}.
Finally, if \ref{enum:dubbelhale} holds, then (\ref{eq:4.21alt}) follows by combining Proposition \ref{prop:4.16volledig} with Lemmas \ref{cor4.7:bewijs1} and 
\ref{cor4.7:bewijs3}. \\

Turning to the Fredholm indices, 
we remark that the right-hand side of (\ref{eq:4.13}) is now constant in $\tau$. Since both $\ind(\pi^+_{P(\tau)})$ and $\ind(\pi^-_{Q(\tau)})$ are upper semi-continuous by Proposition \ref{theorem4.6:indices}, both these factors must be constant as well. By Theorem \ref{thm:4.3equiv} the inclusions $\widehat{P}(\tau)\subset P(\tau)$ and $\widehat{Q}(\tau)\subset Q(\tau)$ have constant codimension $\dim B(\tau)=\dim\mathcal{B}$. Hence the indices $\ind(\pi^+_{\widehat{P}(\tau)})$ and $\ind(\pi^-_{\widehat{Q}(\tau)})$ are also constant. Moreover, these four subspaces vary continuously in $\tau$. Finally, the identity (\ref{eq:4.22}) follows from (\ref{eq:4.13}) and (\ref{eq:4.21alt}), using the value of $\ind(\Lambda)$ given in Proposition \ref{prop:fredholmproperties}.\qed\\

\textit{Proof of Corollaries \ref{cor:4.7:a} and \ref{cor:4.7:b}} These results follow directly from Proposition \ref{corollary4.7:indices}.\qed\\

\subsection{Cyclic coefficients}
\label{sec:ver:cyc}

In this subsection, we collect several results 
from the literature
concerning the cyclicity of the backwards shift operator and the translation semigroup. In addition, we translate these results into our setting and explore their consequences.

\begin{proposition}[{\cite[Thm. 2.2.4, Rem. 2.2.6]{DOUGSHAPSHIELDS}}]\label{prop:oldresultcyclic} Consider a sequence $\alpha=(\alpha_n)_{n\geq 0}\in \ell^2(\N_0;\C)$ that decays exponentially and write
$f$ for the associated function
\begin{equation}
\begin{array}{lcl}
\label{eq:deg:def:f}
    f(z)&=&\sum\limits_{n=0}^\infty \alpha_n z^n,
\end{array}
\end{equation}
defined on its natural domain in $\C$. Then the sequence $\alpha$ is cyclic for the backwards shift operator \eqref{eq:deg:def:bck:shift} if and only if $f$ is not a rational function. In fact, if $\alpha$ is not cyclic, then $\span\{S^N\alpha:N\geq 0\}$ is finite dimensional in $\ell^2(\N_0;\C)$. \end{proposition}

\begin{lemma}\label{prop:oldresultcyclictranslation}
For any $T>0$ and any function $f\in L^1\big([0,\infty);\C\big)$ that is cyclic for the translation  group $(S_t)_{t\geq 0}$ defined in \eqref{def:semigroup}, the shifted function $s\mapsto f(s+T)$ is also cyclic for $(S_t)_{t\geq 0}$. In addition, for any $\tilde{\eta} > 0$, there exists
a function $f\in L^1_{\tilde{\eta}}\big([0,\infty);\C\big)$ that is cyclic for the translation  group $(S_t)_{t\geq 0}$. In particular, there exists a convolution kernel that satisfies both (\asref{aannamesconvoluties}{\text{H}}) and (\asref{aannamesatomicconv2}{\text{h}}).
 
\end{lemma}
\textit{Proof.} The first statement follows directly from \cite[Lem. 1]{matsui2001supercyclic}.
Turning to the existence claim,
we fix $\tilde{\eta}>0$ and let $(T_t)_{t\geq 0}$ be the translation semigroup on $L^1_{\tilde{\eta}}\big([0,\infty);\C\big)$. It follows from \cite[Thm. 1(i)]{matsui2003supercyclic} that there exists $f\in L^1_{\tilde{\eta}}\big([0,\infty);\C\big) $ that is supercyclic for $(T_t)_{t\geq 0}$, which means that $\{\lambda S(t)f:t\geq 0,\lambda\in\R\}$ is dense in $L^1_{\tilde{\eta}}\big([0,\infty);\C\big)$. Such a function is clearly also cyclic for $(T_t)_{t\geq 0}$ (with respect to the norm $\nrm{\cdot}_{\tilde{\eta}}$). We write
\begin{equation}
    \begin{array}{lclcl}
         \mathcal{D}&=& \span\{T(t)f:t\geq 0\}&=&\span\{S(t)f:t\geq 0\}.
    \end{array}
\end{equation}
Since $L^1_{\tilde{\eta}}\big([0,\infty);\C\big)$ contains all compactly supported functions, we see that $L^1_{\tilde{\eta}}\big([0,\infty);\C\big)$ is dense in $L^1\big([0,\infty);\C\big)$ with respect to the usual norm $\nrm{\cdot}_{L^1}$. Hence it is sufficient to show that $\mathcal{D}$ is dense in $L^1_{\tilde{\eta}}\big([0,\infty);\C\big)$ with respect to $\nrm{\cdot}_{L^1}$. Fix any $g\in L^1_{\tilde{\eta}}\big([0,\infty);\C\big)$ and let $\{g_n\}_{n\geq 1}$ be a sequence in $\mathcal{D}$ with
\begin{equation}
    \begin{array}{lcl}
         \lim\limits_{n\rightarrow\infty}\nrm{g_n-g}_{\tilde{\eta}}&=&0.
    \end{array}
\end{equation}
For $n\in\N$ we can compute
\begin{equation}
    \begin{array}{lcl}
         \nrm{g_n-g}_{\tilde{\eta}}&=&\int_{\R}e^{\tilde{\eta}|\xi|}|g_n(\xi)-g(\xi)|d\xi\\[0.2cm]
         &\geq &\int_{\R}|g_n(\xi)-g(\xi)|d\xi\\[0.2cm]
         &=&\nrm{g_n-g}_{L^1},
    \end{array}
\end{equation}
which immediately implies that also $g_n\rightarrow g$ in $L^1\big([0,\infty);\C\big)$, as desired. Hence $f$ is cyclic for the translation  group $(S_t)_{t\geq 0}$.
In particular, the convolution kernel
\begin{equation}
    \begin{array}{lcl}
         \mathcal{K}(\xi;t)&=&f(|\xi|) 
    \end{array}
\end{equation}
satisfies both (\asref{aannamesconvoluties}{\text{H}}) and (\asref{aannamesatomicconv2}{\text{h}}).
\qed\\

\begin{lemma}\label{lemma:cyclicvectors} Let $\{D_n\}_{n\geq 0}$ be an exponentially decaying sequence of $M\times M$ diagonal matrices. Then the following statements are equivalent.
\begin{enumerate}[label=(\roman*)]
    \item\label{enum:cyclicperp}
     There exists a non-zero sequence $y\in \ell^2(\N_0;\C^M)$ that satisfies \begin{equation}
     \begin{array}{lcl}
       \label{eq:dev:inn:prod:y:D:zero}
         \sum\limits_{n=0}^\infty y_n^\dagger D_{n+N}&=&0
     \end{array}
     \end{equation} 
     for each $N\in\Z_{\geq 0}$.

    \item\label{enum:cyclicdefinition} There exists at least one $1\leq k\leq M$ for which the sequence $(D_n^{(k,k)})_{n\geq 0}$ is not cyclic for the backwards shift operator on $\ell^2(\N_0;\C)$.
\end{enumerate}
In addition, if these statements hold, then the sequence $y$ in \ref{enum:cyclicperp} can be chosen to decay exponentially. Finally, if these conditions do not hold, then they also do not hold for
the shifted sequence $\{D_n\}_{n\geq N}$,
for any $N\in\Z_{\geq 0}$.
%
\end{lemma}
\textit{Proof.} 
As a preparation, we introduce the sequences
\begin{equation}
\begin{array}{lcl}
    \alpha^{(k);N} &=& (\alpha^{(k);N}_n)_{n \ge 0} = \big( D^{(k,k)}_{n+N} \big)_{n \ge 0}
    \in \ell^2(\N_0;\C)
\end{array}
\end{equation}
for any $N \ge 0$ and any $1 \leq k \le M$. 
In addition, we define the associated subspaces
\begin{equation}
\label{eq:deg:defn:subsp:D:k}
    \begin{array}{lcl}
        \mathcal{D}^{(k)}&=&\span\{\alpha^{(k);N}\text{ }|\text{ }N\geq 0\}
    \end{array}
\end{equation}
for $1 \le k \le M$.\\

Let us first assume that \ref{enum:cyclicperp} holds, but that \ref{enum:cyclicdefinition} fails. Then
the subspaces $\mathcal{D}^{(k)}$
are all dense in $\ell^2(\N_0;\C)$. In addition,
our diagonality assumption together with
\eqref{eq:dev:inn:prod:y:D:zero} implies that
\begin{equation}
    \begin{array}{lcl}
         \Big\langle y^{(k)},\alpha^{(k);N} \Big\rangle_{\ell^2(\N_0;\C)}&=&0 
    \end{array}
\end{equation} 
for any $N\in\Z_{\geq 0}$ and $1 \le k \le M$ and thus 
\begin{equation}
    \begin{array}{lcl}
         \big\langle y^{(k)},d\big\rangle_{\ell^2(\N_0;\C)}&=&0
    \end{array}
\end{equation} 
for any $d\in\mathcal{D}^{(k)}$ and $1 \le k \le M$.
Together these two properties yield
the contradiction $y=0$. \\

Let us now assume that \ref{enum:cyclicdefinition}
holds.
Then Proposition \ref{prop:oldresultcyclic} implies
there exists $1 \le k_0 \le M$ for which the subspace
$\mathcal{D}^{(k_0)}$ defined in \eqref{eq:deg:defn:subsp:D:k} is finite dimensional, with a basis that consists of exponentially decaying sequences.
In particular, we can pick an exponentially decaying sequence $\psi\in \ell^2(\N_0;\C)$ that satisfies
$\ip{\psi,d}_{\ell^2(\N_0;\C)}=0$ for any $d\in\mathcal{D}^{(k_0)}$.
Upon writing $y=(0,...,0,\psi,0,...,0)\in \ell^2(\N_0;\C^M)$, where $\psi$ takes the $k_0^{\mathrm{th}}$ position, we hence see that
\eqref{eq:dev:inn:prod:y:D:zero} is satisfied
by construction.\\

The final statement follows 
from the characterization in Proposition \ref{prop:oldresultcyclic}, which implies that (non)-cyclicity is preserved under translation. Indeed,
if the function $f$ defined in \eqref{eq:deg:def:f}
is not a rational function, then the function 
\begin{equation}
    \begin{array}{lcl}
         f_N(z)&=&z^{-N}\big[f(z)-\sum\limits_{n=0}^{N-1} \alpha_n z^n\big]
    \end{array}
\end{equation}
associated to the shifted sequence $S^N \alpha$
is also not rational.\qed\\

\textit{Proof of Lemmas \ref{lem:cyclicshifts1} and \ref{lem:cyclicshifts2}.} Both results follow directly from 
Proposition \ref{prop:oldresultcyclic} 
and Lemma \ref{lemma:cyclicvectors}.\qed\\

\subsection{Nondegeneracy of the Hale inner product}\label{subsec:hale}

In this subsection we show how the nondegeneracy of the Hale inner product can be derived from the conditions formulated in \S \ref{subsec:assump}.
In particular, we establish Propositions \ref{prop:4.16:infrangereverse} and \ref{prop:4.16volledig}. \\

As a convenience, we 
first connect the right-nondegeneracy properties
for the system (\ref{ditishetprobleem}) to the left-nondegeneracy properties for the adjoint system (\ref{ditishetadjointprobleem}). 
This will allows us to focus solely 
on the left-nondegeneracy of the Hale inner product with respect to functions in $B^*(\tau)$.

\begin{lemma}\label{leftnonvsrightnon}
Assume that (\asref{aannamesconstanten}{\text{H}}), (\asref{aannamesconvoluties}{\text{H}}) and (\asref{aannameslimit}{\text{H}}) are satisfied. Fix $\tau\in\R$ and $E\subset X$ with $0\in E$. Then the Hale inner product for the system (\ref{ditishetprobleem}) at $\tau$ is right-nondegenerate for functions in $E$ if and only if the Hale inner product for the adjoint system (\ref{ditishetadjointprobleem}) at $\tau$ is left-nondegenerate for functions in $E$. 
\end{lemma}
\textit{Proof.} 
For any $\phi\in X$, $\psi\in Y$ and $\tau\in\R$,
the Hale inner product for the adjoint system (\ref{ditishetadjointprobleem}) is given by
\begin{equation}
    \begin{array}{lcl}
         \ip{\phi,\psi}_\tau^\mathrm{adj}&=&\phi(0)^\dagger \psi(0)+\sum\limits_{j=-\infty}^\infty \int\limits_0^{-r_j} \phi(s+r_j)^\dagger A_j(\tau+s-r_j)^\dagger\psi(s)ds\\[0.2cm]
         &&\qquad +\int\limits_{\R}\int\limits_0^r \phi(s-r)^\dagger \mathcal{K}(s-r;\tau+s-r)^\dagger\psi(s)dsdr . \\[0.2cm]
\end{array}
\end{equation}
A short computation shows that
\begin{equation}
    \begin{array}{lcl}
         \overline{\ip{\phi,\psi}_\tau^\mathrm{adj}}&=&
         \psi(0)^\dagger \phi(0)-\sum\limits_{j=-\infty}^\infty \int\limits_0^{r_j} \psi(s-r_j)^\dagger A_j(\tau+s-r_j)\phi(s)ds\\[0.2cm]
         &&\qquad -\int\limits_{\R}\int\limits_0^r \psi(s-r)^\dagger \mathcal{K}(r;\tau+s-r)\phi(s)dsdr\\[0.2cm]
         &=&\ip{\psi,\phi}_\tau,
    \end{array}
\end{equation}
which directly implies the desired result. \qed\\

We proceed by discussing the cyclicity criteria
in
introduced in {\S}\ref{sec:deg:unb:cmp}
and {\S}\ref{sec:deg:bnd:unbnd:kernel}. The following
preparatory result will help us to link
the discussion in \eqref{sec:ver:cyc}
to the degeneracy properties of the Hale inner product.


\begin{lemma} 
\label{lem:deg:tech:heart:degener}
Assume that (\asref{aannamesconstanten}{\text{H}}), (\asref{aannamesconvoluties}{\text{H}}), (\asref{aannameslimit}{\text{H}}) 
and (\asref{aannamesatomicbasis}{\text{h}})
are satisfied and fix $\tau \in \R$.   Pick any $\psi\in Y$ that does not vanish on $D_Y^+$ and satisfies $\ip{\psi,\phi}_\tau=0$ for every $\phi\in X$. Writing 
\begin{equation}
\begin{array}{lcl}
\sigma=\inf\{s\in D_Y^+\text{ }|\text{ }\psi(s)\neq 0\},
\end{array}
\end{equation}
there exist $\epsilon>0$ and $N_0\in\Z_{\geq  K_{\mathrm{const}}}$
so that the identity
\begin{equation}\label{eq:perpendicularlemma}
    \begin{array}{lcl}
         \sum\limits_{j=0}^{\infty}\psi(s+j)^\dagger \tilde{A}_{-j-N}+\int\limits_{\sigma-s}^{\infty}\psi(s+r)^\dagger \tilde{\mathcal{K}}(-r-N)dr&=&0
    \end{array}
\end{equation}
holds for almost every $s\in (\sigma,\sigma+\epsilon)$ and every 
integer $N \ge N_0$. In addition, if $\tilde{A}_{j}=0$ for each $j\leq -N_0$, then we in fact have
\begin{equation}\label{eq:perpendicularlemma2}
    \begin{array}{lcl}
         \int\limits_{0}^{\infty}\psi(\sigma+r)^\dagger \tilde{\mathcal{K}}(-r- \theta)dr&=&0
    \end{array}
\end{equation}
for all (reals) $\theta \ge N_0 + \epsilon$.
\end{lemma}

\textit{Proof.} 
We first pick an arbitrary $s < 0$ with $s \notin \mathbb{Z}$.
Using
a sequence of functions supported on small intervals that shrink to the singleton $\{s\}$,
we can use \eqref{eq:defHale} to conclude that
\begin{equation}\label{eq:4.38inf}
    \begin{array}{lcl}
\sum\limits_{ j < s }\psi(s-j)^\dagger A_{j}(\tau+s-j)+\int\limits_{-\infty}^{s}\psi(s-r)^\dagger \mathcal{K}(r;\tau+s-r)dr&=&0 .
    \end{array}
\end{equation}
Imposing the further restriction $s \le - K_{\mathrm{const}}$, this can be rephrased as
\begin{equation}
\label{eq:deg:zero:id:for:sum:kernel:psi}
    \begin{array}{lcl}
         \sum\limits_{j < s}\psi(s-j)^\dagger \tilde{A}_{j}+\int\limits_{-\infty}^{s}\psi(s-r)^\dagger \tilde{\mathcal{K}}(r)dr&=&0.
    \end{array}
\end{equation}

We now choose $\epsilon > 0$ to be so small that
$(\sigma, \sigma + \epsilon)$ contains no integers.
Then for any sufficiently large integer $N \gg 1$,
we can combine \eqref{eq:deg:zero:id:for:sum:kernel:psi} 
together with the definition of $\sigma$ to conclude that
\begin{equation}
    \begin{array}{lcl}
         \sum\limits_{j< s - \sigma}\psi(s-j)^\dagger \tilde{A}_{j}+\int\limits_{-\infty}^{s-\sigma}\psi(s-r)^\dagger \tilde{\mathcal{K}}(r)dr&=&0
    \end{array}
\end{equation}
for all $s \in (\sigma - N, \sigma + \epsilon - N)$.
This yields \eqref{eq:perpendicularlemma} upon
introducing new variables 
\begin{equation}
\begin{array}{lcl}
(s',j', r') &=& (s + N, -j-N , -r - N)
\end{array}
\end{equation}
and dropping the primes, noting that $\lceil \sigma - s' \rceil = 0$. The final statement follows
from the fact that we no longer need to rule out integer values of $s'$ above, together with
the replacement $r \mapsto r + \sigma - s$.
\qed\\

\begin{lemma}\label{prop:4.16:infrange} Assume that (\asref{aannamesconstanten}{\text{H}}), (\asref{aannamesconvoluties}{\text{H}}) and (\asref{aannameslimit}{\text{H}}) are satisfied
and fix $\tau \in \R$. Assume moreover that the cyclicity conditions
(\asref{aannamesatomic2:prep}{\text{h}})-(\asref{aannamesatomic2}{\text{h}}) are satisfied. Then the Hale inner product 
at $\tau$ is left-nondegenerate for exponentially decaying functions. 
\end{lemma}
\textit{Proof.} 
Assume that $\psi \in Y$
decays exponentially and has
$\ip{\psi,\phi}_\tau=0$ for every $\phi\in X$. Exploiting symmetry, we
assume further that $\psi$ does not vanish on $D_Y^+$
and set out to find a contradiction.
Recalling the setting
of Lemma \ref{lem:deg:tech:heart:degener}
and remembering that $\tilde{\mathcal{K}}=0$ on $(-\infty,-K_{\mathrm{const}}]$, we obtain from (\ref{eq:perpendicularlemma}) that the identity
\begin{equation}\label{eq:bewijs:shift}
    \begin{array}{lcl}
         \sum\limits_{j=0}^{\infty}\psi(s+j)^\dagger \tilde{A}_{-j-N}&=&0
    \end{array}
\end{equation}
holds for almost every $s\in (\sigma,\sigma+\epsilon)$ and every $N\geq N_0 \ge K_{\mathrm{const}}$. \\

By (\asref{aannamesatomic2}{\text{h}}) and the invariance of cyclicity under translations, the sequences
$(\tilde{A}_{-j}^{(k,k)})_{j\geq N}$
are cyclic for each $1 \le k \le M$. In particular,
Lemma \ref{lemma:cyclicvectors} implies that
that the sequence $\psi(s + \mathbb{N}_0) \in \ell^2(\mathbb{N}_0 ;\C^M)$
and hence also the first coordinate $\psi(s)$
must vanish for all $s\in (\sigma,\sigma+\epsilon)$.
This contradicts the definition of $\sigma$.
\qed\\

\textit{Proof of Proposition \ref{prop:4.16:infrangereverse}.} Assume without loss of generality that the sequence $( \tilde{A}_{-j}^{(k,k)})_{j\geq K_{\mathrm{const}}}$ is not cyclic for the backwards shift operator. 
Lemma \ref{lemma:cyclicvectors} then allows us to pick an exponentially decaying nonzero sequence
\begin{equation}
\begin{array}{lcl}
    y&=&(y_n)_{n\geq 0}\in \ell^2(\N_0;\C^M)
\end{array}
\end{equation}
for which the identity
\begin{equation}
\begin{array}{lcl}
    \sum\limits_{j=0}^\infty y_j^\dagger  \tilde{A}_{-j -N}&=&0
\end{array}
\end{equation}
holds for all integers $N \ge K_{\mathrm{const}}$. \\

We now define a continuous, bounded function $\psi:D_Y\rightarrow \C^M$ by writing
\begin{equation}
\begin{array}{lcl}
     \psi(s) &=& 0, \qquad \qquad  s\in (-\infty, K_{\mathrm{const}} ),
\end{array}
\end{equation}
together with
\begin{equation}
\begin{array}{lclcl}
    \psi(j) &=& 0, \qquad \psi(j + \frac{1}{2} ) &= &
    y_{j- K_{\mathrm{const}}},
    \qquad \qquad  j \in \Z_{\geq  K_{\mathrm{const}}}
\end{array}
\end{equation}
and performing a linear interpolation between these prescribed values. This construction implies that
\begin{equation}\label{eq:chardefpsi}
    \begin{array}{lclcl}
         \sum\limits_{j=N}^\infty \psi\big(\tilde{s}+j+ K_{\mathrm{const}}-N\big)^\dagger \tilde{A}_{-j}&=&\tilde{s} \sum\limits_{j=0}^\infty  y_j^\dagger  \tilde{A}_{-N-j}
         &=&0, 
    \end{array}
\end{equation}
for any integer $N \ge  K_{\mathrm{const}}$ and any $\tilde{s}\in[0,1)$. \\

Let us now consider an arbitrary $s \le 0$
and make the decomposition
\begin{equation}
\begin{array}{lcl}
    s &=& \tilde{s} + K_{\mathrm{const}} - N
\end{array}
\end{equation}
for some integer $N \ge K_{\mathrm{const}}$ and $\tilde{s} \in [0,1)$. Applying \eqref{eq:chardefpsi},
we now compute
\begin{equation}
\label{eq:deg:id:for:s:neg}
\begin{array}{lclcl}
    \sum\limits_{j \ge K_{\mathrm{const}} - s} \psi(s + j) \tilde{A}_{-j}
    &=& \sum\limits_{j \ge N - \tilde{s}} \psi(\tilde{s} + j + K_{\mathrm{const}} - N ) \tilde{A}_{-j}
    &=& 0
\end{array}
\end{equation}
since the final sum in fact ranges over $j \ge N$.\\

For any $\phi\in X$, we note that
\eqref{eq:defHale}
reduces to
\begin{equation}
    \begin{array}{lcl}
         \ip{\psi,\phi}_{\tau}&=&
         - \sum\limits_{j=-\infty}^\infty \int\limits_0^{j} \psi(s-j)^\dagger A_j(\tau +s-j)\phi(s)ds
         \\[0.2cm]
         && \qquad 
           -  \int\limits_{\R}\int\limits_0^r \psi(s-r)^\dagger \mathcal{K}(r;t+s-r)\phi(s)dsdr
         \\[0.2cm]
  \end{array}
\end{equation}
since $\psi(0) = 0$. Exploiting (\asref{aannamesatomicbasis}{\text{h}}),
this can be further simplified and recast as
\begin{equation}
    \begin{array}{lcl}
         \ip{\psi,\phi}_{\tau}
         & = & 
         -\sum\limits_{j=-\infty}^{-K_{\mathrm{const}}} \int\limits_0^{j} \psi(s-j)^\dagger \tilde{A}_j\phi(s)ds
         \\[0.2cm]
         & = & -\sum\limits_{j= K_{\mathrm{const}}}^\infty \int\limits_0^{K_{\mathrm{const}} - j} \psi(s+j)^\dagger  \tilde{A}_{-j}\phi(s)ds. \\[0.2cm]
    \end{array}
\end{equation}
The dominated convergence theorem allows us to interchange the infinite sum and the integral, which yields
\begin{equation}
    \begin{array}{lclcl}
         \ip{\psi,\phi}_{\tau}&=&-\int\limits_{0}^{-\infty} \sum\limits_{j \ge K_{\mathrm{const}} -s} \psi(s + j )^\dagger \tilde{A}_{-j}\phi(s)ds
         & = & 0
    \end{array}
\end{equation}
on account of \eqref{eq:deg:id:for:s:neg}.
\qed\\


\begin{lemma}\label{prop:4.16:convolutie} Assume that (\asref{aannamesconstanten}{\text{H}}), (\asref{aannamesconvoluties}{\text{H}}) and (\asref{aannameslimit}{\text{H}}) are satisfied and fix $\tau \in \R$.  Assume moreover that the cyclicity conditions (\asref{aannamesatomicconv2:prep}{\text{h}})-(\asref{aannamesatomicconv2}{\text{h}}) are satisfied. Then the Hale inner product 
at $\tau$
is left-nondegenerate  for functions in $Y$. 
\end{lemma}
\textit{Proof.}
Assume that $\psi \in Y$
has
$\ip{\psi,\phi}_\tau=0$ for every $\phi\in X$. Exploiting symmetry, we
assume further that $\psi$ does not vanish on $D_Y^+$
and set out to find a contradiction. \\

We pick $1\leq k\leq M$ for which $\psi^{(k)}$ does not vanish on $D_Y^+$. Recalling the setting
of Lemma \ref{lem:deg:tech:heart:degener}
and remembering that $\tilde{A}_j=0$ for each $|j|\geq K_{\mathrm{const}}$, we obtain from (\ref{eq:perpendicularlemma2}) that the identity
\begin{equation}\label{eq:bewijs:conv}
    \begin{array}{lcl}
         \int\limits_{0}^{\infty}\psi(\sigma+r)^\dagger \tilde{\mathcal{K}}(-r - \theta)dr&=&0
    \end{array}
\end{equation}
holds for every $\theta \ge N + \epsilon$. 
We introduce the subspace
\begin{equation}
    \begin{array}{lcl}
         \mathcal{D}&=&\span\big\{t\mapsto \tilde{\mathcal{K}}^{(k,k)}(-t-r)\text{ }|\text{ }r\geq N+\epsilon \big\},
    \end{array}
\end{equation}
which is dense in $L^1\big([0,\infty);\C\big)$ by
(\asref{aannamesatomicconv2}{\text{h}})
and
Lemma \ref{prop:oldresultcyclictranslation}. We therefore have
\begin{equation}\label{eq:perpincyclicconv}
    \begin{array}{lcl}
         \int\limits_{0}^{\infty}\psi^{(k)}(\sigma+r)^*f(r)dr&=&0
    \end{array}
\end{equation}
for every $f\in\mathcal{D}$. \\

We fix any $f\in L^1\big([0,\infty);\C\big)$ and let $\{f_n\}_{n\geq 1}$ be a sequence in $\mathcal{D}^{(k)}$ with $f_n\rightarrow f$. Using (\ref{eq:perpincyclicconv}) we can estimate
\begin{equation}
    \begin{array}{lcl}
         \Big|\int\limits_{0}^{\infty}\psi^{(k)}(\sigma+r)^*f(r)dr\Big|&=&\Big|\int\limits_{0}^{\infty}\psi^{(k)}(\sigma+r)^*\big(f(r)-f_n(r) \big) dr\Big|\\[0.2cm]
         &\leq &\nrm{\psi}_\infty\int\limits_{0}^{\infty}|f(r)-f_n(r)|dr,
         
    \end{array}
\end{equation}
which converges to $0$ as $n\rightarrow \infty$. Hence (\ref{eq:perpincyclicconv}) holds for any $f\in L^1\big([0,\infty);\C\big)$. In particular, we pick $s\in (\sigma,\sigma+\epsilon)$ for which $\psi^{(k)}(s)\neq 0$ and we let $f\in L^1\big([0,\infty);\C\big)$ be a sufficiently small peak function, centered around $s-\sigma$. This immediately yields
\begin{equation}
    \begin{array}{lcl}
         \int\limits_{0}^{\infty}\psi^{(k)}(\sigma+r)^*f(r)dr&\neq&0,
    \end{array}
\end{equation}
which contradicts (\ref{eq:perpincyclicconv}).\qed\\

\begin{lemma}\label{prop:4.16:positive} Assume that (\asref{aannamesconstanten}{\text{H}}), (\asref{aannamesconvoluties}{\text{H}}) and (\asref{aannameslimit}{\text{H}}) are satisfied and fix $\tau \in \R$.  Assume moreover that the positivity condition (\asref{aannamesatomicpositivestronger}{\text{h}}) is satisfied. Then the Hale inner product 
at $\tau$ is left-nondegenerate 
for nonnegative functions. 
\end{lemma}
\textit{Proof.}
Assume that $\psi \in Y$
is nonnegative and has
$\ip{\psi,\phi}_\tau=0$ for every $\phi\in X$. Exploiting symmetry, we
assume further that $\psi$ does not vanish on $D_Y^+$
and set out to find a contradiction. 
Recalling the setting
of Lemma \ref{lem:deg:tech:heart:degener},
we obtain from (\ref{eq:perpendicularlemma}) that the identity
\begin{equation}\label{eq:positivedefcontracition}
    \begin{array}{lcl}
         \sum\limits_{j=0}^{\infty}\psi(s+j)^\dagger \tilde{A}_{-j-N}+\int\limits_{\sigma - s}^{\infty}\psi(s+r)^\dagger \tilde{\mathcal{K}}(-r -N )dr&=&0
    \end{array}
\end{equation}
holds for almost every $s\in (\sigma,\sigma+\epsilon)$ and every 
$N\geq N_0 \ge K_{\mathrm{const}}$. 
In addition, the definition of $\sigma$ allows us to conclude $\psi(s) > 0$ for $s \in (\sigma, \sigma + \epsilon)$.\\

Since the matrices
\eqref{eq:deg:coeffs:tilde:a:k} are all positive semidefinite, we have 
\begin{equation}\label{eq:positivedef:sh}
\begin{array}{lclcl}
\big(\psi(s+j)^\dagger \tilde{A}_{-j-N}\big)^{(k)}&\geq &0,
\qquad \qquad s \in (\sigma, \sigma + \epsilon)
\end{array}
\end{equation}
for all $j \ge 0$, $1\leq k\leq M$ and $N \ge N_0$, together with
\begin{equation}\label{eq:positivedef:conv}
\begin{array}{lclcl}
 \big(\psi(s+r)^\dagger \tilde{\mathcal{K}}(-r-N)\big)^{(k)}&\geq &0,
 \qquad \qquad 
 s \in (\sigma, \sigma + \epsilon)
\end{array}
\end{equation}
for all $r \ge \sigma - s$, $1\leq k\leq M$ and all $N \ge N_0$.  
On the other hand, fixing $j = 0$ and $r = 0$, 
item (a) and (b) in (\asref{aannamesatomicpositivestronger}{\text{h}})
allow us to find $N \ge N_0$ 
for which one or both of the inequalities
\eqref{eq:positivedef:sh}-\eqref{eq:positivedef:conv}
are strict. This immediately contradicts
(\ref{eq:positivedefcontracition}).\qed\\

\begin{lemma}
\label{prop:4.16} Assume that (\asref{aannamesconstanten}{\text{H}}), (\asref{aannamesconvoluties}{\text{H}}) and (\asref{aannameslimit}{\text{H}}) are satisfied. Assume moreover that the atomic condition (\asref{aannamesatomicbdd}{\text{h}}) is satisfied at some point $\tau\in\R$. Then the Hale inner product 
at $\tau$
is left-nondegenerate 
for functions in $Y$. 
\end{lemma}
\textit{Proof.} 
The proof is identical to that of \cite[Prop. 4.16]{MPVL} and, as such, will be omitted.\qed\\

\textit{Proof of Proposition \ref{prop:4.16volledig}.} The statements \ref{enum:halebdd}-\ref{enum:halepos}
follow from
Lemmas \ref{leftnonvsrightnon},  \ref{prop:4.16:infrange}, \ref{prop:4.16:convolutie},
\ref{prop:4.16:positive} and
\ref{prop:4.16}.
%
%
The final statement follows from the representation (\ref{eq:4.12}),
applying Proposition \ref{prop:fredholmproperties} for \ref{enum:halecycsh}-\ref{enum:halecycconv} or 
using the nonnegative $B^*(\tau)$-basis
for \ref{enum:halepos}. \qed\\

\begin{remark}
The conclusion in Lemma \ref{prop:4.16:infrange} that the Hale inner product is nondegenerate for exponentially decaying functions cannot easily be generalized to bounded functions. Indeed, the key argument
is that the sequence $\psi(s+\N_0)^{(k)}$ is perpendicular to a dense subspace of $\ell^2(\N_0;\C)$. This sequence is in $\ell^2$ itself on account of the exponential decay of $\psi$ and must therefore vanish. However, it \textit{is} possible for non-trivial $\ell^\infty$ sequences to be perpendicular to a dense subspace of $\ell^2(\N_0;\C)$;
see the discussion at \cite{StackExPerp}.
In a similar fashion, we do not expect Lemma \ref{prop:4.16:positive} to be easily generalizable. 
\end{remark}

\subsection{The nontriviality condition \texorpdfstring{(\asref{aannames0ophalflijn}{\text{H}})}{HKer}}\label{subsec:HAker}
In this final subsection we show how the nontriviality condition (\asref{aannames0ophalflijn}{\text{H}}) can be verified. 
\begin{lemma}\label{lemma:assumptions1} Assume that (\asref{aannamesconstanten}{\text{H}}), (\asref{aannamesconvoluties}{\text{H}}) and (\asref{aannameslimit}{\text{H}}) are satisfied. Suppose that the atomic condition (\asref{aannamesatomicbdd}{\text{h}}) holds for the system (\ref{ditishetprobleem}) at each $\tau\in\R$. Then the non-triviality condition (\asref{aannames0ophalflijn}{\text{H}}) 
is also satisfied. 
\end{lemma}
\textit{Proof.} By symmetry and the fact
the adjoint system (\ref{ditishetadjointprobleem}) also satisfies (\asref{aannamesatomicbdd}{\text{h}}),
it suffices to show that any nonzero $d\in\mathcal{B}$ cannot vanish on $(-\infty,0]$. Arguing by contradiction, we assume that $d=0$ identically on $(-\infty,0]$. Defining $\sigma=\inf\{s\in\R:d(s)\neq 0\}$,
we have $0\leq\sigma<\infty$ by construction. \\

Recalling the constant $\delta>0$ from (\asref{aannamesatomicbdd}{\text{h}}), we pick $0<\epsilon<\delta$ sufficiently small to have $d(\sigma+\epsilon)\neq 0$ and $r_j+\epsilon<\rmax$ for any $j\in\Z$ with $r_j\neq \rmax$. Evaluating (\ref{ditishetprobleem}) at $t=\sigma+\epsilon-\rmax$ now yields
\begin{equation}
\begin{array}{lcl}
0&=&-\dot{d}(t)+\sum\limits_{j=1}^\infty A_{j}(t)d(t+r_j)+\int\limits_\R\mathcal{K}(\xi;t)d(t+\xi)d\xi\\[0.2cm]
         &=&A_{j_{\mathrm{max}}}(\sigma+\epsilon-\rmax)d(\sigma+\epsilon).
    \end{array}
\end{equation}
Since the matrix $A_{j_{\mathrm{max}}}(\sigma+\epsilon-\rmax)$ is nonsingular, we obtain the desired contradiction
$d(\sigma+\epsilon)=0$.
\qed\\

\begin{lemma}\label{lemma:assumptions2} Assume that (\asref{aannamesconstanten}{\text{H}}), (\asref{aannamesconvoluties}{\text{H}}) and (\asref{aannameslimit}{\text{H}}) are satisfied,
together with the
cyclicity conditions (\asref{aannamesatomic2:prep}{\text{h}})-(\asref{aannamesatomic2}{\text{h}}).
Then the non-triviality condition (\asref{aannames0ophalflijn}{\text{H}}) also holds.
\end{lemma}
\textit{Proof.} By symmetry and the fact
the adjoint system (\ref{ditishetadjointprobleem}) also satisfies (\asref{aannamesatomic2:prep}{\text{h}})-(\asref{aannamesatomic2}{\text{h}}),
it suffices to show that any nonzero $d\in\mathcal{B}$ cannot vanish on $(-\infty,0]$. Writing $\sigma=\inf\{s\in\R:d(s)\neq 0\}$,  
we have $0\leq \sigma<\infty$  by construction. Recalling the constant $ K_{\mathrm{const}}\in\Z_{\geq 0}$ from (\asref{aannamesatomicbasis}{\text{h}}),
we use (\ref{ditishetprobleem})
to conclude that
\begin{equation}\label{eq:shiftimpker}
\begin{array}{lcl}
0&=&-\dot{d}(s)+\sum\limits_{j\in\Z}A_{j}(s)d(s+j)+\int\limits_{\R}\mathcal{K}(\xi;s)d(s+\xi)d\xi\\[0.2cm]
&=&\sum\limits_{j\geq \sigma-s}\tilde{A}_{j}d(s+j)
\end{array}
\end{equation}
for any $s\in(-\infty,- K_{\mathrm{const}}]$. \\

We now pick an integer $N_0$ and a constant $\epsilon > 0$ in such a way that  $N_0>\sigma + K_{\mathrm{const}}$ and  $d(\sigma+\epsilon)\neq 0$
both hold. Then for any integer $N \ge N_0$
and any $s\in(\sigma,\sigma+\epsilon)$,
we can use \eqref{eq:shiftimpker} to conclude
\begin{equation}
    \begin{array}{lcl}
         \sum\limits_{j=0}^\infty d(s+j)^\dagger  \tilde{A}_{j+N}&=&0 ,
    \end{array}
\end{equation}
which closely resembles (\ref{eq:bewijs:shift}).
We can hence follow the proof of Lemma \ref{prop:4.16:infrange} to obtain the contradiction $d=0$. The other statements follow in a similar fashion, noting that the adjoint system (\ref{ditishetadjointprobleem}) also satisfies (\asref{aannamesatomic2:prep}{\text{h}})-(\asref{aannamesatomic2}{\text{h}}).
 \qed\\

\begin{lemma}\label{lemma:assumptions3} Assume that (\asref{aannamesconstanten}{\text{H}}), (\asref{aannamesconvoluties}{\text{H}}) and (\asref{aannameslimit}{\text{H}}) are satisfied. Suppose that the cyclicity conditions (\asref{aannamesatomicconv2:prep}{\text{h}})-(\asref{aannamesatomicconv2}{\text{h}}) are satisfied for the system (\ref{ditishetprobleem}). Then the non-triviality condition (\asref{aannames0ophalflijn}{\text{H}}) is satisfied for the system (\ref{ditishetprobleem}).\end{lemma}
\textit{Proof.} By symmetry and the fact
the adjoint system (\ref{ditishetadjointprobleem}) also satisfies (\asref{aannamesatomicconv2:prep}{\text{h}})-(\asref{aannamesatomicconv2}{\text{h}}),
it suffices to show that any nonzero $d\in\mathcal{B}$ cannot vanish on $(-\infty,0]$. We can follow the proof of Lemmas 
\ref{prop:4.16:convolutie} and
\ref{lemma:assumptions2} to arrive at a contradiction. The other statements follow in a similar fashion, noting that the adjoint system (\ref{ditishetadjointprobleem}) also satisfies (\asref{aannamesatomicconv2}{\text{h}}).\qed\\

\begin{lemma}\label{lemma:assumptions4} Assume that (\asref{aannamesconstanten}{\text{H}}), (\asref{aannamesconvoluties}{\text{H}}) and (\asref{aannameslimit}{\text{H}})
and (\asref{aannamesatomicpositivestronger}{\text{h}})
are satisfied. Suppose furthermore that the positivity condition (\asref{aannamesatomicpositivestronger2}{\text{h}}) holds for both the system (\ref{ditishetprobleem}) and the adjoint system (\ref{ditishetadjointprobleem}). Then the non-triviality condition (\asref{aannames0ophalflijn}{\text{H}}) is also satisfied.
\end{lemma}
\textit{Proof.} By symmetry,
it suffices to show that any nonzero, nonnegative $d\in\mathcal{B}$ cannot vanish on $(-\infty,0]$. Write $\sigma=\inf\{s\in\R:d(s)\neq 0\}$ and recall the constant $ K_{\mathrm{const}}\in\Z_{\geq 0}$ from (\asref{aannamesatomicbasis}{\text{h}}). 
Using (\ref{ditishetprobleem})
we see that
\begin{equation}
\begin{array}{lcl}
0&=&-\dot{d}(s)+\sum\limits_{j\in\Z}A_{j}(s)d(s+j)+\int\limits_\R\mathcal{K}(\xi;s)d(s+\xi)d\xi\\[0.2cm]
&=&\sum\limits_{j\geq \sigma-s}\tilde{A}_{j}d(s+j)+\int\limits_{\sigma-s}^\infty \tilde{\mathcal{K}}(\xi)d(s+\xi)d\xi
\end{array}
\end{equation}
for any $s\in(-\infty,- K_{\mathrm{const}}]$. \\

We now pick an integer $N_0$ and a constant $\epsilon > 0$ in such a way that  
$d(\sigma+\delta)\neq 0$ for each $0<\delta<\epsilon$ and
$N_0>\sigma + K_{\mathrm{const}} + \epsilon$  
both hold. If (a) holds in (\asref{aannamesatomicpositivestronger}{\text{h}}),
we pick $N \ge N_0$ in such a way that $\tilde{A}_{N}$ is positive definite.
Picking $s=\sigma+\epsilon-N \in(-\infty,- K_{\mathrm{const}}]$, we arrive at the contradiction
\begin{equation}
\begin{array}{lclcl}
0
&\geq &\big(\tilde{A}_{N} d(\sigma+\epsilon)\big)^{(k)}
&>&0\end{array}
\end{equation}
for some $1\leq k\leq M$. On the other hand, if (b) holds in (\asref{aannamesatomicpositivestronger}{\text{h}}),
we pick $\theta \ge N_0$ in such a way that
$\tilde{\mathcal{K}}_{\theta+ \delta}$ is positive definite whenever $|\delta|\leq \frac{\epsilon}{4}$.
Picking $s=\sigma+\frac{\epsilon}{2}-\theta \in(-\infty,- K_{\mathrm{const}}]$, we obtain

\begin{equation}
    \begin{array}{lclcl}
0
&\geq &C \inf\limits_{t\in[\frac{\epsilon}{4},\frac{3\epsilon}{4}]}\{d(\sigma+t)^{(k)}\}
&>&0\end{array}
\end{equation}
for some constant $C>0$ and some $1\leq k\leq M$, a contradiction.
The other statements follow similarly.
\qed\\

\textit{Proof of Proposition \ref{assumptionimplications}.} This follows directly from Lemmas \ref{lemma:assumptions1}-\ref{lemma:assumptions4}.\qed\\

\bibliographystyle{plain}
\bibliography{ref}

\end{document}